\RequirePackage{fix-cm}
\documentclass{article}

\usepackage{arxiv}

\usepackage{placeins}
\usepackage[italicComments=true,indLines=true,noEnd=false]{algpseudocodex}
\usepackage{float}
\usepackage[english]{babel}
\usepackage[utf8]{inputenc} 
\usepackage[T1]{fontenc}    
\usepackage{hyperref}       
\hypersetup{colorlinks=true, citecolor=blue}
\usepackage{url}            
\usepackage{booktabs}       
\usepackage{stackrel,amssymb}
\usepackage{subfigure}
\usepackage{amsfonts,amsmath,amstext,amsbsy,amsthm}
\usepackage{nicefrac}       
\usepackage[nopatch=eqnum]{microtype}      
\usepackage{graphicx}
\usepackage[numbers,sort&compress]{natbib}
\usepackage{doi}

\newtheorem{lemma}{Lemma}

\usepackage{xcolor}
\definecolor{abstract_background}{RGB}{235,235,235}%
\usepackage{caption}
\usepackage[title]{appendix}

\usepackage[bottom]{footmisc}

\title{Slow Invariant Manifolds of Singularly Perturbed Systems via Physics-Informed Machine Learning}

\author{
\textbf{Dimitrios G. Patsatzis\textcolor{teal}{$^{1}$}, Gianluca Fabiani\textcolor{teal}{$^{1}$}, Lucia Russo\textcolor{teal}{$^{2}$}, Constantinos Siettos\textcolor{teal}{$^{3,}$}\thanks{Corresponding author, email: \texttt{constantinos.siettos@unina.it}}}
{}\\
\textcolor{teal}{$^{(1)}$} Modelling Engineering Risk and Complexity, \emph{Scuola Superiore Meridionale}, Naples 80138, Italy \hspace{1cm}\\
\textcolor{teal}{$^{(2)}$} Institute of Science and Technology for Energy and Sustainable Mobility, \emph{Consiglio Nazionale}\\ \hspace{0.38cm}\emph{ delle Ricerche}, Naples 80125, Italy\\
\textcolor{teal}{$^{(3)}$} Dipartimento di Matematica e Applicazioni ‘‘Renato Caccioppoli", \emph{Universit\`a degli Studi di Napoli}\\ \hspace{0.38cm}\emph{Federico II}, Naples 80126, Italy\\
}

\begin{document}
\maketitle
\begin{abstract}
\colorbox{abstract_background}{\begin{minipage}{1\linewidth}
We present a physics-informed machine-learning (PIML) approach for the approximation of slow invariant manifolds (SIMs) of singularly perturbed systems, providing  functionals in an explicit form that facilitate the construction and numerical integration of reduced order models (ROMs). The proposed scheme solves the partial differential equation corresponding to the invariance equation (IE) within the Geometric Singular Perturbation Theory (GSPT) framework.~For the solution of the IE, we used two neural network structures, namely feedforward neural networks (FNNs), and random projection neural networks (RPNNs), with symbolic differentiation for the computation of the gradients required for the learning process.~The efficiency of our PIML method is assessed via three benchmark problems, namely the Michaelis-Menten, the target mediated drug disposition reaction mechanism, and the 3D Sel'kov model. We show that the proposed PIML scheme provides approximations, of equivalent or even higher accuracy, than those provided by other traditional GSPT-based methods, and importantly, for any practical purposes, it is not affected by the magnitude of the perturbation parameter. This is of particular importance, as there are many systems for which the gap between the fast and slow timescales is not that big, but still 
ROMs can still be constructed.~A comparison of the computational costs between symbolic, automatic and numerical approximation of the required derivatives in the learning process is also provided.
\end{minipage}
}
\end{abstract}


\section{Introduction} 
The construction of reduced-order models (ROMs) for the multiscale mathematical modelling and numerical analysis of stiff ODEs, DAEs, PDEs and complex systems is an open and challenging problem.~A fundamental hypothesis is that the effective/long-term/slow/emergent dynamics evolve on low-dimensional invariant topological spaces (slow invariant manifolds (SIMs)) that can be parametrized by a few variables \citep{valorani2001explicit,kevrekidis2003equation,kuehn2015multiple}.

Over the years, several analytical and numerical analysis methods have been developed for the approximation of SIMs, in the context of Singular Perturbation Theory (SPT) \citep{kuehn2015multiple,verhulst2005methods} and Geometric Singular Perturbation Theory (GSPT) \citep{tikhonov1952systems,fenichel1979geometric,jones1995geometric,wechselberger2020geometric}. For singularly perturbed dynamical systems characterized by an explicit timescale splitting, expressed by a perturbation parameter, $\epsilon$, considered to be ``sufficiently small'', the main GSPT-based approaches include analytical methods based on the invariance equation (IE) \cite{verhulst2005methods,kuehn2015multiple}, the Method of Invariant Manifolds (MIM) \citep{gorban2003method} and the Rousel-Fraser method \citep{roussel1991geometry,fraser1988steady}. For fast-slow dynamical systems without an explicit timescale splitting, the key approaches include the celebrated Computational Singular Perturbation (CSP) \citep{lam1989understanding,goussis1992study,valorani2005higher}, the Intrinsic Low-Dimensional Manifold (ILDM) \citep{maas1992simplifying},  the Zero-Derivative Principle (ZDP)  \citep{gear2005projecting,zagaris2009analysis}, the Flow Curvature (FC) \citep{ginoux2008slow}, and the Tangential Stretching Rate (TSR) \citep{valorani2015dynamical} methods. These methods are specifically designed for high-dimensional systems, in which iterative numerical computations for the discovery of SIMs are required \citep{kuehn2015multiple}. However, they can also be applied to singularly perturbed systems with explicit timescale splitting, resulting in equivalent SIM approximations as those provided by analytical GSPT methods \citep{zagaris2004analysis,kuehn2015multiple}.~For a detailed review of the GSPT methods used to approximate SIMs see, e.g., \citep{goussis2006efficient,ginoux2021slow}.~Finally, SIMs can be computed by traditional approximations, such as the Quasi Steady-State Approximation (QSSA) and the Partial Equilibrium Approximation (PEA), which can be recovered by GSPT methods \citep{bowen1963singular,zagaris2009analysis,goussis2012quasi,patsatzis2023algorithmic}.~At this point, we note, that, the GSPT-based analytical methods, for high order of accuracy, require complicated calculations, which can be intractable as the dimension of the system increases \citep{kuehn2015multiple}; in fact, they have been applied just for low-dimensional systems. Furthermore, such methods, relying on the approximation of SIMs using asymptotic series expansions, are efficient for very small values of the perturbation parameter. On the other hand, the GPST-based numerical methods, such as the CSP, which have been developed to cope with the intractability limitation
, require point-by-point in time numerical estimations and, usually, they don't provide expressions in an explicit form for higher order approximations (see e.g., in the Appendix \ref{app:MM_SIM_CSPexp}).\par
Machine Learning has been also used to construct surrogate nonlinear ROMs from data, including fuzzy systems \cite{siettos2002semiglobal},
nonlinear manifold learning, such as ISOMAP \citep{balasubramanian2002isomap,bollt2007attractor}, Local Linear Embedding \citep{roweis2000nonlinear,papaioannou2022time}, and Diffusion Maps (DMs) \citep{coifman2005geometric,singer2009detecting,lee2020coarse,dsilva2016data,patsatzis2023data}, Autoencoders \citep{chen2018molecular,vlachas2022multiscale}, Koopman operator \cite{williams2015data,bollt2018matching,lusch2018deep,santos2021reduced}, and deep-learning \cite{linot2020deep}.
The above machine-learning based methodologies are data-driven, thus they don't encounter/exploit the knowledge of the physics that can be available in the form of differential equations. Nor they provide explicit forms of SIMs. Blended methods using machine learning to construct ROMs by fitting slow manifold/closures to coarse-grained models based on high-fidelity simulations or by approximating projection basis vectors via regression have been also proposed. Within this framework, in \cite{wan2018machine}, the authors used long-short term memory (LSTMs) neural networks to provide a first-order approximation of the slow manifold describing the kinematics of finite-size spherical particles in arbitrary fluid flows. In \cite{chen2021physics}, the authors used a PIML to construct ROMs for PDEs, based on the Proper Orthogonal Decomposition (POD) method, where the reduced basis is constructed via high-fidelity simulations using FNNs. In \cite{galassi2022adaptive}, the authors used FNNs to learn from data a surrogate neural network model for the projection basis vectors produced by the implementation of the CSP method  for a stiff chemical kinetics problem. Recently, in \cite{lee2022learning}, FNNs and Gaussian Processes were used to provide closures between the fast and the slow variables, in order to construct ROMs in the form of PDEs from high-fidelity Monte Carlo chemotactic simulations.\par 
Here, we address a physics-informed machine learning \cite{raissi2019physics,karniadakis2021physics} (PIML) approach, within the GSPT framework, for deriving analytical functionals of SIMs in an explicit form.~We consider the class of singular perturbed dynamical systems with an explicit timescale splitting, for which the dimension of the SIM and the variables associated with the fast dynamics do not change during the desired, for the construction of the ROM, timeframe. In contrast to the above-mentioned ML schemes that construct surrogate models via regression, our approach provides a functional form that solves a partial differential equation corresponding to the invariance equation; hence does not require input data for the fast variables. The proposed PIML approach: (i) is not limited to approximations that are theoretically valid only locally, as with the traditional GSPT-based asymptotic series expansions, because of the universal approximation property of neural networks, (ii) does not require extended complicated analytical calculations, which may become intractable as the dimension of the system increases, and (iii) results in functionals that can be evaluated at any set of values of the slow variables and the parameter $\epsilon$, thus avoiding the point-by-point numerical estimations that more sophisticated numerical methods (such as CSP) usually require. For the implementation of the proposed scheme, we considered two ML structures, namely a single hidden layer FNN, and a single-hidden layer RPNN. The performance of the proposed scheme was assessed via three benchmark singularly perturbed dynamical systems, namely the Michaelis-Menten (MM) enzyme reaction scheme, the pharmacokinetic/pharmacodynamic Target Mediated Drug Disposition (TMDD) mechanism, and, the 3D Sel'kov model of glycolytic oscillations.~In order to have a straightforward comparison with other traditional GSPT methods, we also derived the analytical approximations of the SIMs in an explicit form in two ways: (a) by taking the regular asymptotic series expansion of the SIMs and then using the invariance equation to determine its terms, and (b) based on the CSP method with one iteration.~We note that for the particular benchmark problems, one-iteration of the CSP procedure provides analytically explicit forms of the SIMs; this may not be the case for other systems and for more iterations of the CSP (see for example in the Appendix \ref{app:MM_SIM_CSPexp}).Finally, we compared the computational costs when implementing symbolic, automatic and numerical differentiation for the quantities required for the learning process. 

The structure of the paper is as follows. We first present the proposed PIML approach using both SLFNNs and RPNNs, and provide the analytical expressions for the derivatives of both schemes that are required for the learning process, i.e. the solution of the IE. 
Then, we present, the three benchmark problems namely the MM and TMDD reaction mechanisms and the 3D Sel'kov model, and we provide the analytical expressions of the corresponding SIMs, based on both GSPT, and CSP with one iteration. Then, we provide the numerical results based on the proposed scheme,  assess and compare its performance with the other schemes. A comparison of the computational costs of symbolic, automatic, and numerical differentiation schemes is also given. Finally, we discuss the pros and cons of the proposed method in comparison with the other traditional approaches and provide some future directions of research. 

\section{Methodology}
\label{sec:Methods}
For the completeness of the presentation, we first provide some basic elements of perturbation theory, relevant to the proposed ML methodology.~Let us first consider the autonomous system of ODEs (that can also result from the discretization of PDEs):
\begin{equation}
    \dfrac{d \mathbf{z}}{dt} = \mathbf{F}(\mathbf{z}), 
    \qquad \mathbf{z}(t_0) = \mathbf{z_0},
    \label{eq:gen}
\end{equation}
where $\mathbf{z}\in \mathbb{R}^N$ is the $N$-dim. state vector, $\mathbf{F}:\mathbb{R}^N \rightarrow \mathbb{R}^N$ is a smooth vector field, and $\mathbf{z_0}\in \mathbb{R}^N$ is the vector of initial conditions.~In the presence of multiple timescales, various non-dimensional forms of \eqref{eq:gen} exist, for describing its evolution in different timescales.~In the framework of the Singular Perturbation Theory (SPT) \citep{verhulst2005methods,kevorkian2013perturbation, kuehn2015multiple}, the identification of a small parameter $\epsilon=\tau_f/\tau_s\ll1$, expressing the gap between the characteristic fast, say $\tau_f$, and slow, say $\tau_s$, timescales, is firstly required in order to obtain the so-called \emph{fast-slow} subsystems.~Then, the state variables $\mathbf{z}\in\mathbb{R}^N$ in Eq.~\eqref{eq:gen} are linearly transformed into $M$ fast variables (the ones mostly associated with the fast timescales) and $N-M$ remaining slow ones, each denoted as $\mathbf{x}\in \mathbb{R}^M$ and $\mathbf{y}\in \mathbb{R}^{N-M}$, respectively.~Following this decomposition, the vector field $\mathbf{F}(\mathbf{z})$ is transformed into the fast and slow vector fields $\mathbf{f}: \mathbb{R}^N \times I \rightarrow \mathbb{R}^M$ and $\mathbf{g}: \mathbb{R}^N \times I \rightarrow \mathbb{R}^{N-M}$, respectively, where $I\subset \mathbb{R}$ is an interval containing $\epsilon=0$.~Given the above transformations, the original system in Eq.~\eqref{eq:gen} can be cast to its \emph{fast} subsystem form:
\begin{equation}
    \dfrac{d\mathbf{x}}{dt_f} = \mathbf{f}(\mathbf{x},\mathbf{y},\epsilon), \qquad \dfrac{d \mathbf{y}}{dt_f} = \epsilon   \mathbf{g}(\mathbf{x},\mathbf{y},\epsilon), 
    \label{eq:SPfast}
\end{equation}
that is appropriate for describing the system's \emph{fast}
dynamics in the range of the fast timescale $\tau_f$.~Setting $t_s = \epsilon t_f$, one obtains the \emph{slow} subsystem form:
\begin{equation}
    \epsilon \dfrac{d \mathbf{x}}{dt_s} =  \mathbf{f}(\mathbf{x},\mathbf{y},\epsilon), \qquad
    \dfrac{d \mathbf{y}}{dt_s} =   \mathbf{g}(\mathbf{x},\mathbf{y},\epsilon),
\label{eq:SPslow}
\end{equation}
that is appropriate for describing the system's  \emph{slow} dynamics in the range of the slow timescale $\tau_s$.~We emphasize  that the above singularly perturbed systems describe the fast/slow dynamics of the original system, when the variables, associated with the fast timescales, as well as their number $M$, do not change during the timeframe of interest; i.e., in the range of $\tau_f$ or $\tau_s$.~In the case where one of the above assumptions is violated, the fast/slow subsystems are no longer valid (for reproducing the dynamics of the original system) and new transformations are required.


For the implementation of the SPT, the \emph{inner}/\emph{outer} solutions of the fast/slow subsystems in Eqs.~\eqref{eq:SPfast}/\eqref{eq:SPslow} are extracted as asymptotic regular expansions of desired order of accuracy, with their coefficients being powers of $\epsilon$.~Note that the inner solution, frequently referred as \emph{boundary layer}, depends on the initial conditions.~For obtaining a global asymptotic solution of the system, one matches the inner and outer solutions in the timeframe when the two solutions overlap (i.e., when $\tau_f \le t \le \tau_s$).~We highlight here that the SPT is a very cumbersome procedure to apply, especially when the system under study is high-dimensional and complicated.

\subsection{Geometric Singular Perturbation Theory and the Invariance Equation}
In an attempt to overcome the limitations of the SPT, the \emph{Geometric Singular Perturbation Theory} (GSPT) was developed \citep{tikhonov1952systems,fenichel1979geometric} aiming to take into account the dynamical geometric structures arising in the phase space from systems in the form of Eqs.~(\ref{eq:SPfast}, \ref{eq:SPslow}) and their properties; a detailed review and description of GSPT can be found in \citep{kaper1999introduction,verhulst2005methods,kuehn2015multiple,jones1995geometric}.\par
In the asymptotic limit $\epsilon= 0$, the slow subsystem in Eq.~\eqref{eq:SPslow} is confined to evolve onto a \emph{critical manifold} $C_0=\{ (\mathbf{x},\mathbf{y})\in \mathbb{R}^N: \mathbf{f}(\mathbf{x},\mathbf{y},0)=0 \}$, which includes all the equilibrium points of the fast subsystem in Eq.~\eqref{eq:SPfast}.~In addition, a subset $S_0\subset C_0$ is said to be \emph{normally hyperbolic}, if the Jacobian matrix $\nabla_x \mathbf{f}(\mathbf{x},\mathbf{y},0)$ has no eigenvalues with zero real part for every $(\mathbf{x},\mathbf{y}) \in S_0$.~Under this assumption, the Implicit Function Theorem \citep{rudin1976principles} implies the existence of a map $\mathbf{h}_0: \mathbb{R}^{N-M}\rightarrow \mathbb{R}^M$, such that the \emph{compact normally hyperbolic submanifold} $S_0$ can be locally written as 
$$S_0 = \{ (\mathbf{x},\mathbf{y})\in \mathbb{R}^N: \mathbf{x} = \mathbf{h}_0(\mathbf{y}) \}.$$
In the cases of small, yet non-zero, $0<\epsilon \ll1$, according to the Fenichel-Tikhonov theorem \citep{tikhonov1952systems,fenichel1979geometric}, it exists a \emph{locally invariant and normally hyperbolic manifold} $S_\epsilon$ that is diffeomorphic to $S_0 \subset C_0$.~The so-called \emph{Slow Invariant Manifold} (SIM) $S_\epsilon$ has the following properties: (i) has a Hausdorff distance $\mathcal{O}(\epsilon)$ from $S_0$, (ii) the flow on it converges to the slow flow of Eq.~\eqref{eq:SPslow} on $S_0$ as $\epsilon\rightarrow0$, and (iii) is normally hyperbolic with the same stability properties with respect to the fast variables as $S_0$.
$S_\epsilon$ can be locally described as \citep{fenichel1979geometric,jones1995geometric}:
\begin{equation}
S_\epsilon = \{ (\mathbf{x},\mathbf{y})\in \mathbb{R}^N: \mathbf{x} = \mathbf{h}(\mathbf{y},\epsilon) \}.
\label{eq:SIM00}
\end{equation}
The map $\mathbf{x}=\mathbf{h}(\mathbf{y},\epsilon)$ can be plugged in the differential equations of the slow variables $\mathbf{y}$ in Eq.~\eqref{eq:SPslow} for describing the slow dynamics of the full system on the SIM; i.e.,  to get: $d\mathbf{y}/dt=\mathbf{g}(\mathbf{h}(\mathbf{y},\epsilon),\mathbf{y},\epsilon)$.

For the discovery of analytical SIM approximations in the form of Eq.~\eqref{eq:SIM00}, a wide variety of GSPT methods exploit the local invariance property, implying that the map $\mathbf{x}=\mathbf{h}(\mathbf{y},\epsilon) \in \mathbb{R}^{N-M}\times \mathbb{R}  \rightarrow \mathbb{R}^M$ satisfies the \emph{invariance property} corresponding to a (partial) differential equation (PDE) \cite{guckenheimer2013nonlinear,ginoux2021slow}:
\begin{equation}
   \epsilon \nabla_\mathbf{y}\mathbf{h}(\mathbf{y},\epsilon) \mathbf{g}(\mathbf{h}(\mathbf{y},\epsilon),\mathbf{y},\epsilon)=\mathbf{f}(\mathbf{h}(\mathbf{y},\epsilon),\mathbf{y},\epsilon),
    \label{eq:Inv}
\end{equation} 
with appropriate boundary conditions at the boundary, say, $\partial \Omega$, of the manifold \cite{guckenheimer2013nonlinear,ginoux2021slow}. 
As the solution of the above PDE is usually a difficult task (see the discussion in  \cite{guckenheimer2013nonlinear,ginoux2021slow}), approximations of SIMs are frequently derived by asymptotic series expansions of  $\mathbf{h}(\mathbf{y},\epsilon)$ around $S_0$, as:
\begin{equation}
\mathbf{x}=\mathbf{h}(\mathbf{y},\epsilon) = \mathbf{h}_0(\mathbf{y}) + \epsilon \mathbf{h}_1(\mathbf{y}) + \epsilon^2 \mathbf{h}_2(\mathbf{y}) + \ldots + \epsilon^q \mathbf{h}_q(\mathbf{y}) + \mathcal{O}(\epsilon^{q+1}).
    \label{eq:SIMrexp}
\end{equation}
As it has been shown \cite{fenichel1979geometric,guckenheimer2013nonlinear}, such a Taylor series expansion, can approximate the SIM arbitrarily close, around the equilibria of the system for $\epsilon\ll 1$; i.e., around $S_0$.
The most common technique for obtaining $\mathbf{h}_k(\mathbf{y})$ for $k=1,\ldots,q$ involves \citep{fenichel1979geometric,kuehn2015multiple}: first the substitution of Eq.~\eqref{eq:SIMrexp} into the slow subsystem in Eq.~\eqref{eq:SPslow}, then the expansion of the vector fields around $\epsilon=0$, so that $\mathbf{f}(\mathbf{y},\epsilon)=\sum_{k=0}^q \epsilon^k \mathbf{f}_k(\mathbf{y})$ and $\mathbf{g}(\mathbf{y},\epsilon)=\sum_{k=0}^q \epsilon^k \mathbf{g}_k(\mathbf{y})$, and finally matching order-by-order the terms of the invariance equation~\eqref{eq:Inv} to determine  $\mathbf{h}_k(\mathbf{y})$ for $k=1,\ldots,q$.~A detailed presentation of the above procedure is presented in Appendix~\ref{app:A}. 
%
We re-iterate that as the dimension of the system increases the analytical calculation become intractable, thus limiting the application of the above methods to low-dimensional systems \citep{kuehn2015multiple}.\par
SIM approximations can be also derived by sophisticated computational methods in the context of GSPT, such as the computational singular perturbation (CSP) \citep{lam1989understanding,goussis1992study,valorani2005higher}, the invariant low-dimensional manifold (ILDM) \citep{maas1992simplifying}, the zero-derivative principle (ZDP) \citep{gear2005projecting,zagaris2009analysis} and the tangential stretching rate (TSR) \citep{valorani2015dynamical} methods.~The above methods can be applied directly to the original form of the system in Eq.~\eqref{eq:gen} to compute numerical approximations by following different iterative procedures (based on the invariance equation, such as CSP and ZDP, or on the decomposition of the fast and slow subspaces resolving the tangent space, such as CSP, ILDM and TSR).~However, due to their numerical nature, these methods require point-by-point estimations along the trajectory of the system in Eq.~\eqref{eq:gen} \citep{goussis2006efficient,kuehn2015multiple}. 
While originally developed for high-dimensional systems in the form of Eq.~\eqref{eq:gen}, the above methods can deal with singularly perturbed systems in the slow subsystem form of Eq.~\eqref{eq:SPslow}, leading to SIM approximations of desired order of accuracy, depending on the number of iterations \citep{zagaris2004analysis,goussis2006efficient,kuehn2015multiple}.~With a low number of iterations, the computational GSPT methods may result, especially for systems with a low number of variables, to analytic SIM approximations.~Nevertheless, there is no guarantee that these expressions can be written in an explicit form, rather in an implicit one (i.e., $\mathbf{h}(\mathbf{x},\mathbf{y},\epsilon)=\mathbf{0}$).~In this work, the CSP method was applied using one iteration, which results in explicit SIM approximations for the systems under study.~The CSP algorithmic procedure for the derivation of SIM approximations is presented in detail in Appendix \ref{app:B}. 

\subsection{The proposed Physics-Informed Machine Learning (PIML) methodology}
\label{sub:PIMLgen}
Here, we propose a physics-informed Machine Learning (PIML) approach for the discovery of SIM approximations, that are explicitly expressed in terms of the fast variables (i.e., in the form of Eq.~\eqref{eq:SIM00}) by solving the invariance equation via the proposed PIML scheme. 
We begin by assuming a valid fast-slow system in the form of Eq.~\eqref{eq:SPslow} exhibiting a slow evolution on a $(N-M)$-dim. SIM $S_{\epsilon}$, which can be locally approximated by the map $\mathbf{x} = \mathbf{h}(\mathbf{y},\epsilon)$ in Eq.~\eqref{eq:SIM00}.~Let's now take a set of $n_y$ points of the slow variables $\mathbf{y} \in \Omega \subset \mathbb{R}^{N-M}$, and a set of $n_{\epsilon}$ points $\epsilon \in I \subset \mathbb{R}$ in $I=[\epsilon_0,\epsilon_{end}]$ domain.~Then, the numerical approximation of the SIM can be obtained via PIML for the solution of the IE as a minimization problem of the form:
\begin{equation}
    \min_{\mathbf{P},\mathbf{Q}} E (\mathbf{P},\mathbf{Q}):= \sum_{i=1}^{n_y} \sum_{j=1}^{n_{\epsilon}} \big\lVert \mathbf{f}(\mathcal{N}(\mathbf{y}_{i},\epsilon_{j},\mathbf{P},\mathbf{Q}),\mathbf{y}_{i},\epsilon_{j}) - \epsilon_{j} ~ \nabla_\mathbf{y}\mathcal{N}(\mathbf{y}_{i},\epsilon_{j},\mathbf{P},\mathbf{Q}) ~ \mathbf{g}(\mathcal{N}(\mathbf{y}_{i},\epsilon_{j},\mathbf{P},\mathbf{Q}),\mathbf{y}_{i},\epsilon_{j})  \big\rVert^2,
    \label{eq:minFunPI}
\end{equation}
where $\mathcal{N}(\cdot):=\mathcal{N}(\mathbf{y},\epsilon,\mathbf{P},\mathbf{Q}): \mathbb{R}^{N-M} \times \mathbb{R} \rightarrow \mathbb{R}^M$ approximates the output $\mathbf{x} = \mathbf{h}(\mathbf{y},\epsilon)$ of the $M$ algebraic equations of the SIM in the domain $\Omega \times I$.\par 
$\mathcal{N}(\cdot)$ contains the parameters $\mathbf{P}$ of the ML structure (e.g., for a Neural Network, the weights and biases of the layers) and the hyper-parameters $\mathbf{Q}$ such as the type and parameters of the activation function, the learning rate, the number of epochs, etc.~Note that the solution of the optimization problem in Eq.~\eqref{eq:minFunPI}, requires the derivatives of $\mathcal{N}(\mathbf{y},\epsilon,\mathbf{P},\mathbf{Q})$ with respect to the slow variables $\mathbf{y}$ and the parameters $\mathbf{P}$, which can be obtained either numerically (e.g. using finite differences), symbolic or automatic differentiation \citep{baydin2018automatic,lu2021deepxde}.~Among the various structures, we used single-hidden layer FNNs (SLFNNs) and random projection neural networks (RPNNs) \cite{fabiani2022parsimonious} with sigmoid activation functions.

\subsubsection{Solution of the invariance equation with SLFNNs}
\label{subsub:SLFNN_PIML}
Since the inputs to the neural network contain the slow variables $\mathbf{y}\in\mathbb{R}^{N-M}$, and the parameter $\epsilon\in\mathbb{R}$, the input dimension of the SLFNN is $D=N-M+1$.~For $L$ neurons in the hidden layer, the $m$-th output of the SLFNN for $m=1,\ldots,M$ can be written as:
\begin{equation}
\mathcal{N}^{(m)}(\mathbf{y},\epsilon,\mathbf{P}^{(m)}) = \sum_{l=1}^L w^{o(m)}_{l} \mathcal{\phi}_l\left(\sum_{d=1}^{D-1} w^{(m)}_{ld} y_d + w_{lD}^{(m)} \epsilon + b^{(m)}_l\right)+ b^{o(m)} = \mathbf{w}^{o(m)\top} \phi \left( \mathbf{W}^{(m)} \begin{bmatrix} \mathbf{y} \\ \epsilon \end{bmatrix}   + \mathbf{b}^{(m)} \right)  + b^{o(m)},
    \label{eq:SFLNNsumf}
\end{equation}
where the learning parameters $\mathbf{P}^{(m)}$ of the $m$-th output of the SLFNN are: (i) the vector of the output weights $\mathbf{w}^{o(m)} = [w^{o(m)}_1, \ldots, w^{o(m)}_L]^\top\in \mathbb{R}^{L}$ between the neurons in the hidden layer and the output layer, (ii) the bias $b^{o(m)} \in \mathbb{R}$ of the output layer, (iii) the matrix of the internal weights $\mathbf{W}^{(m)}\in\mathbb{R}^{L\times D}$ between the input layer and the hidden layer, whose columns are the vectors $\mathbf{w}^{(m)}_l = [w_{l1}^{(m)}, \ldots, w_{lD}^{(m)}]^\top\in \mathbb{R}^{D}$ corresponding to the weights between the nodes of the input layer and the $l$-th neuron in the hidden layer, and (iv) the vector of the internal biases  $\mathbf{b}^{(m)} = [b_1^{(m)}, \ldots, b_L^{(m)}]^\top\in\mathbb{R}^{L}$ of the neurons in the hidden layer.~Note that the first $D-1$ elements of  the internal weights $\mathbf{w}^{(m)}_l$ correspond to the slow variables in $\mathbf{y}=[y_1,\ldots,y_d,\ldots,y_{N-M}]^\top$, while the last element $w^{(m)}_{lD}$ corresponds to $\epsilon$.~Here, the activation function $\phi_l(\cdot) = \phi(\mathbf{w}_l^{(m)\top} (\mathbf{y}, \epsilon)^\top  + b_l^{(m)})$ is the logistic sigmoid function that can facilitate symbolic differentiation of the required derivatives of the SLFNN for the optimization problem.\par
For the solution of the optimization problem in Eq.~\eqref{eq:minFunPI}, we first collect the $M$ outputs of the SLFNN in Eq.~\eqref{eq:SFLNNsumf} in the column vector:
\begin{equation}
    \mathcal{N}(\mathbf{y}_i,\epsilon_j,\mathbf{P})=\begin{bmatrix} \mathcal{N}^{(1)}(\mathbf{y}_i,\epsilon_j,\mathbf{P}^{(1)}) & \ldots & \mathcal{N}^{(m)}(\mathbf{y}_i,\epsilon_j,\mathbf{P}^{(m)}) & \ldots & \mathcal{N}^{(M)}(\mathbf{y}_i,\epsilon_j,\mathbf{P}^{(M)}) \end{bmatrix}^\top 
    \label{eq:SLFNNcol}
\end{equation}
for every input point $(\mathbf{y}_i,\epsilon_j)$ with $i=1,\ldots,n_y$ and $j=1,\ldots,n_{\epsilon}$.\par 
Then, the PIML optimization problem reduces to the minimization of the loss function:
\begin{equation}
    \mathcal{L}(\mathbf{P}) = \sum_{i=1}^{n_y} \sum_{j=1}^{n_{\epsilon}} \big\lVert \mathbf{f}(\mathcal{N}(\mathbf{y}_i,\epsilon_j,\mathbf{P}),\mathbf{y}_{i},\epsilon_{j}) - \epsilon_{j} ~ \nabla_{\mathbf{y}}\big(\mathcal{N}(\mathbf{y}_i,\epsilon_j,\mathbf{P})\big) ~ \mathbf{g}(\mathcal{N}(\mathbf{y}_i,\epsilon_j,\mathbf{P}),\mathbf{y}_{i},\epsilon_{j}) \big\rVert^2,
     \label{eq:SLFNN_PIresM}
\end{equation}
with respect to all the parameters $\mathbf{P}^{(m)}=[\mathbf{w}^{o(m)},b^{o(m)},\mathbf{W}^{(m)},\mathbf{b}^{(m)}]^\top\in\mathbb{R}^{L(D+2)+1}$ of the $m$-th SLFNN output contained in $\mathbf{P}$.\par 
To minimize of the loss function in Eq.~\eqref{eq:SLFNN_PIresM}, the minimization of the $M \times n_y \times  n_{\epsilon}$ non-linear residuals: 
\begin{equation}
    \mathcal{F}_q(\mathbf{P}) = f_m(\mathcal{N}(\mathbf{y}_i,\epsilon_j,\mathbf{P}),\mathbf{y}_i,\epsilon_j) -  \epsilon_j \sum_{d=1}^{N-M} \dfrac{\partial \mathcal{N}^{(m)}(\mathbf{y}_i,\epsilon_j,\mathbf{P}^{(m)})}{ \partial y_d} g_d(\mathcal{N}(\mathbf{y}_i,\epsilon_j,\mathbf{P}),\mathbf{y}_i,\epsilon_j)
    \label{eq:SLFNN_PIresM2}
\end{equation}
is required, where for every pair $(i,j)$ of input points $m=1,\ldots,M$ residuals are formed, such that $q=m+(i-1+(j-1)n_y)M$.~The terms $f_m(\cdot)$ and $g_d(\cdot)$ in Eq.~\eqref{eq:SLFNN_PIresM2}  denote the $m$-th and $d$-th components of the  analytically known fast and slow vector fields $\mathbf{f}(\cdot)=[f_1(\cdot), \ldots, f_m(\cdot),\ldots,f_M(\cdot)]^\top$ and $\mathbf{g}(\cdot) = [g_1(\cdot), \ldots, g_d(\cdot),\ldots,g_{N-M}(\cdot)]^\top$ in Eq.~\eqref{eq:SPslow}, respectively.\par 
The formation of the residuals in Eq.~\eqref{eq:SLFNN_PIresM2} additionally requires the calculation of the $M\times(N-M)$ derivatives $\partial \mathcal{N}^{(m)}(\mathbf{y}_i,\epsilon_j,\mathbf{P}^{(m)})/ \partial y_d$, which can be calculated through symbolic differentiation, given the derivative of the sigmoid function $\phi'=\phi(1-\phi)$.~According to Eq.~\eqref{eq:SFLNNsumf}, the derivative of the $m$-th SLFNN output w.r.t. the $d$-th slow variable is:
\begin{equation}
    \dfrac{\partial \mathcal{N}^{(m)}(\mathbf{y}_i,\epsilon_j,\mathbf{P}^{(m)})}{ \partial y_d} = \sum_{l=1}^L  w^{o(m)}_l w^{(m)}_{ld}\left[ \phi_l(\cdot) ( 1 - \phi_l (\cdot ) ) \right], \qquad \phi_l(\cdot) = \phi_l(\mathbf{w}_{l}^{(m)\top} \begin{bmatrix} \mathbf{y}_i \\ \epsilon_j \end{bmatrix} + b^{(m)}_l),
    \label{eq:SFLNNder_sumf}
\end{equation}
where $d=1,\ldots,N-M$.\par 
Given Eqs.~(\ref{eq:SFLNNsumf}, \ref{eq:SLFNNcol}, \ref{eq:SFLNNder_sumf}), all the components for the calculation of the non-linear residuals in Eq.~\eqref{eq:SLFNN_PIresM2} are now available.~In general, the PIML optimization problem is overdetermined, since the number of residuals $F_q$ is larger than the number of the tunable parameters of the SLFNN; i.e., $Mn_y n_{\epsilon} \gg M(L(D+2)+1)$.~Here, for the minimization of the non-linear residuals in Eq.~\eqref{eq:SLFNN_PIresM2}, we used the Levenberg-Marquardt iterative algorithm \citep{hagan1994training} for the update of the learnable parameters $\mathbf{P}^{(m)}$, implemented in MATLAB R2022b, which is described below. 

The implementation of the Levenberg-Marquardt algorithm requires the Jacobian matrix of the residuals with respect to the SLFNN weights and biases, which can be calculated with symbolic differentiation as follows.~Stacking the residuals into the column vector $\mathbf{F}(\mathbf{P})=[\mathcal{F}_1(\mathbf{P}), \ldots, \mathcal{F}_q(\mathbf{P}),\ldots,\mathcal{F}_{Mn_yn_{\epsilon}}(\mathbf{P})]^\top$ and the learnable parameters into $\mathbf{P}=[\mathbf{P}^{(1)},\ldots,\mathbf{P}^{(r)},\ldots,\mathbf{P}^{(M)}]^\top\in\mathbb{R}^{M(L(D+2)+1)}$, the elements of the Jacobian matrix $\nabla_{\mathbf{P}} \mathbf{F}\in \mathbb{R}^{Mn_yn_\epsilon \times M(L(D+2)+1)}$ over any learnable parameter $p\in\mathbf{P}$ are: 
\begin{equation}
    \dfrac{\partial F_q}{\partial p} = \dfrac{\partial f_m (\cdot)}{\partial p}  - \epsilon_j \sum_{d=1}^{N-M} \left( \dfrac{\partial^2 \mathcal{N}^{(m)}(\mathbf{y}_i,\epsilon_j,\mathbf{P}^{(m)})}{ \partial p \partial y_d} g_d(\cdot) +  \dfrac{\partial \mathcal{N}^{(m)}(\mathbf{y}_i,\epsilon_j,\mathbf{P}^{(m)})}{ \partial y_d} \dfrac{\partial g_d (\cdot)}{\partial p}\right),
    \label{eq:SLFNN_NLres_grad}
\end{equation}
where for the calculation of $f_m(\cdot)=f_m(\mathcal{N}(\mathbf{y}_i,\epsilon_j,\mathbf{P}))$ and $g_d(\cdot)=g_d(\mathcal{N}(\mathbf{y}_i,\epsilon_j,\mathbf{P}))$, all the $\mathcal{N}^{(m)}(\mathbf{y}_i,\epsilon_j,\mathbf{P}^{(m)})$ outputs of the SLFNN for $m=1,\ldots,M$ are, in general, required.~Since the analytical expressions of $f_m(\cdot)$, $g_d(\cdot)$ are known from Eq.~\eqref{eq:SPslow}, the first order derivatives in Eq.~\eqref{eq:SLFNN_NLres_grad} can be calculated, as:
\begin{equation}
\dfrac{\partial f_m(\cdot)}{\partial p} = \dfrac{\partial f_m(\cdot)}{\partial x_r}  \dfrac{\partial \mathcal{N}^{(r)}(\mathbf{y}_i,\epsilon_j,\mathbf{P}^{(r)})}{\partial p}, \qquad
\dfrac{\partial g_d(\cdot)}{\partial p} = \dfrac{\partial g_d(\cdot)}{\partial x_r} \dfrac{\partial \mathcal{N}^{(r)}(\mathbf{y}_i,\epsilon_j,\mathbf{P}^{(r)})}{\partial p}, \label{eq:SLFNN_NLres_grad1}
\end{equation}
where $\partial f_m/\partial x_r$ and $\partial g_d/\partial x_r$ denote the derivatives of the system in Eq.~\eqref{eq:SPslow} w.r.t. the $r$-th fast variable $x_r$ in $\mathbf{x}=[x_1,\ldots,x_r,\ldots,x_M]^\top$ with $r=1,\ldots,M$.\par
Furthermore, using Eq.~\eqref{eq:SFLNNder_sumf}, the mixed derivative term involved in Eq.~\eqref{eq:SLFNN_NLres_grad} is given by:
\begin{equation}
\dfrac{\partial^2 \mathcal{N}^{(m)}(\mathbf{y}_i,\epsilon_j,\mathbf{P}^{(m)})}{ \partial p \partial y_d} = \dfrac{\partial\sum_{l=1}^L  w^{o(m)}_l w^{(m)}_{ld}\left[ \phi_l(\cdot) ( 1 - \phi_l (\cdot ) ) \right]}{\partial p}, \qquad \phi_l(\cdot) = \phi_l(\mathbf{w}_{l}^{(m)\top} \begin{bmatrix} \mathbf{y}_i \\ \epsilon_j \end{bmatrix} + b^{(m)}_l),
\label{eq:SLFNN_NLres_grad2}
\end{equation}
where $d=1,\ldots,N-M$.\par 
For calculating the derivatives w.r.t. the SLFNN weights and biases, symbolic differentiation and the derivatives of the logistic sigmoid function $\phi'=\phi(1-\phi)$ and $\phi''=\phi(1-\phi)(1-2\phi)$ are used, implying: 
\begin{itemize}
\item for $p=w_l^{o(r)}$:
\begin{equation}
\dfrac{\partial \mathcal{N}^{(r)}(\mathbf{y}_i,\epsilon_j,\mathbf{P}^{(r)})}{\partial w^{o(r)}_l} = \phi_l(\cdot),
\quad \dfrac{\partial^2 \mathcal{N}^{(m)}(\mathbf{y}_i,\epsilon_j,\mathbf{P}^{(m)})}{ \partial w_l^{o(r)} \partial y_d} = 
    \begin{cases}
        w_{ld}^{(r)}\left[\phi_l(\cdot)(1-\phi_l(\cdot))\right], & \text{if } r = m\\
        0, & \text{if } r \neq m
    \end{cases}     
\end{equation}
\item for $p=b^{o(r)}$:
\begin{equation}
\dfrac{\partial \mathcal{N}^{(r)}(\mathbf{y}_i,\epsilon_j,\mathbf{P}^{(r)})}{\partial b^{o(r)}} = 1,
\quad \dfrac{\partial^2 \mathcal{N}^{(m)}(\mathbf{y}_i,\epsilon_j,\mathbf{P}^{(m)})}{ \partial b^{o(r)} \partial y_d} = 0 
\end{equation}
\item for $p=w_{lh}^{(r)}$ with $h=1,\ldots,N-M$ (i.e., the internal weights for the slow variables $\mathbf{y}$):
\begin{align}
\dfrac{\partial \mathcal{N}^{(r)}(\mathbf{y}_i,\epsilon_j,\mathbf{P}^{(r)})}{\partial w_{lh}^{(r)}} & = w_l^{o(r)} y_h \left[\phi_l(\cdot)(1-\phi_l(\cdot))\right] \nonumber \\
\dfrac{\partial^2 \mathcal{N}^{(m)}(\mathbf{y}_i,\epsilon_j,\mathbf{P}^{(m)})}{\partial w_{lh}^{(r)}  \partial y_d} & = w_l^{o(r)} \left[\phi_l(\cdot)(1-\phi_l(\cdot))\right] \cdot 
    \begin{cases}
         w_{ld}^{(r)}y_h(1-2\phi_l(\cdot))+1, & \text{if } r=m\ , \ h=d  \\
         w_{ld}^{(r)}y_h(1-2\phi_l(\cdot)), & \text{if } r=m\ , \ h \neq d \\
         0, & \text{if } r \neq m
    \end{cases}
\end{align}
and for $p=w_{lD}^{(r)}$ (i.e., the internal weights for $\epsilon$):
\begin{align}
\dfrac{\partial \mathcal{N}^{(r)}(\mathbf{y}_i,\epsilon_j,\mathbf{P}^{(r)})}{\partial w_{lD}^{(r)}} & = w_l^{o(r)} \epsilon \left[\phi_l(\cdot)(1-\phi_l(\cdot))\right] \nonumber \\
\dfrac{\partial^2 \mathcal{N}^{(m)}(\mathbf{y}_i,\epsilon_j,\mathbf{P}^{(m)})}{\partial w_{lD}^{(r)}  \partial y_d} & =       \begin{cases}
        w_l^{o(r)} w_{ld}^{(r)}\epsilon \left[\phi_l(\cdot)(1-\phi_l(\cdot))(1-2\phi_l(\cdot))\right], & \text{if } r=m \\
         0, & \text{if } r \neq m
    \end{cases}
\end{align}
\item for $p=b_l^{(r)}$:
\begin{align}
\dfrac{\partial \mathcal{N}^{(r)}(\mathbf{y}_i,\epsilon_j,\mathbf{P}^{(r)})}{\partial b^{(r)}_l} & = w_l^{o(r)}  \left[\phi_l(\cdot)(1-\phi_l(\cdot))\right] \nonumber \\
\dfrac{\partial^2 \mathcal{N}^{(m)}(\mathbf{y}_i,\epsilon_j,\mathbf{P}^{(m)})}{ \partial b^{(r)}_l \partial y_d} & = \begin{cases}
        w_l^{o(r)} w_{ld}^{(r)} \left[\phi_l(\cdot)(1-\phi_l(\cdot))(1-2\phi_l(\cdot))\right], & \text{if } r=m \\
         0, & \text{if } r \neq m
    \end{cases} 
    \label{eq:LM_Jac}
\end{align}
\end{itemize}
where $r,m=1,\ldots,M$, $l=1,\ldots,L$, $d=1,\ldots,N-M$ and $\phi_l(\cdot) = \phi_l(\mathbf{w}_{l}^{(r)\top} \cdot \left[ \mathbf{y}_i, \ \epsilon_j \right]^\top + b^{(r)}_l)$.~Finally, substitution of Eqs.~(\ref{eq:SLFNN_NLres_grad1}-\ref{eq:LM_Jac}) into Eq.~\eqref{eq:SLFNN_NLres_grad} results in the formulation of the Jacobian matrix $\nabla_{\mathbf{P}} \mathbf{F}$.~We note here, that the Jacobian matrix can be alternatively calculated with finite differences schemes or automatic differentiation. 

For learning $\mathbf{P}$ with the Levenberg-Marquardt iterative algorithm, we began with a random initial guess of the parameters $\mathbf{P}^0$ and the damping factor set to $\lambda_0=0.01$.~At the $\nu$-th iteration, the residual vector $\mathbf{F}(\mathbf{P}^{\nu})$ and the Jacobian matrix $\nabla_{\mathbf{P}^{\nu}} \mathbf{F}$ are computed through Eqs.~\eqref{eq:SLFNN_PIresM2} and \eqref{eq:SLFNN_NLres_grad}, respectively.~Using the estimation $\left(\nabla_{\mathbf{P}^{\nu}} \mathbf{F}\right)^\top \nabla_{\mathbf{P}^{\nu}} \mathbf{F}$ for the $M(L(D+2)+1)\times M(L(D+2)+1)$ Hessian matrix, the Levenberg-Marquardt algorithm computes the search direction (a $M(L(D+2)+1)$-dim. $\mathbf{d}^\nu$ vector) obtained from the solution of the linearized system: 
\begin{equation}
\left( \left(\nabla_{\mathbf{P}^{\nu}} \mathbf{F}\right)^\top \nabla_{\mathbf{P}^{\nu}} \mathbf{F} + \lambda_\nu diag(\left(\nabla_{\mathbf{P}^{\nu}} \mathbf{F}\right)^\top \nabla_{\mathbf{P}^{\nu}} \mathbf{F})\right) \mathbf{d}^\nu = - \left( \nabla_{\mathbf{P}^{\nu}} \mathbf{F}\right)^\top \mathbf{F}(\mathbf{P}^{\nu})
    \label{eq:LMsd}
\end{equation}
where $diag(\cdot)$ denotes the diagonal matrix of the Hessian matrix approximation.~Then, the learnable parameters $\mathbf{P}^{\nu+1}$ and the dumping factor $\lambda_{\nu+1}$ at the next iteration are updated as follows: 
\begin{itemize}
    \item if $\lVert \mathbf{F}(\mathbf{P}^\nu+\mathbf{d}^\nu)\rVert_{l^2}<\lVert\mathbf{F}(\mathbf{P}^\nu)\rVert_{l^2}$ (successful step), then $\mathbf{P}^{\nu+1} = \mathbf{P}^\nu+\mathbf{d}^\nu$ and $\lambda_{\nu+1}=\lambda_{\nu}/10$, or
    \item if $\lVert \mathbf{F}(\mathbf{P}^\nu+\mathbf{d}^\nu)\rVert_{l^2}\ge\lVert\mathbf{F}(\mathbf{P}^\nu)\rVert_{l^2}$ (unuccessful step), then $\mathbf{P}^{\nu+1} = \mathbf{P}^\nu$ and $\lambda_{\nu+1}=10\lambda_{\nu}$.
\end{itemize}
Convergence of the algorithm is achieved when the stopping criterion $\lVert \mathbf{F}(\mathbf{P^{\nu+1}}) \rVert_{l^2}<tol$ is met, where $tol$ depends on the number $L$ of neurons in the hidden layer.

\subsubsection{Solution of the invariance equation with RPNNs}
\label{subsub:RPNN_PIML}
Random Projection Neural Networks (RPNNs) is a general class of neural networks introduced to provide a low computational cost alternative to other ML schemes (for a brief review see \cite{fabiani2022parsimonious}). 
Here, we build RPNNs by considering randomly parameterized activation functions in the hidden layer, as implemented in \citep{fabiani2021numerical,fabiani2022parsimonious}.
The $m$-th output of the RPNN for $m=1,\ldots,M$ is written as:
\begin{equation}
    \mathcal{N}^{(m)}(\mathbf{y},\epsilon,\mathbf{w}^{o(m)}) = \sum_{l=1}^L w^{o(m)}_{l} \mathcal{\phi}_l\left(\sum_{d=1}^{D-1} \alpha^{(m)}_{ld} y_d + \alpha_{lD}^{(m)} \epsilon + \beta^{(m)}_l\right) = \mathbf{w}^{o(m)\top} \boldsymbol{\phi} \left( \mathbf{A}^{(m)} \begin{bmatrix} \mathbf{y} \\ \epsilon \end{bmatrix}   + \boldsymbol{\beta}^{(m)}\right), 
    \label{eq:RPNNsumf}
\end{equation}
where $\boldsymbol{\alpha}^{(m)}_l = [\alpha^{(m)}_{l1},\ldots,\alpha^{(m)}_{lD}]^\top \in \mathbb{R}^{D}$ are the internal weights of the $l$-th neuron for $l=1,\ldots,L$, all together forming the matrix $\mathbf{A}^{(m)} = [\boldsymbol{\alpha}^{(m)\top}_1, \ldots, \boldsymbol{\alpha}_L^{(m)\top}]\in\mathbb{R}^{L\times D}$, and  $\boldsymbol{\beta}^{(m)}=[\beta_1^{(m)},\ldots,\beta_L^{(m)}]^\top \in\mathbb{R}^{L}$ are the internal biases. The RPNN contains also the output weights $\mathbf{w}^{o(m)} \in \mathbb{R}^{L}$ for each $m=1,\ldots,M$, while the output biases are set to zero.~However, the internal weights and biases of the RPNN are random fixed parameters (parsimoniously chosen, as described next), thus leading to the compact expression for the $m$-th output of the RPNN:
\begin{equation}
    \mathcal{N}^{(m)}(\mathbf{y}_i,\epsilon_j,\mathbf{w}^{o(m)})  = \mathbf{w}^{o(m)\top} \boldsymbol{\Phi}^{(m)}, \quad m =1,\ldots,M
    \label{eq:RPNNcol}
\end{equation}
for any input point $\mathbf{y}_i$ for $i=1,\ldots,n_y$ and $\epsilon_j$ for $j=1,\ldots,n_{\epsilon}$, where the only learnable parameters are the output weights $\mathbf{w}^{o(m)}$.~As shown in Eq.~\eqref{eq:RPNNcol}, the output of the RPNN is linearly related to $\mathbf{w}^{o(m)}$, since the random projection matrix $\boldsymbol{\Phi}^{(m)}\in \mathbb{R}^{L \times n_y n_{\epsilon}}$ of the activated outputs in the hidden layer contains no learnable parameters.~In particular, its elements depend only on the input points, such that: 
\begin{equation}
\Phi^{(m)}_{l,k} = \phi(\boldsymbol{\alpha}_l^{(m)\top} \begin{bmatrix} \mathbf{y}_i \\ \epsilon_j \end{bmatrix}  + \boldsymbol{\beta}^{(m)}),  
\label{eq:RPmatrix}
\end{equation}
where $k=i+(j-1)n_y$.
\paragraph{The sampling process}
Before further proceeding with the methodology for solving the PIML optimization problem, we first discuss the random sampling of the internal weights $\boldsymbol{\alpha}_l^{(m)}$ and biases $\boldsymbol{\beta}^{(m)}$.~Although the existing theoretical framework \citep{barron1993universal,igelnik1995stochastic,huang2006extreme,rahimi2007random} suggests that any random choice in the $[-1,1]$ interval should be good enough, in practice a parsimonious selection from appropriate uniform distributions is more convenient \cite{fabiani2021numerical,calabro2021extreme,galaris2022numerical,dong2022computing,dong2022numerical}
~In particular, we considered the logistic sigmoid activation function (as for the SLFNNs) $\phi(y)=1/(1+exp(-y))$, the inflection point of which is at $y=0$.~Given that the input points $(\mathbf{y}_i,\epsilon_j)^\top$ lie in the subset $\Omega \times I\subset \mathbb{R}^D$, we randomly sampled $L$ points, say $\mathbf{c}_{l}\in \Omega \times I$ (for the $\epsilon$ dimension, we sampled logarithmically spaced points in the interval $I$), which are called \emph{centers} of the activation function $\phi$, requiring them to be the activation function's inflection points: 
\begin{equation}
    \boldsymbol{\alpha}_l^{(m)\top} \mathbf{c}_{l} + \beta_l^{(m)} = 0 \Rightarrow \sum_{d=1}^D \alpha_{ld}^{(m)} c_{ld} + \beta_l^{(m)} = 0.
    \label{eq:AFinf}
\end{equation}
We then sampled the values of $\boldsymbol{\alpha}_l^{(m)}$ in the $[-1,1]$ interval.~However, instead of sampling the values of biases  $\beta_l^{(m)}$ in the same way, we determined  them using Eq.~\eqref{eq:AFinf}, in order to ensure that the centers are the inflection points of the activation function.
A detailed rationalization for the above sampling procedure can be found in \citep{galaris2022numerical}.

For the solution of the PIML optimization problem in Eq.~\eqref{eq:minFunPI} with RPNNs, we first collect the $M$ outputs of the RPNN in Eq.~\eqref{eq:RPNNcol} in the column vector
\begin{equation}
    \mathcal{N}(\mathbf{y}_i,\epsilon_j,\mathbf{W}^{o})=\begin{bmatrix} \mathbf{w}^{o(1)\top} \boldsymbol{\Phi}^{(1)} & \ldots & \mathbf{w}^{o(m)\top} \boldsymbol{\Phi}^{(m)} & \ldots & \mathbf{w}^{o(M)\top} \boldsymbol{\Phi}^{(M)}) \end{bmatrix}^\top,
    \label{eq:RPNNcol2}
\end{equation}
where $\mathbf{W}^{o} = [\mathbf{w}^{o(1)}, \ldots, \mathbf{w}^{o(M)}] \in \mathbb{R}^{ML}$ collects all the learnable parameters $\mathbf{w}^{o(m)}$ for $m=1,\ldots,M$ of the RPNN.\par
\paragraph{The solution of the IE}Based on the above, the  solution of the IE reduces to the minimization of the loss function:
\begin{equation}
    \mathcal{L}(\mathbf{W}^{o}) = \sum_{i=1}^{n_y} \sum_{j=1}^{n_{\epsilon}} \big\lVert \mathbf{f}(\mathcal{N}(\mathbf{y}_i,\epsilon_j,\mathbf{W}^{o}),\mathbf{y}_{i},\epsilon_{j}) - \epsilon_{j} ~ \nabla_{\mathbf{y}}\big(\mathcal{N}(\mathbf{y}_i,\epsilon_j,\mathbf{W}^{o})\big) ~ \mathbf{g}(\mathcal{N}(\mathbf{y}_i,\epsilon_j,\mathbf{W}^{o}),\mathbf{y}_{i},\epsilon_{j}) \big\rVert^2,
     \label{eq:RPNN_PIresM}
\end{equation}
with respect to the parameters $\mathbf{W}^{o}$, with the random projection matrices $\boldsymbol{\Phi}^{(m)}$ fixed for $m=1,\ldots,M$.\par
The minimization of the loss function in Eq.~\eqref{eq:RPNN_PIresM} requires the minimization of the $M \times n_y \times n_{\epsilon}$ non-linear residuals $\mathcal{F}_q$: 
\begin{equation}
    \mathcal{F}_q(\mathbf{W}^o) = f_m(\mathcal{N}(\mathbf{y}_i,\epsilon_j,\mathbf{W}^{o}),\mathbf{y}_i,\epsilon_j) -  \epsilon_j \sum_{d=1}^{N-M} \dfrac{\partial \mathcal{N}^{(m)}(\mathbf{y}_i,\epsilon_j,\mathbf{w}^{o(m)})}{ \partial y_d} g_d(\mathcal{N}(\mathbf{y}_i,\epsilon_j,\mathbf{W}^{o}),\mathbf{y}_i,\epsilon_j),
    \label{eq:RPNN_PIresM2}
\end{equation}
where $q=m+(i-1+(j-1)n_y)M$ with $m=1,\ldots,M$ and $f_m(\cdot)$ and $g_d(\cdot)$ denote the $m$-th and $d$-th components of the  analytically known fast and slow vector fields $\mathbf{f}(\cdot)$ and $\mathbf{g}(\cdot)$ in Eq.~\eqref{eq:SPslow}, respectively.\par
The formation of the residuals in Eq.~\eqref{eq:RPNN_PIresM2} additionally requires the calculation of the $M \times N-M$ derivatives $\partial \mathcal{N}^{(m)}(\mathbf{y}_i,\epsilon_j,\mathbf{w}^{o(m)})/ \partial y_d$, which can be calculated by symbolic differentiation.~Thus, according to Eq.~\eqref{eq:RPNNsumf}, the derivative of the $m$-th RPNN output w.r.t. the $d$-th slow variables reads:
\begin{equation}
    \dfrac{\partial \mathcal{N}^{(m)}(\mathbf{y}_i,\epsilon_j,\mathbf{w}^{o(m)})}{ \partial y_d} = \sum_{l=1}^L  w^{o(m)}_l \alpha^{(m)}_{ld}\left[ \phi_l(\cdot) ( 1 - \phi_l (\cdot ) ) \right], \qquad \phi_l(\cdot) = \phi_l(\boldsymbol{\alpha}_{l}^{(m)\top} \begin{bmatrix}        \mathbf{y}_i \\ \epsilon_j   \end{bmatrix} + \boldsymbol{\beta}^{(m)}),
    \label{eq:RPNNder_y}
\end{equation}
where $d=1,\ldots,N-M$.\par
Equation~\eqref{eq:RPNNder_y} enables the calculation of the non-linear residuals $F_q$ in Eq.~\eqref{eq:RPNN_PIresM2}.~As in the case of SLFNNs, the PIML optimization problem is, generally, overdetermined.~However, in contrast to SLFNNs, the RPNN output is linearly related to the learnable parameters; see Eq.~\eqref{eq:RPNNcol}.~Hence, in this case, the minimization of the non-linear residuals in Eq.~\eqref{eq:RPNN_PIresM2} can be obtained with Newton-type iterative schemes that allow faster convergence \cite{fabiani2021numerical,fabiani2022parsimonious}.

For the implementation of the Newton-type iterative scheme, the Jacobian matrix w.r.t the RPNN output weights $\mathbf{W}^{o}$ is required.~Collecting the residuals in the column vector $\mathbf{F}(\mathbf{W}^{o})=[\mathcal{F}_1(\mathbf{W}^o), \ldots, \mathcal{F}_q(\mathbf{W}^o),\ldots,\mathcal{F}_{Mn_yn_{\epsilon}}(\mathbf{W}^o)]^\top$, the elements of the Jacobian matrix $\nabla_{\mathbf{W}^{o}} \mathbf{F}\in \mathbb{R}^{M n_y n_{\epsilon} \times ML}$ are calculated, through symbolic differentiation of  Eq.~\eqref{eq:RPNN_PIresM2}, as:
\begin{equation}
    \dfrac{\partial F_q}{\partial W^{o}_p} = \dfrac{\partial f_m (\cdot)}{\partial w^{o(r)}_l}  - \epsilon_j \sum_{d=1}^{N-M} \left( \dfrac{\partial^2 \mathcal{N}^{(m)}(\mathbf{y}_i,\epsilon_j,\mathbf{w}^{o(m)})}{ \partial w^{o(r)}_l \partial y_d} g_d(\cdot) +  \dfrac{\partial \mathcal{N}^{(m)}(\mathbf{y}_i,\epsilon_j,\mathbf{w}^{o(m)})}{ \partial y_d} \dfrac{\partial g_d (\cdot)}{\partial w^{o(r)}_l}\right),
    \label{eq:NR_Jac}
\end{equation}
where $p = l+(r-1)L$ with $r=1,\ldots,M$ and $l=1,\ldots,L$.~Since the analytic expressions of $f_m(\cdot)$/$g_d(\cdot)$ in Eq.~\eqref{eq:SPslow} are known, the first-order derivatives involved in Eq.~\eqref{eq:NR_Jac} can be calculated analytically using Eq.~\eqref{eq:RPNNsumf} as: 
\begin{equation}
\dfrac{\partial f_m(\cdot)}{\partial w^{o(r)}_l} =  \dfrac{\partial f_m(\cdot)}{\partial x_r} \phi_l(\boldsymbol{\alpha}_{l}^{(r)\top} \begin{bmatrix}        \mathbf{y}_i \\ \epsilon_j   \end{bmatrix} + \boldsymbol{\beta}^{(r)}), \qquad
\dfrac{\partial g_d(\cdot)}{\partial w^{o(r)}_l} = \dfrac{\partial g_d(\cdot)}{\partial x_r} \phi_l(\boldsymbol{\alpha}_{l}^{(r)\top} \begin{bmatrix}        \mathbf{y}_i \\ \epsilon_j   \end{bmatrix} + \boldsymbol{\beta}^{(r)}), \label{eq:NR_Jac1}
\end{equation}
where $\partial f_m/\partial x_r$ and $\partial g_d/\partial x_r$ denote the derivatives of the system in Eq.~\eqref{eq:SPslow} w.r.t. the $r$-th fast variable $x_r$.~Similarly, from Eq.~\eqref{eq:RPNNder_y}, the mixed derivative term involved in Eq.~\eqref{eq:NR_Jac} reads:
\begin{equation}
    \dfrac{\partial^2 \mathcal{N}^{(m)}(\mathbf{y}_i,\epsilon_j,\mathbf{w}^{o(m)})}{\partial w_l^{o(r)} \partial y_d} = 
    \begin{cases}
        \alpha^{(r)}_{ld}\left[ \phi_l(\cdot) ( 1 - \phi_l (\cdot ) ) \right], & \text{if } r = m\\
        0, & \text{if } r\neq m
    \end{cases} 
    , \qquad \phi_l(\cdot) = \phi_l(\boldsymbol{\alpha}_{l}^{(r)} (\mathbf{y}_i, \epsilon_j)^T + \boldsymbol{\beta}^{(r)}),
    \label{eq:NR_Jac2}
\end{equation}
where $r,m=1,\ldots,M$, $l=1,\ldots,L$ and $d=1,\ldots,N-M$.~Finally, the Jacobian matrix $\nabla_{\mathbf{W}^{o}} \mathbf{F}$ is formulated with the use of Eqs.~(\ref{eq:NR_Jac1}, \ref{eq:NR_Jac2}) into Eq.~\eqref{eq:NR_Jac}.

For training the RPNN with the Newton iterative scheme, we begin with a random initial guess of the output weights $\mathbf{W}^{o(0)}$.~At the $\nu$-th iteration, the residual vector $\mathbf{F}(\mathbf{W}^{o(\nu)})$ and the Jacobian matrix $\nabla_{\mathbf{W}^{o(\nu)}} \mathbf{F}$ are computed through Eqs.~(\ref{eq:RPNN_PIresM2}) and \eqref{eq:NR_Jac}, respectively.~Then, the update $d\mathbf{W}^{o(\nu)} \in \mathbb{R}^{ML}$ at the $\nu$-th  iteration is computed via the solution of the linearized system:
\begin{equation}
    \nabla_{\mathbf{W}^{o(\nu)}} \mathbf{F} ~~ d\mathbf{W}^{o(\nu)} = - \mathbf{F}(\mathbf{W}^{o(\nu)}). 
    \label{eq:NR_wRPNN}
\end{equation}
Since the Jacobian is, in general,  expected to be ill-defined, the SVD decomposition can be used for computing the pseudo-inverse of the Jacobian, which is then used for obtaining the $\nu$-th update of $d\mathbf{W}^{o(\nu)}$ as:
\begin{equation}
    d\mathbf{W}^{o(\nu)} = - \left( \nabla_{\mathbf{W}^{o(\nu)}} \mathbf{F} \right)^\dagger \mathbf{F}(\mathbf{W}^{o(\nu)}), \qquad  \left( \nabla_{\mathbf{W}^{o(\nu)}} \mathbf{F} \right)^\dagger = \mathbf{U}_{p} \boldsymbol{\Sigma}_{p}^\dagger \mathbf{V}_{p}^\top,
    \label{eq:NR_SVD}
\end{equation}
where $\mathbf{U}_{p}\in \mathbb{R}^{ML \times ML}, \mathbf{V}_{p}^\top\in \mathbb{R}^{Mn_yn_{\epsilon} \times Mn_yn_{\epsilon}}$ are the unitary matrices containing the left and right singular vectors and  $\boldsymbol{\Sigma}_{p}\in \mathbb{R}^{ML \times Mn_yn_{\epsilon}}$ is the diagonal matrix with the singular values resulting from SVD.~For the implementation of the Newton-Raphson iterative scheme in Eq.~\eqref{eq:NR_SVD}, we used as stopping criterion $\lVert \mathbf{F}(\mathbf{W^{o(\nu+1)}}) \rVert_{l^2}<tol$, where $tol$ depends on the number $L$ of neurons in the hidden layer.

\section{The benchmark problems}
The efficiency of the proposed PIML approach is demonstrated via three benchmark models, namely the Michaelis-Menten (MM) $2$-dim. enzyme reaction mechanism \citep{michaelis1913kinetik}, the Target Mediated Drug Disposition (TMDD) $3$-dim. pharmacokinetic/pharmacodynamic mechanism \citep{levy1994pharmacologic,mager2001general} and the 3D Sel'kov model of glycolytic oscillations \citep{sel1968self,kourdis2013glycolysis}.~The MM and TMDD mechanisms have been extensively studied in the context of SPT and GSPT (see e.g. \citep{segel1989quasi,patsatzis2019new,schnell2000enzyme} and \citep{peletier2012dynamics,patsatzis2016asymptotic,kristiansen2019geometric}, respectively), thus allowing us to compare our framework with well documented SPT and GSPT analytical results.~The TMDD mechanism exhibits the interesting feature of two SIMs emerging in the phase space; the trajectories approach the first SIM, then exit from it through its boundaries, and finally approach the second SIM, through which they reach the stable equilibrium of the system.~We focus on the first SIM to demonstrate that our framework can approximate SIMs that do not include the stable equilibrium, that is the case of the MM mechanism.~Finally, the 3D Sel'kov model exhibits limit cycles in specific parameter regimes \citep{pye1966sustained}. In what follows, we describe the three benchmark problems and provide the analytical expressions of the corresponding SIMs, based on GSPT analytic calculations using the invariance equation and CSP with one iteration. As discussed, taking more iterations of the CSP, result in implicit forms of the SIMs (see for example in the Appendix \ref{app:MM_SIM_CSPexp}) which are not directly comparable with the explicit expression that we seek. 

\subsection{The Michaelis-Menten mechanism}
\label{sub:MM1}
The Michaelis-Menten (MM) reaction scheme describes the basic mechanism of enzyme action \citep{michaelis1913kinetik}, according to which, an enzyme $E$ reversibly binds to a substrate $S$ for the formation of a complex $C$ which is in turn decays irreversibly to form a product $P$ and the same enzyme $E$.~Using the law of mass action and the conservation laws for the enzyme and the substrate, the MM mechanism is formulated in the form of Eq.~\eqref{eq:gen} for the concentrations $s$ and $c$ as:
\begin{equation}
   \dfrac{d}{dt} \begin{bmatrix} s \\ c \end{bmatrix} = \begin{bmatrix} -k_{1f} (e_0-c) s + k_{1b} c~~~~~~~~~ \\ ~~~~~k_{1f} (e_0-c) s - k_{1b} c - k_2 c \end{bmatrix}, \quad s(0) = s_0, \quad c(0) = c_0<e_0,
   \label{eq:MM}
\end{equation}
where $k_{1f}$, $k_{1b}$ and $k_2$ are the formation, dissociation and catalysis rate constants, respectively and $e_0$ is the concentration of the enzyme.\par
\begin{figure}[!h]
    \centering
    \includegraphics[scale=0.3]{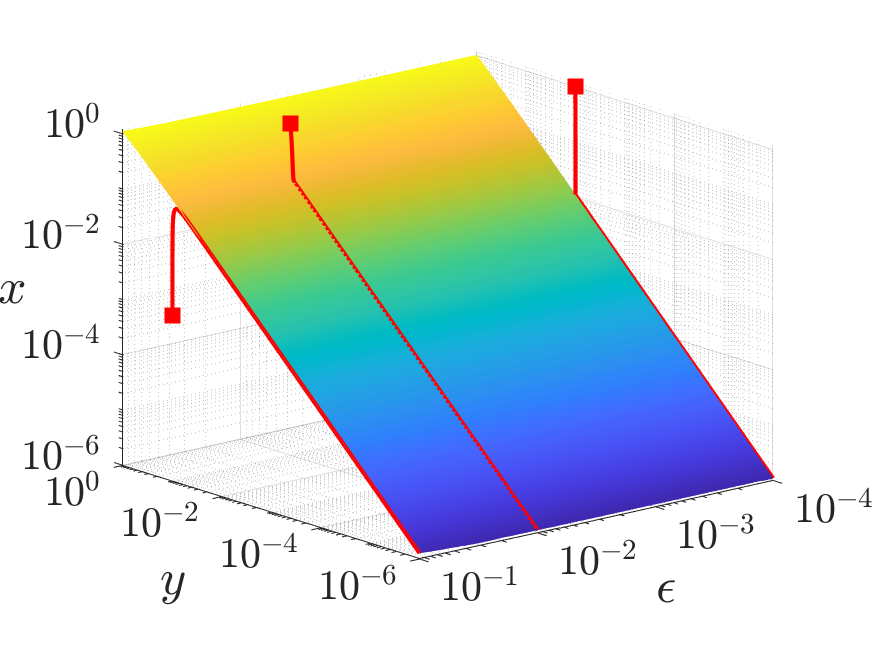}
    \caption{MM system (\ref{eq:MMss}). SIM surface of the slow system of Eq.~\eqref{eq:MMss} in the domain $y\in [10^{-6}, 1]$ for various values of $\epsilon \in [10^{-4}, 10^{-1}]$.~The red trajectories (starting at red squares) are attracted to the SIM and then evolve on it.}
    \label{fig:MM_SIM}
\end{figure}

The evolution of the MM mechanism exhibits fast/slow timescale separation almost everywhere in the parameter space \citep{patsatzis2019new}, thus admitting a $N-M=1$-dim. SIM governing its slow evolution. Various fast-slow subsystems in the form of Eqs.~(\ref{eq:SPfast},\ref{eq:SPslow}) have been extracted in the literature, depending on the regions of the parameter and phase space in which the identification of the SIM is sought; e.g., see \citep{segel1989quasi,patsatzis2019new,schnell2000enzyme,patsatzis2023algorithmic}.~For our illustrations, we consider the form proposed by Segel and Slemrod \citep{segel1989quasi} that is appropriate for investigating the cases where the complex $c$ is the fast variable. The introduction of the rescaled variables $x = c (\kappa +1)/e_0$ and $y=s/s_0$ in the rescaled time $\tau=t k_{1f} e_0 \sigma/(\kappa + 1)$ casts the system in Eq.~\eqref{eq:MM} to its slow subsystem form in Eq.~\eqref{eq:SPslow}, as:
\begin{equation}
    \epsilon \dfrac{dx}{d\tau} = y - \dfrac{\kappa + y}{\kappa + 1} x, \qquad \dfrac{dy}{d \tau} = \sigma^{-1} (-(\kappa +1) y + (\kappa-\sigma + y)x)
    \label{eq:MMss}
\end{equation}
where $\kappa = K_M/s_0$, $\sigma = K/s_0$ and $\epsilon = (\sigma e_0/s_0)/(\kappa +1)^2\ll1$; $K_M = (k_{1b}+k_2)/k_{1f}$ is the Michaelis-Menten constant and $K=k_2/k_{1f}$ the Van Slyke-Culen constant.~A schematic representation of the SIM arising in the phase space $(x,y)$ is shown in Fig.~\ref{fig:MM_SIM} for various values of $\epsilon$, with the parameter values set to $\kappa =10$, $\sigma =100$.~As shown, the trajectories are attracted to the SIM and then evolve on it, towards reaching the fixed point $(0,0)$ of the MM mechanism.

On the basis of SPT, Segel and Slemrod \citep{segel1989quasi} derived the SPT $\mathcal{O}(\epsilon)$ regular asymptotic expansion of the SIM (see also in \cite{kaper2002asymptotic}):
\begin{equation}
    x = h_0(y) + \epsilon h_1(y) 
    = \dfrac{\kappa+1}{\kappa+y} y + \epsilon \dfrac{\kappa (\kappa+1)^2}{\sigma (\kappa+y)^3} \Bigg[ \dfrac{2 \sigma y}{\kappa + y} - y +\dfrac{(\kappa-\sigma)y}{\kappa} ln\left( \dfrac{\kappa+y}{(\kappa+1)y}\right)\Bigg] 
    \label{eq:MManalSIM}
\end{equation}
For obtaining SIM approximations on the basis of GSPT, we additionally employed (a) the basic analytic technique, using the invariance equation~\eqref{eq:Inv} and (b) the CSP methodology with one iteration \citep{lam1989understanding,goussis1992study,valorani2005higher}, resulting in the following lemma (for the proof, see Appendix \ref{app:MM_SIM_GSPTexp} and \ref{app:MM_SIM_CSPexp}): 
\begin{lemma}
Consider the MM slow subsystem in Eq.~\eqref{eq:MMss} with $\epsilon \ll 1$.~The analytic SIM approximation, derived on the basis of the invariance equation~\eqref{eq:Inv}, is given by the $\mathcal{O}(\epsilon^2)$ regular asymptotic expansion:
\begin{equation}
    x = h_0(y)+\epsilon h_1(y) + \epsilon^2 h_2(y) = \dfrac{\kappa+1}{\kappa+y} y + \epsilon \dfrac{\kappa (\kappa+1)^3 y}{(\kappa+y)^4} - \epsilon^2 \dfrac{\kappa (\kappa+1)^5 y(\kappa^2+3\sigma y+\kappa (y-2\sigma))}{\sigma(\kappa+y)^7}
    \label{eq:MM_GSPT}
\end{equation}
In addition, the explicit analytic SIM approximation, derived on the basis of CSP with one iteration, is:
\begin{equation}
x=h(y,\epsilon)=\dfrac{\sigma(\kappa+y)^2+\epsilon(\kappa+1)^2(\kappa-\sigma+2y)-\sqrt{\left(\epsilon(\kappa+1)^2(\kappa-\sigma)+\sigma(\kappa+y)^2 \right)^2+4\epsilon(\kappa+1)^2\sigma^2y(\kappa+y)}}{2 \epsilon (\kappa+1)(\kappa-\sigma+y)}
    \label{eq:MM_CSP}
\end{equation}
\end{lemma}
The analytical SIM approximations in Eqs.~(\ref{eq:MManalSIM}-\ref{eq:MM_CSP}) recover, to the leading order $x=h_0(y)$, the QSSA for the complex $c$, well-known as \emph{standard} QSSA (sQSSA).~All the above SPT/GSPT approximations are accurate within the regions of the parameter and phase space where $\epsilon \ll 1$ and $\epsilon (\kappa+1)(\kappa+1-x) \ll \sigma (\kappa + y)$; introduced as $e_0-c \ll K_M+s$ for the MM original system in \citep{patsatzis2019new,patsatzis2023algorithmic}. 
Finally, note that two iterations of the CSP result in an implicit form for the SIM (see in the Appendix \ref{app:MM_SIM_CSPexp}), that is not directly comparable to the explicit expression of the SIM obtained by the proposed PIML approach.

\subsection{The Target Mediated Drug Disposition mechanism}
\label{sub:TMDD1}
The Target Mediated Drug Disposition (TMDD) mechanism is a pharmacokinetic/pharmacodynamic (PKPD) reaction scheme that describes the action of specific drugs, such as monoclonal antibodies, the disposition of which affects the pharmacodynamic properties of their pharmacological targets \citep{levy1994pharmacologic,mager2001general}.~The simplest form of such interaction is the one-compartmental TMDD mechanism, 
which can be formulated by a $N=3$-dim. system of nonlinear ODEs, describing the evolution of a ligand $L$, a receptor $R$ and a ligand-target complex $RL$ concentrations as:
\begin{equation}
    \dfrac{d}{dt}\begin{bmatrix} L~~ \\ R~~ \\ RL
    \end{bmatrix} = \begin{bmatrix} -k_{on}L.R+k_{off}RL-k_{el} L \qquad \quad ~~~\\ 
    -k_{on}L.R+k_{off}RL+k_{syn}-k_{deg}R \\
    k_{on}L.R-k_{off}RL - k_{int}RL\quad ~~~~
    \end{bmatrix}, \quad L(0)=L_0, \quad R(0)=R_0 \quad RL(0)=0.
    \label{eq:TMDD}
\end{equation}
$k_{on}$ and $k_{off}$ are the formation and dissociation rate constants, $k_{el}$ is the ligand elimination one, $k_{syn}$ and $k_{deg}$ are the receptor synthesis and degradation ones, and $k_{int}$ is the internalization one.~Following \citep{peletier2012dynamics}, we set the parameter values to $k_{on}=0.091$, $k_{off}=0.001$, $k_{el}=0.0015$, $k_{syn}=0.11$, $k_{deg}=0.0089$ and $k_{int}=0.003$ and the initial condition of the receptor to $R_0=k_{syn}/k_{deg}=12$ (for simulating the administration of an intravenous injection of the drug $L_0$ when the system is in equilibrium).

\begin{figure}[!h]
    \centering
    \subfigure[$\epsilon=10^{-1}$]{
    \includegraphics[width=0.3\textwidth]{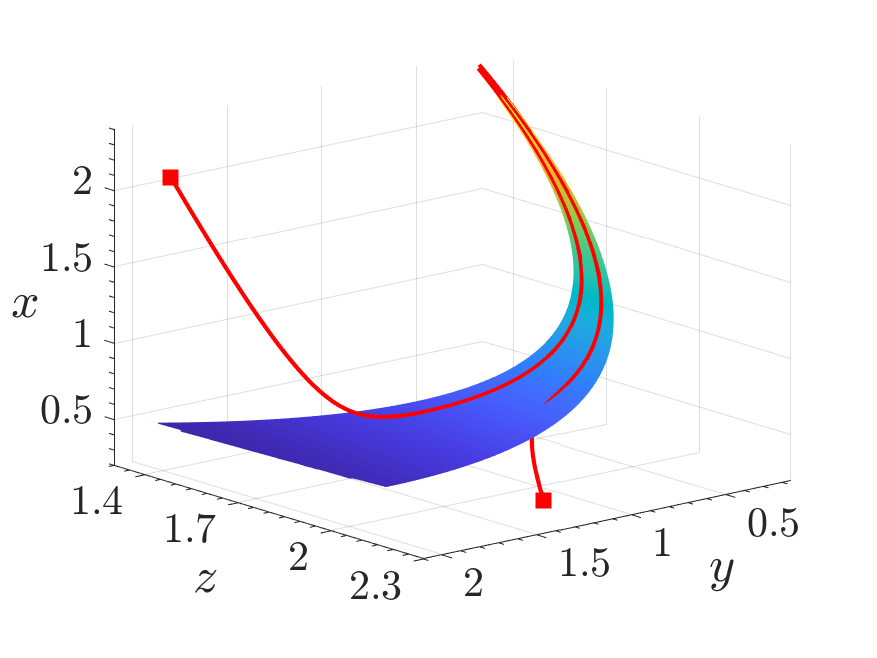}}
    \subfigure[$\epsilon=10^{-2}$]{
    \includegraphics[width=0.3\textwidth]{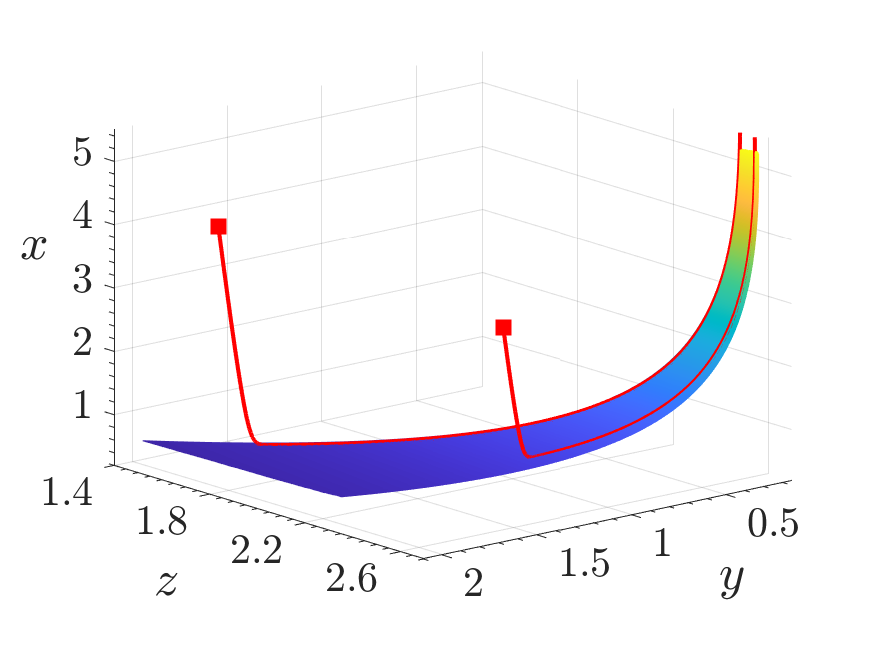}}
    \subfigure[$\epsilon=10^{-4}$]{
    \includegraphics[width=0.3\textwidth]{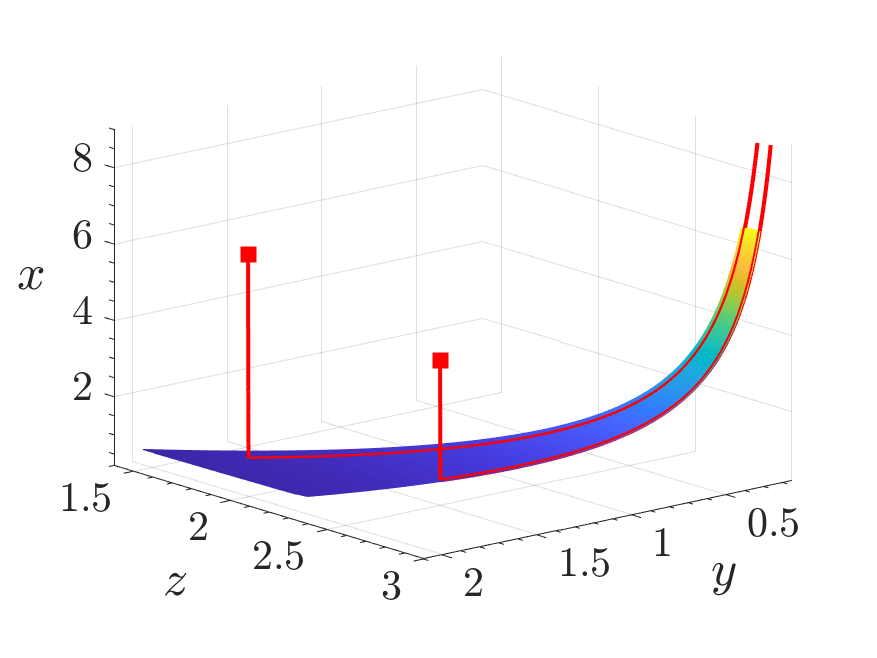}}
    \caption{TMDD system (\ref{eq:TMDD_SPA}).~SIM surface of the slow system of Eq.~\eqref{eq:TMDD_SPA} in the domain $(y,z)\in [0.2, 2.0]\times[1.3,2.9]$ for three indicative values of  $\epsilon$.~The red trajectories (starting at red squares) are attracted to the SIM and then evolve on it until eventually exiting from its boundaries.}
    \label{fig:TMDD_SIM}
\end{figure}

The dynamics of the TMDD model in Eq.~\eqref{eq:TMDD} exhibits fast/slow timescale separation in different regions of the phase space and various SIM approximations have been proposed by employing either the QSSA for $L$ \citep{peletier2012dynamics}, $R$ \citep{aston2011mathematical} or $RL$ \citep{van2016topics}, or the \emph{partial equilibrium approximation} (PEA) \citep{mager2001general,peletier2009dynamics}.
~Recent systematic analysis in the context of GSPT demonstrated that the TMDD model initially evolves along a $N-M=2$-dim. SIM, then degenerates from it and subsequently approaches another SIM that leads the system to its equilibrium \citep{patsatzis2016asymptotic}.~Here, we focus on the former period, during which 
the adoption of the rescaled variables $x=R~(k_{on}L_0)/k_{syn}$, $y=L/L_0$ and $z=RL~k_{deg}/k_{syn}$ in the rescaled time $\tau=k_{int} t$ casts the system in Eq.~\eqref{eq:TMDD} to its slow subsystem form in Eq.~\eqref{eq:SPslow}, as: 
\begin{equation}
    \epsilon \dfrac{dx}{d\tau} = -xy+k_1z+1-\epsilon k_2 x, \quad
    \dfrac{dy}{d\tau}=k_3(-xy+k_1z)-k_4y, \quad
    \dfrac{dz}{d \tau}=k_2(xy-k_1z)-z,
    \label{eq:TMDD_SPA}
\end{equation}
where $k_1=k_{off}/k_{deg}$, $k_2=k_{deg}/k_{int}$, $k_3=k_{syn}/(k_{int}L_0)$, $k_4=k_{el}/k_{int}$ and $\epsilon = k_{int}/(k_{on}L_0)\ll 1$.~Note that the above TMDD slow subsystem accurately reflects the dynamics of the TMDD original one in Eq.~\eqref{eq:TMDD} during the ``slow first-order disposition" period, when $\epsilon/y \ll 1$ \citep{patsatzis2016asymptotic}.~Considering the parameter set proposed in \citep{peletier2012dynamics}, a schematic representation of the SIM is shown in Fig.~\ref{fig:TMDD_SIM} for 3 indicative values of $\epsilon$.~As shown, the trajectories are attracted towards the SIM and then evolve on it until $\epsilon/y \ll 1$, when the trajectories degenerate from this SIM due to the invalidity of the slow subsystem in Eq.~\eqref{eq:TMDD_SPA}.~It is further shown in Fig.~\ref{fig:TMDD_SIM} that for different values of $\epsilon$, the SIM surface changes its orientation and curvature, as expected.

The approximation of the SIM provided by the proposed PIML scheme is compared with the analytic GSPT approximations for the TMDD slow subsystem in Eq.~\eqref{eq:TMDD_SPA}. Due to the complexity of the system, an SPT $\mathcal{O}(\epsilon)$ SIM approximation is not tractable.~Thus, we derived SIM approximations by employing: (a) the basic analytic technique using the invariance equation \eqref{eq:Inv}, and, (b) CSP methodology with one iteration; for the proof of the following lemma, see Appendix \ref{app:TMDD_SIM_GSPTexp} and \ref{app:TMDD_SIM_CSPexp}.
\begin{lemma}
Consider the TMDD slow subsystem in Eq.~\eqref{eq:TMDD_SPA} with $\epsilon \ll 1$.~The analytic SIM approximation, derived on the basis of the invariance equation~\eqref{eq:Inv}, is given by the $\mathcal{O}(\epsilon^2)$ regular asymptotic expansion:
\begin{align}
    & x = h_0(y,z)+\epsilon h_1(y,z) + \epsilon^2 h_2(y,z) = \dfrac{1+k_1z}{y} - \epsilon \dfrac{k_3(1+k_1)z+y(k_2+k_1k_2+k_4+k_1z(-1+k_2+k_4))}{y^3} + \nonumber \\
    &  \dfrac{\epsilon^2}{y^5} \left( k_3^2 (k_1 z+1) (k_1 z + 4) + y^2 ((k_2 + k_4) (k_2 + 2 k_4) + k_1 k_2 (3 k_2 + 4 k_4-1) + k_1 (k_2 + k_4-1) (k_2 + 2 k_4-1) z+ \right.  \nonumber \\
    & \left. k_1^2 k_2 (k_2 + ( k_2 + k_4-1) z)) + k_3 y (6 k_4 + k_1 z (7 k_4 + k_1 (-1 + k_4) z-4) + k_2 (4 + k_1 (5 + z (5 + k_1 (2 + z))))) \right)
    \label{eq:TMDD_GSPT}
\end{align}
In addition, the explicit analytic SIM approximation, derived on the basis of CSP with one iteration, is:
\begin{equation}
x = h(y,z,\epsilon) = -\dfrac{(y+\epsilon k_2)^2+\epsilon(y(k_1k_2+k_4)-k_1k_3z)}{2\epsilon k_3 y} \left(1-\sqrt{1+\dfrac{4 \epsilon k_3 y (\epsilon k_2 + y + k_1z (\epsilon  (1 + k_2 + k_1 k_2) + y))}{((y+\epsilon k_2)^2+\epsilon(y(k_1k_2+k_4)-k_1k_3z))^2}} \right)
\label{eq:TMDD_CSP}
\end{equation}
\end{lemma}
The SIM approximations in Eqs.~(\ref{eq:TMDD_GSPT}, \ref{eq:TMDD_CSP}) recover, to the leading order, the same QSSA approximation for the receptor $x$, $x=h_0(y,z)$.~Note that the SIM approximation provided by the CSP is accurate within the region of validity of the slow subsystem in Eq.~\eqref{eq:TMDD_GSPT} (i.e., when $\epsilon/y \ll 1$), while the ones provided by the QSSA, $\mathcal{O}(\epsilon)$ and $\mathcal{O}(\epsilon^2)$ regular expansions are accurate when the conditions $k_2 \epsilon \ll y$ and $k_3 \epsilon x \ll y$ are additionally satisfied \citep{patsatzis2016asymptotic}.

\subsection{The 3D Sel'kov model of glycolytic oscillations}
\label{sub:TLC1}
\begin{figure}[!h]
    \centering
    \subfigure[$\epsilon=10^{-1}$]{
    \includegraphics[width=0.3\textwidth]{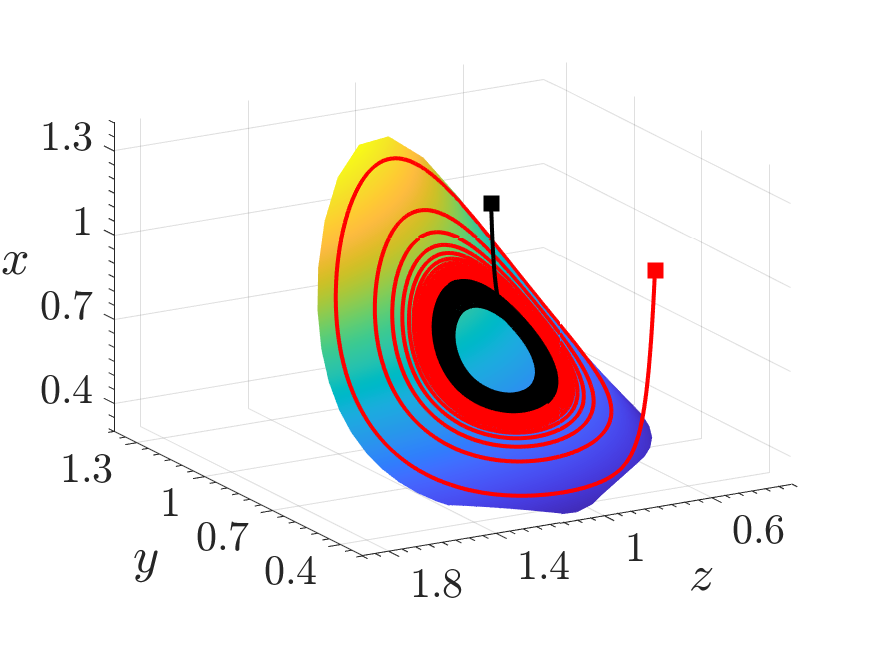}}
    \subfigure[$\epsilon=10^{-2}$]{
    \includegraphics[width=0.3\textwidth]{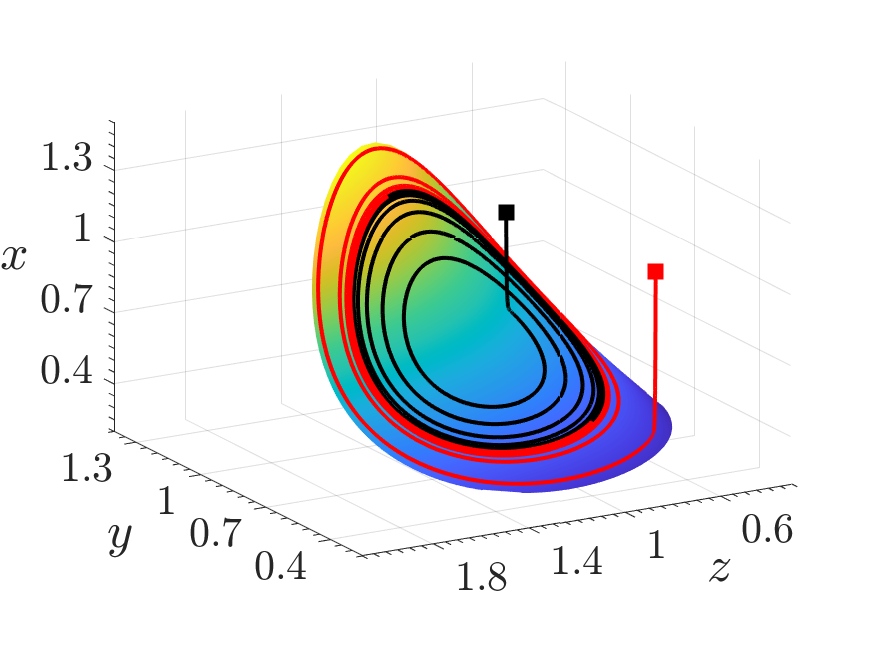}}
    \subfigure[$\epsilon=10^{-4}$]{
    \includegraphics[width=0.3\textwidth]{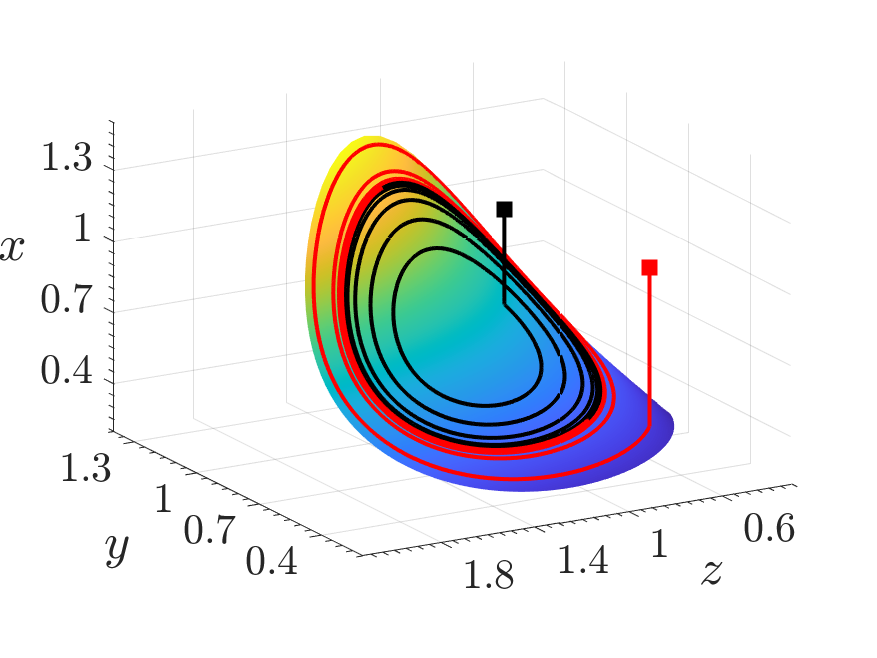}}
    \caption{3D Sel'kov system~\eqref{eq:ToyLC_SPA}.~SIM surface of arising in the phase space along the domain of slow variables $(y,z)\in [0.2, 1.4]\times[0.3,2.1]$ for three indicative values of  $\epsilon$.~The red/black trajectories (starting at red/black squares) are attracted to the SIM and then evolve on it, in the exterior/interior of the limit cycle, towards reaching it.}
    \label{fig:TLC_SIM}
\end{figure}

Here, we consider an extension of the Sel'kov kinetics model of glycolysis, which under specific parameter regimes may exhibit limit cycles, relaxation-oscillations or fixed points \citep{roy2011periodic,sel1968self}.
In particular, we include the addition of a fast variable \citep{kourdis2013glycolysis}, so that the resulting 3D model incorporates, in the slow timescale, the dynamics of the 2D model.~The resulting 3D Sel'kov model in non-dimensional form reads:  
\begin{equation}
    \epsilon \dfrac{dx}{d\tau} = y^2z - kxy, \quad
    \dfrac{dy}{d\tau}=a z + y^2z -y +\epsilon x, \quad
    \dfrac{dz}{d \tau}=-a z - y^2z +b,
    \label{eq:ToyLC_SPA}
\end{equation}
where $a$, $b$ and $k$ are parameters of the original Sel'kov model and $\epsilon\ll 1$.~Being interested in limit cycles, we choose the parameters in the stable limit cycle regime \citep{roy2011periodic}; $a=0.1$, $b=0.6$ and $k=1$.~As evident, the model in Eq.~\eqref{eq:ToyLC_SPA} is already written in the slow subsystem form of Eq.~\eqref{eq:SPslow}; hence, its slow dynamics evolves on a $N-M=2$-dim. SIM with $x$ being the fast variable.~In particular, for the selected parameter values, the $N-M=2$-dim. SIM includes and converges to the stable limit cycle.~A representation of the SIM emerging in the phase space $(x,y,z)$ is shown in Fig.~\ref{fig:TLC_SIM} for three indicative values of $\epsilon$, where trajectories are shown to be attracted to the SIM and then evolving on it towards reaching the limit cycle either from its exterior or interior.~As $\epsilon$ varies, the limit cycle changes, so as the SIM surface around it.

For deriving explicit analytic SIM approximations on the basis of GSPT, we employed: (a) the basic analytic technique, using the invariance equation~\eqref{eq:Inv}, and, (b) CSP with one iteration, resulting in the following Lemma (for the proof, see in Appendix \ref{app:TLC_SIM_GSPTexp} and \ref{app:TLC_SIM_CSPexp}).
\begin{lemma}
Consider the slow subsystem of the 3D Sel'kov model in Eq.~\eqref{eq:ToyLC_SPA} with $\epsilon \ll 1$.~The analytic SIM approximation, derived on the basis of the invariance equation~\eqref{eq:Inv}, is given by the $\mathcal{O}(\epsilon^2)$ regular asymptotic expansion:
\begin{align}
    & x = h_0(y,z)+\epsilon h_1(y,z) +\epsilon^2 h_2(y,z) = \dfrac{y z}{k} + \epsilon
    \left( \dfrac{z-b}{k^2} + \dfrac{z(a+y^2)(y-z)}{k^2y} \right) + \dfrac{\epsilon^2}{k^3 y^3} \left(b y (-y (1 + a + y^2) + \right.  \nonumber \\
     & \left. 2 (a + y^2) z) + z (a y (y + 2 y^3 + z - 6 y^2 z) + a^2 (y^2 - 2 y z - z^2) + y^3 (-2 z + y (3 + y^2 - 4 y z + z^2))) \right)
    \label{eq:TLC_GSPT}
\end{align}
In addition, the explicit analytic SIM approximation, derived on the basis of CSP with one iteration, is:
\begin{equation}
x = h(y,z,\epsilon) = \left( \dfrac{yz}{k}+\dfrac{y-z(a+y^2)}{2\epsilon}+\dfrac{ky^2}{2\epsilon^2}\right) \left(1 - \sqrt{ 1 + \dfrac{4 \epsilon^2 k y (\epsilon b y +z (2 \epsilon z (a + y^2) -y (k y + \epsilon (2 + a + y^2))))}{(k y (k y+\epsilon) + \epsilon z(2 \epsilon y - k (a + y^2)))^2}} \right)
\label{eq:TLC_CSP}
\end{equation}
\end{lemma}
The leading order term of the regular asymptotic expansion in Eq.~\eqref{eq:TLC_GSPT} $x=h_0(y,z)$ corresponds to the QSSA expression for the fast scaled variable $x$.

\section{Numerical results}
\label{sec:NumR}
We assessed the efficiency of the proposed PIML scheme, thus providing a comparison between symbolic (SD), automatic (AD), and, numerical differentiation using forward finite differences (FD), with respect to the computational cost. Furthermore, we provide a comparative analysis with
the analytically derived GSPT approximations given in the previous section.~The numerical accuracy of all the explicit SIM approximations (either PIML or GPST-derived ones) were assessed on the basis of the actual SIMs, along which the numerical solution of the singularly perturbed system evolves, computed by numerical integration of the ODEs.

To learn the approximations of the SIMs via the proposed PIML schemes, we have  collected data (serving as collocation points for the solution of the IE) in the domain of $\mathbf{y} \in \Omega \subset \mathbb{R}^{N-M}$ and $\epsilon \in I \subset \mathbb{R}$ ($I=[10^{-4}, 10^{-1}]$).~Since there is no guarantee that the SIM exists for every grid point in the $\Omega\times I$ subdomain, we collected these collocation points from numerically derived trajectories, within the regions of validity of the slow subsystems in the form of Eq.~\eqref{eq:SPslow}.~In particular, for all the problems considered, we varied $\epsilon$ in the interval $I$ (logarithmically equally spaced $\epsilon_j\in I$ for $j=1,\ldots,n_\epsilon $).~For every $\epsilon_j$, we generated  trajectories that sufficiently cover the domain $\Omega\times I$ by using a number of random initial conditions close to the boundary $\partial \Omega$.~From the resulting trajectories, we collected only the values of the slow variables to get $20$ equidistant in time points per trajectory. These sets form the collocation points $[\mathbf{y}_i,\epsilon_j]^\top\in  \mathbb{R}^{D\times n_y n_\epsilon}$ on which the IE is solved via PIML. 

The PIML schemes were trained using an 80\% uniformly random sample of the above data sets, while the rest 20\% of the points were used for validation purposes.~For all problems under study, the number of neurons in the hidden layer of the SLFNNs was set to $L = 20$.~To enable straightforward comparison, for the RPNN we selected the same number of learnable parameters as in the SLFNN, resulting to $L=81$ for the MM mechanism and $L=101$ for the 3D Sel'kov system.~For the TMDD mechanism, to achieve a high approximation accuracy, $L=400$ neurons were required; the significant change of the SIM orientation and curvature for different values of $\epsilon$ requires denser sampling of the basis functions.~The logistic sigmoid function was used as activation function for both SLFNNs and RPNNs.~The tolerance $tol = 10^{-3}$ was set as stopping criterion for both the Levenberg-Marquardt and Newton-Raphson iterative schemes for the PIML optimization problem with SLFNNs and RPNNs, respectively.

To assess the numerical accuracy of the explicit SIM approximations, we constructed test sets consisting of data \textit{lying exclusively on the SIM} at the domain of interest $\Omega \times I$. To achieve this, we integrated the slow subsystem given by Eq.~\eqref{eq:SPslow}, within the regimes of their validity, and kept data from the trajectories only after a transient period, set as $t>10 \epsilon$; $\epsilon$ was randomly sampled from a uniform distribution in the interval $I$. Again, for every $\epsilon$, we considered a number of random initial conditions outside $\Omega$, ensuring, as explained above, that each trajectory evolves on the SIM in the domain of interest $\Omega\times I$. From the resulting time series, we collected $100$ equidistant- in time- points per trajectory to form the test sets consisting of $[\mathbf{y}_i,\epsilon_j]^\top\in  \mathbb{R}^{D\times n n_\epsilon}$, where $n$ is the total number of points per $\epsilon_j$ and the corresponding values of the fast variables $\mathbf{x}_{i,j}\in \mathbb{R}^{M\times n n_\epsilon}$ as obtained by the numerical integration of the ODEs. The numerical approximation accuracy of the PIML and GSPT schemes is measured in terms of the $l^2$, $l^\infty$ and MSE of $\lVert \mathbf{x}_{i,j}-\mathcal{N}(\mathbf{y}_i,\epsilon_j)\rVert$, and $\lVert \mathbf{x}_{i,j}-\mathbf{h}(\mathbf{y}_i,\epsilon_j)\rVert$ ($i=1,\ldots,n$ and $j=1,\ldots,n_\epsilon$), respectively.

All simulations were carried out with a CPU Intel(R) Xeon(R) CPU E5-2630 v4 @ 2.20GHz (2 processors), RAM 64.0 GB using MATLAB R2022b.

\subsection{The Michaelis-Menten mechanism}

Here, we computed the SIM in the domain $y\in [10^{-6}, 1]$ for $\epsilon \in [10^{-4}, 10^{-1}]$.~In particular, we considered $n_\epsilon =13$ logarithmically spaced values of $\epsilon_j \in [10^{-4}, 10^{-1}]$.~For each $\epsilon_j$, $10$ trajectories were generated with random initial conditions varying in $x(0) \in [0,2]$ and $y(0) \in [1,2]$.~From the resulting trajectories, $n_y=200$ points were sampled ($20$ equidistant in time points per trajectory; integrations stopped for $y<10^{-6}$) and the values of the slow variables $y_i$ for $i=1,\ldots,n_y$ were collected. Hence, for the training set: $[y_i,\epsilon_j]^\top \in \mathbb{R}^{2\times2600}$. 

Table~\ref{tab:MM_train} summarizes the comparison results between the differentiation schemes (SD, AD, FD). In particular, we report the loss function $\lVert \boldsymbol{\mathcal{F}} \rVert^2_2$  for the training and validation sets, as well as the corresponding computational costs, on the basis of 10 runs with different randomly sampled training and validation sets.~As shown for SLFNNs, SD results in smaller loss functions and in $\sim 3\times$ faster times than FD, which is in turn $\sim 3\times$ faster than AD. As it is shown, the RPNNs result in similar loss functions while they are more $\sim100\times$ faster than the SLFNNs.
\begin{table}[!h]
    \centering
    \footnotesize
    \begin{tabular}{l| c c | c c c}
    \toprule
    & \multicolumn{2}{c}{Loss Function $\lVert \boldsymbol{\mathcal{F}} \rVert^2_2$}  & \multicolumn{3}{|c}{Computational times (s)} \\
    PIML scheme & Training & Validation & mean & min & max \\
    \midrule
    SLFNN AD &	3.01E$-$06	&	7.65E$-$07	&	6.52E$+$01	&	3.87E$+$01	&	7.51E$+$01	\\  
    SLFNN FD &	3.90E$-$06	&	1.07E$-$06	&	3.08E$+$01	&	1.91E$+$01	&	4.43E$+$01	\\
    SLFNN SD &	6.18E$-$08	&	2.62E$-$08	&	1.05E$+$01	&	4.72E$+$00	&	1.25E$+$01	\\
    RPNN SD &	6.79E$-$08	&	4.86E$-$08	&	1.06E$-$01	&	7.66E$-$02	&	4.20E$-$01	\\
    \bottomrule
    \end{tabular}
    \caption{MM system (\ref{eq:MMss}). Loss function $\lVert \boldsymbol{\mathcal{F}} \rVert^2_2$ of the PIML schemes for the training and validation sets using Automatic Differentiation (AD), Finite Differences (FD) and Symbolic Differentiation (SD). The corresponding  computational times (in seconds) are also given. The results are obtained by averaging over 10 runs.}
    \label{tab:MM_train}
\end{table}

To build the test set ($[y_i,\epsilon_j]^\top$ and the corresponding values of the fast variables $x_{i,j}$), we considered $500$ trajectories of the MM slow subsystem with $50$ randomly varied values of $\epsilon$, initialized with $10$ random initial conditions each, in $x(0) \in [0,2]$ and $y(0) \in [2,3]$.~We then kept the data after $t=10\epsilon$ in order for the trajectory to lie on the SIM and recorded $100$ equidistant-in time-points per trajectory (integrations stopped for $y<10^{-6}$), including only data in the desired domain $y\in [10^{-6}, 1]$.

Table~\ref{tab:MM_acc} enlists the \emph{overall, with respect to all values of $\epsilon$}, $l^2$, $l^{\infty}$ and MSE $\lVert x_{i,j} - \mathcal{N}(y_i, \epsilon_j) \rVert$ approximation errors as obtained by the PIML schemes, and the corresponding $\lVert \mathbf{x}_{i,j}-\mathbf{h}(\mathbf{y}_i,\epsilon_j)\rVert$ approximation errors as obtained by the sQSSA, SPT, GSPT and CSP-derived SIM approximations in Eqs.~(\ref{eq:MManalSIM})-(\ref{eq:MM_CSP}).~As shown, the PIML schemes provide SIM approximations with a high accuracy (the overall $l^{\infty}$ is of the order of $1E-04$). Furthermore, the proposed PIML schemes are more accurate with respect to all norms, than the sQSSA, SPT $\mathcal{O}(\epsilon)$, GSPT $\mathcal{O}(\epsilon)$, GSPT  $\mathcal{O}(\epsilon^2)$, and CSP with one iteration, approximations of the SIM.
\begin{table}[!h]
    \centering
    \resizebox{\textwidth}{!}{\begin{tabular}{l| c c c c | c c c c c}
    \toprule
    & \multicolumn{4}{|c}{PIML} &  \multicolumn{5}{| c}{analytic SPT/GSPT approximations} \\
    \midrule
    Error  &  SLFNN AD &  SLFNN FD  & SLFNN SD  & RPNN SD &   sQSSA &  SPT $\mathcal{O}(\epsilon)$ & GSPT $\mathcal{O}(\epsilon)$ & GSPT $\mathcal{O}(\epsilon^2)$ & CSP \\
    \midrule 
    $l^2$ &	9.72E$-$03	&	1.08E$-$02	&	1.89E$-$03	&	1.50E$-$03	&	1.29E$+$00	&	1.14E$+$00	&	1.44E$-$01	&	2.10E$-$02	&	8.00E$-$02	\\
    $l^{\infty}$    &	2.97E$-$04	&	3.30E$-$04	&	1.19E$-$04	&	6.68E$-$05	&	8.90E$-$02	&	7.56E$-$02	&	1.09E$-$02	&	1.51E$-$03	&	5.78E$-$03	\\
    MSE &	2.10E$-$09	&	2.41E$-$09	&	8.52E$-$11	&	4.49E$-$11	&	3.34E$-$05	&	2.62E$-$05	&	4.15E$-$07	&	8.79E$-$09	&	1.28E$-$07	\\
    \bottomrule
    \end{tabular}}
    \caption{MM system (\ref{eq:MMss}). SIM approximation accuracy \emph{over all simulations for all values of $\epsilon$},  $l^2$, $l^{\infty}$ and MSE approximation errors resulting from the PIML schemes (SLFNN and RPNN) and the analytic SPT/GSPT approximations: sQSSA, SPT, GSPT $\mathcal{O}(\epsilon)$ and $\mathcal{O}(\epsilon^2)$ asymptotic expansions, and the CSP with one iteration in Eqs.~(\ref{eq:MManalSIM}-\ref{eq:MM_CSP}), respectively.~The numerical accuracy of each SIM approximation is compared with the numerical solution $x_{i,j}$ of the MM slow subsystem in Eq.~\eqref{eq:MMss}.}
    \label{tab:MM_acc}
\end{table}

\begin{figure}[!h]
    \centering
    \subfigure[PIML SLFNN]{
    \includegraphics[width=0.3\textwidth]{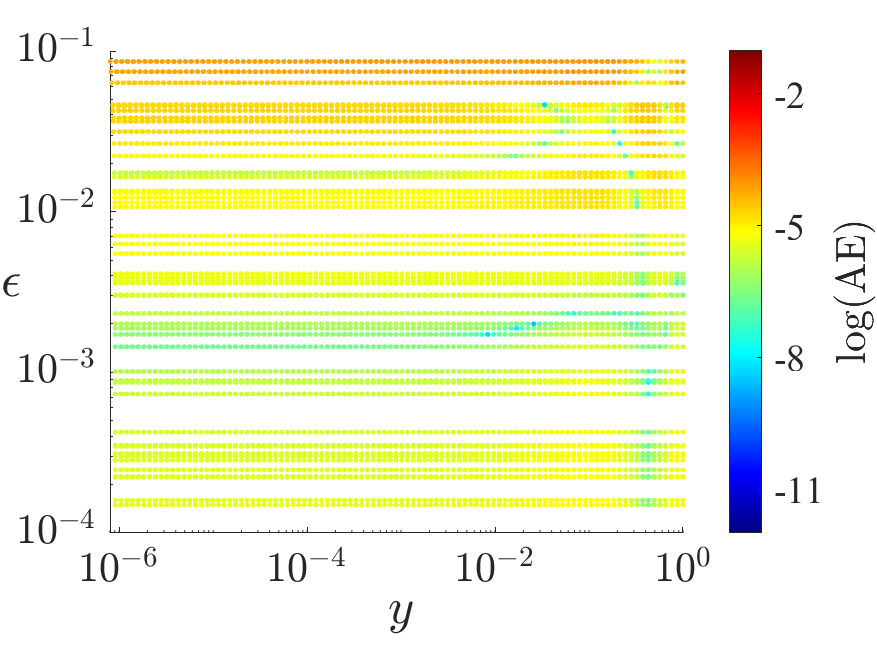}}
    \subfigure[PIML RPNN]{
    \includegraphics[width=0.3\textwidth]{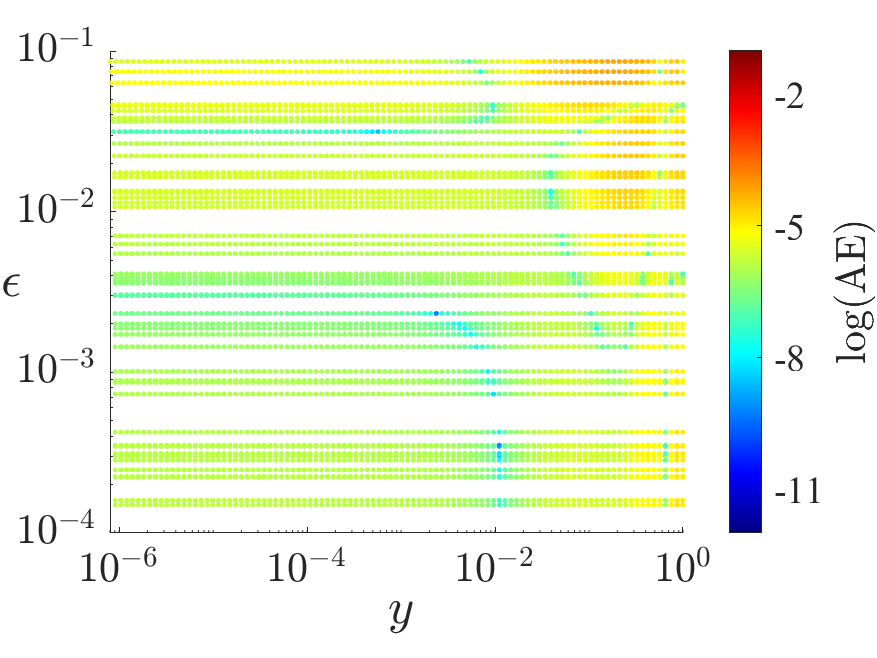}}
    \subfigure[sQSSA]{
    \includegraphics[width=0.3\textwidth]{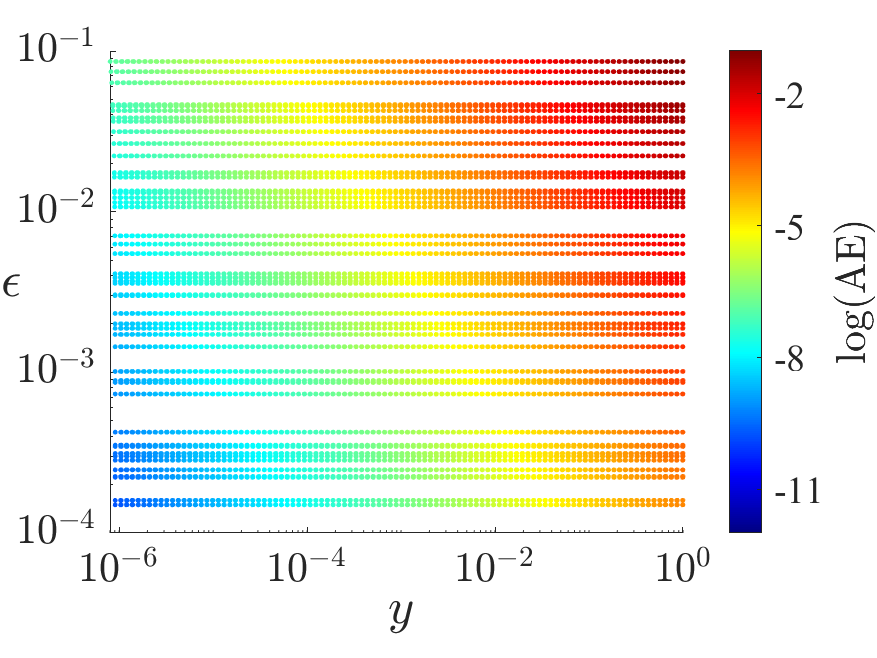}}
    \subfigure[GSPT $\mathcal{O}(\epsilon)$]{
    \includegraphics[width=0.3\textwidth]{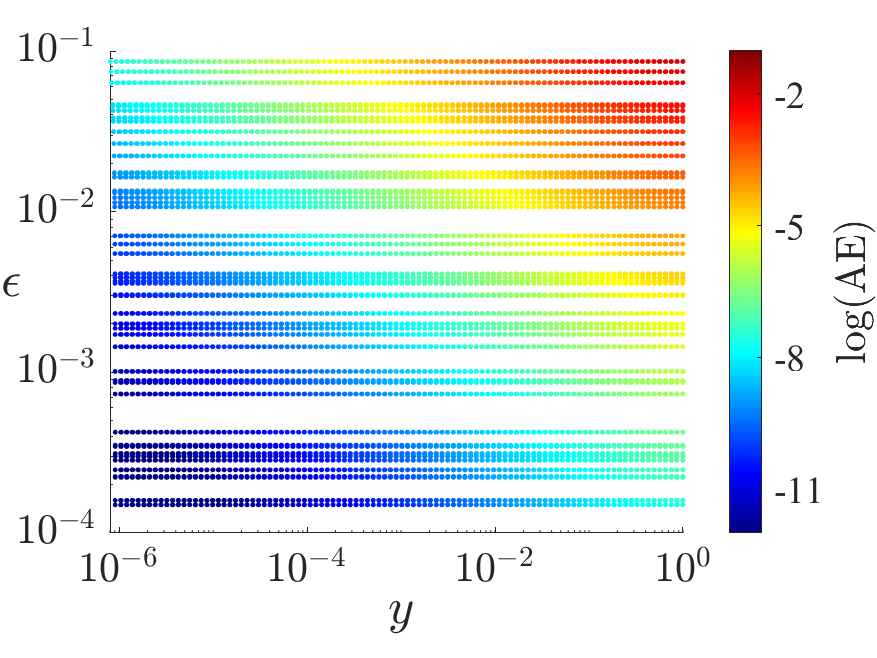}}
    \subfigure[GSPT $\mathcal{O}(\epsilon^2)$]{
    \includegraphics[width=0.3\textwidth]{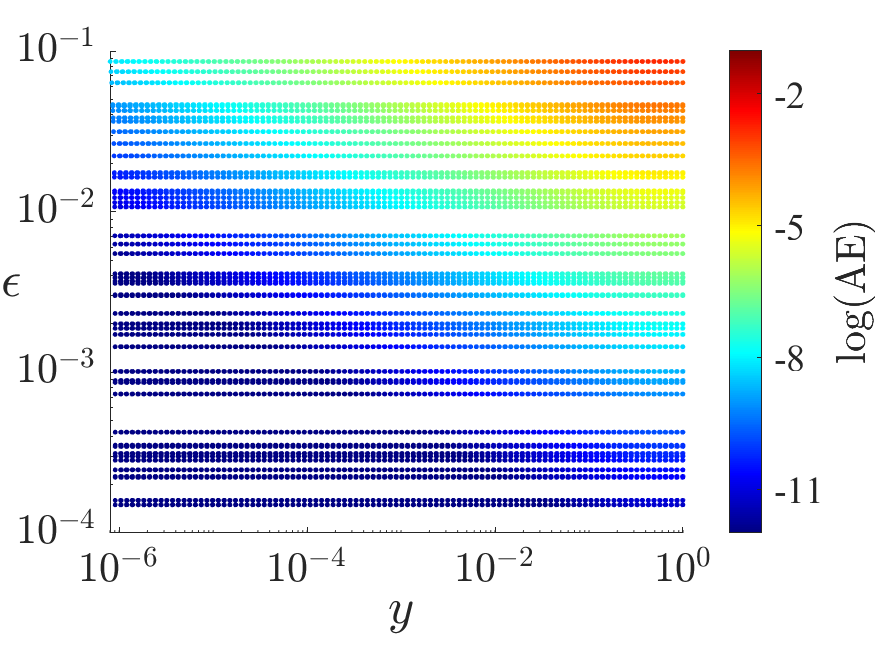}}
    \subfigure[CSP]{
    \includegraphics[width=0.3\textwidth]{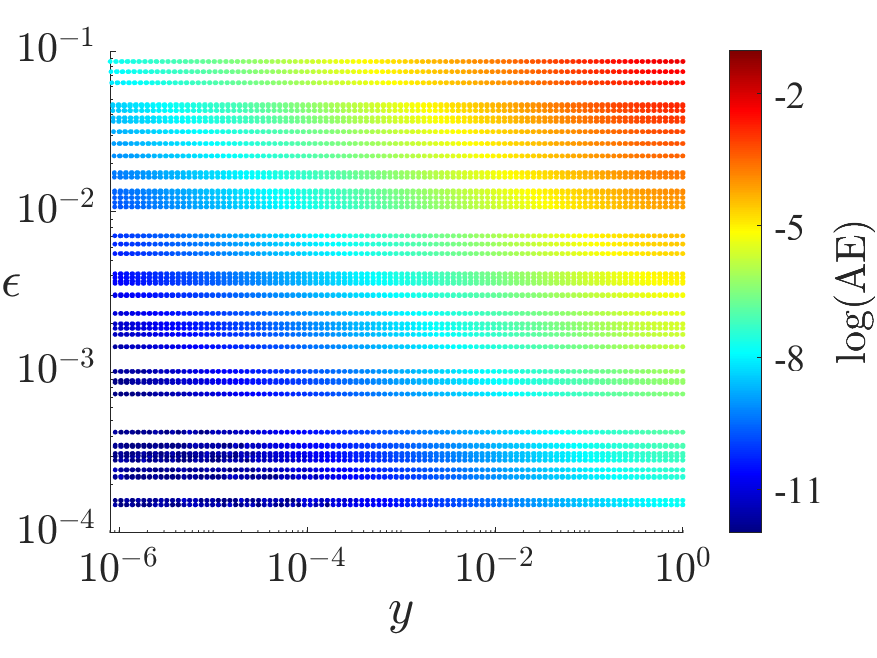}}
    \caption{MM system (\ref{eq:MMss}). Absolute errors (AE) of the SIM approximation in comparison to the numerical solution $x_{i,j}$ of the MM slow subsystem in Eq.~\eqref{eq:MMss}.~Panels (a) and (b) depict the $\lvert x_{i,j} - \mathcal{N}(y_i,\epsilon_j) \rvert$ of the PIML schemes, while panels (c), (d), (e) and (f) depict the $\lvert x_{i,j} - h(y_i,\epsilon_j) \rvert$ of the sQSSA, GSPT $\mathcal{O}(\epsilon)$ and $\mathcal{O}(\epsilon^2)$ asymptotic expansions and the CSP with one iteration approximation in Eqs.~(\ref{eq:MManalSIM}-\ref{eq:MM_CSP}), respectively.}
    \label{fig:MM_AE}
\end{figure}

Fig.~\ref{fig:MM_AE}, depicts the approximation accuracy in terms of $|x_{i,j}-\mathcal{N}(y_{i},\epsilon_j)|$ (for the PIML) and $|x_{i,j}-h(y_{i},\epsilon_j) |$ (for sQSSA, SPT, GSPT and CSP with one iteration).~As shown, the PIML schemes provide a high approximation accuracy, that, for any practical purposes, is not affected by the magnitude of $y$ and $\epsilon$. Furthermore, the SIM approximations provided by the proposed PIML schemes, show a significantly higher numerical accuracy when compared to the sQSSA, GSPT $\mathcal{O}(\epsilon)$ and $\mathcal{O}(\epsilon^2)$, and CSP with one iteration approximations, for large values of $\epsilon$ and the slow variable $y$. The opposite holds for smaller values of $\epsilon$, as expected by the approximation accuracy of the regular asymptotic series expansion.

\subsection{The Target Mediated Drug Disposition mechanism}
For our illustration, we computed the SIM in the domain $(y,z) \in \Omega = [0.2, 2.0]\times[1.3,2.9]$ for $\epsilon \in [10^{-4}, 10^{-1}]$; training, validation and test sets were sampled in this domain.~Following the procedure described in Section~\ref{sec:NumR}, we considered $n_\epsilon =13$ logarithmically spaced values of $\epsilon_j \in [10^{-4}, 10^{-1}]$.~For each $\epsilon_j$, we generated a grid of $5\times5$ initial conditions  $(y(0),z(0))$ selected randomly from a uniform distribution in $[2,2.4]\times [1.3,2.3]$, so that the resulting trajectories sufficiently cover the domain $\Omega$. For each pair of $(y(0),z(0))$, the initial condition $x(0)$ was selected randomly from a uniform distribution in $[0,2]$. From the resulting trajectories, we sampled $n_y=500$ points ($20$ equidistant-in time- points per trajectory). The values of the slow variables $\mathbf{y}_i=[y_i,z_i]^\top$ for $i=1,\ldots,n_y$ were collected until the trajectories were about to exit the boundary $\partial \Omega$ (at $\epsilon/y<5 \epsilon <1$). 

Table~\ref{tab:TMDD_train}  depicts the loss function $\lVert \boldsymbol{\mathcal{F}} \rVert^2_2$ for the training and validation sets of the proposed PIML schemes for all three differentiation schemes (SD, AD, FD).~The corresponding computational costs are also given. Results are obtained by averaging over 10 runs of different randomly sampled training and validations sets.~As shown, the training of SLFNNs using the SD scheme resulted in a better convergence, though of the same order when compared to the other differentiation schemes.
\begin{table}[!h]
    \centering
    \footnotesize
    \begin{tabular}{l| c c | c c c}
    \toprule
    & \multicolumn{2}{|c}{Loss Function $\lVert \boldsymbol{\mathcal{F}} \rVert^2_2$}  & \multicolumn{3}{|c}{Computational times (s)} \\
    PIML scheme & Training & Validation & mean & min & max \\
    \midrule 
    SLFNN AD    &	5.73E$-$04	&	1.49E$-$04	&	6.11E$+$02	&	6.06E$+$02	&	6.15E$+$02	\\
    SLFNN FD    &	8.40E$-$04	&	2.12E$-$04	&	8.97E$+$01	&	8.85E$+$01	&	9.22E$+$01	\\
    SLFNN SD    &	1.54E$-$04	&	4.07E$-$05	&	2.41E$+$01	&	1.15E$+$01	&	4.05E$+$01	\\
    RPNN SD &	1.94E$-$02	&	2.51E$-$02	&	4.30E$+$00	&	3.90E$+$00	&	4.42E$+$00	\\
    \bottomrule
    \end{tabular}
    \caption{TMDD system \eqref{eq:TMDD_SPA}. Loss function $\lVert \boldsymbol{\mathcal{F}} \rVert^2_2$ of the ML schemes for the training and validation sets for the PIML schemes using Automatic Differentiation (AD), Finite Differences (FD) and Symbolic Differentiation (SD). The corresponding  computational times (in seconds) are also given. The results are obtained by averaging over 10 runs.}
    \label{tab:TMDD_train}
\end{table}
For the training process of the SLFNNs, the SD scheme results in $\sim 3\times$ faster times than the FD scheme (computed using parallel computations), which in turn is $\sim 2\times$ faster than the AD scheme.~Regarding RPNNs, the training process resulted in worse convergence of the loss functions in comparison to SLFNNs, albeit a significantly larger number of neurons was used.~However, the computational cost of RPNNs is still lesser, since for their training  $\sim6\times$ faster times than the faster SLFNN scheme were required.

For the test set ($[\mathbf{y}_i,\epsilon_j]^\top$ and the corresponding values of the fast variables $x_{i,j}$), we considered 100 trajectories of the TMDD slow subsystem with $50$ random values of $\epsilon$ varying in $\epsilon \in [10^{-4}, 10^{-1}]$.~For each $\epsilon_j$, $5\times 5$ trajectories were generated with random initial conditions uniformly distributed; $x(0) \in [0, 2]$, $y(0) \in [3,4]$ and $z(0) \in [0,1]$.~From the resulting trajectories, we kept the data after $t=10 \epsilon$ in order for the trajectory to lie on the SIM and recorded 100 equidistant- in time - points per trajectory, including only data in the desired domain $(y,z) \in \Omega = [0.2, 2.0]\times[1.3,2.9]$.~Note that the condition $\epsilon/y<5\epsilon<1$ was taken into account to record data before the trajectories exit the SIM.

In Table~\ref{tab:TMDD_acc}, we report the \emph{overall, with respect to all values of $\epsilon$}, $l^2$, $l^{\infty}$ and MSE, $\lVert x_{i,j} - \mathcal{N}(\mathbf{y}_i, \epsilon_j) \rVert$  approximation errors as obtained by the PIML schemes. We also report the corresponding $\lVert x_{i,j} - h(\mathbf{y}_i, \epsilon_j) \rVert$ errors as obtained by the analytical QSSA, GSPT $\mathcal{O}(\epsilon)$ and $\mathcal{O}(\epsilon^2)$ asymptotic expansions and the CSP with one iteration SIM approximation in Eqs.~(\ref{eq:TMDD_GSPT}, \ref{eq:TMDD_CSP}).
\begin{table}[!h]
    \centering
    \footnotesize
    \begin{tabular}{l| c c c c | c c c c}
    \toprule
    & \multicolumn{4}{|c}{PIML} &  \multicolumn{4}{| c}{analytic GSPT approximations} \\
    \midrule
    Error  &  SLFNN AD &  SLFNN FD  & SLFNN SD  & RPNN SD &   QSSA & GSPT $\mathcal{O}(\epsilon)$ & GSPT $\mathcal{O}(\epsilon^2)$ & CSP \\
    \midrule 
    $l^2$  &	2.72E$-$01	&	3.22E$-$01	&	1.12E$-$01	&	8.40E$-$01	&	8.25E$+$01	&	8.06E$+$01	&	2.19E$+$02	&	1.03E$+$00	\\
    $l^{\infty}$    &	1.11E$-$02	&	1.26E$-$02	&	3.67E$-$03	&	7.21E$-$02	&	3.13E$+$00	&	7.04E$+$00	&	2.83E$+$01	&	3.23E$-$02	\\
    MSE &	6.18E$-$07	&	8.46E$-$07	&	1.19E$-$07	&	5.95E$-$06	&	5.45E$-$02	&	5.20E$-$02	&	3.85E$-$01	&	8.44E$-$06	\\
    \bottomrule
    \end{tabular}
    \caption{TMDD system \eqref{eq:TMDD_SPA}. SIM approximation accuracy \emph{over all simulations for all values of $\epsilon$}, $l^2$, $l^{\infty}$ and MSE approximation errors resulting from the PIML schemes (SLFNN and RPNN) and the analytic GSPT approximations; QSSA, GSPT $\mathcal{O(\epsilon)}$ and $\mathcal{O}(\epsilon^2)$ asymptotic expansions, and the CSP with one iteration approximation in Eqs.~(\ref{eq:TMDD_GSPT}, \ref{eq:TMDD_CSP}), respectively.~The numerical accuracy of each SIM approximation is compared with the numerical solution $x_{i,j}$ of the TMDD slow subsystem in Eq.~\eqref{eq:TMDD_SPA}.}
    \label{tab:TMDD_acc}
\end{table}
As shown, the PIML schemes provide SIM approximations of high accuracy, especially with the use of SLFNNs; again, the overall $l^\infty$ norm is of the order of $1E-02$.~Furthermore, the proposed PIML schemes provide much higher approximation accuracy than the QSSA, GSPT $\mathcal{O(\epsilon)}$ and $\mathcal{O}({\epsilon}^2)$ SIM approximations, and slightly higher (for SLFNNs) or similar (for RPNNs) accuracy than the CSP with one iteration approximation. \par
\begin{figure}[!h]
    \centering
    \subfigure[PIML SLFNN]{
    \includegraphics[width=0.3\textwidth]{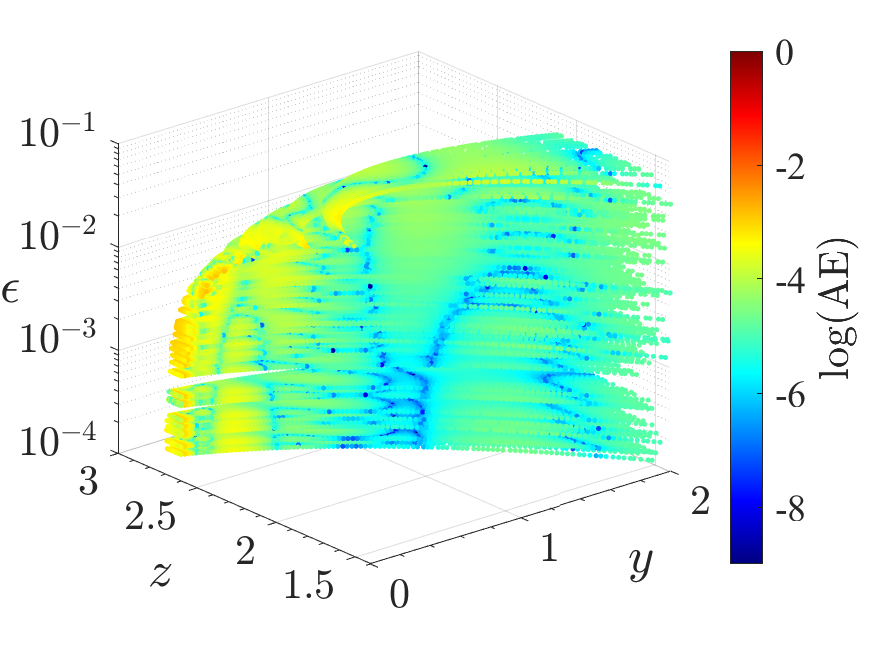}}
    \subfigure[PIML RPNN]{
    \includegraphics[width=0.3\textwidth]{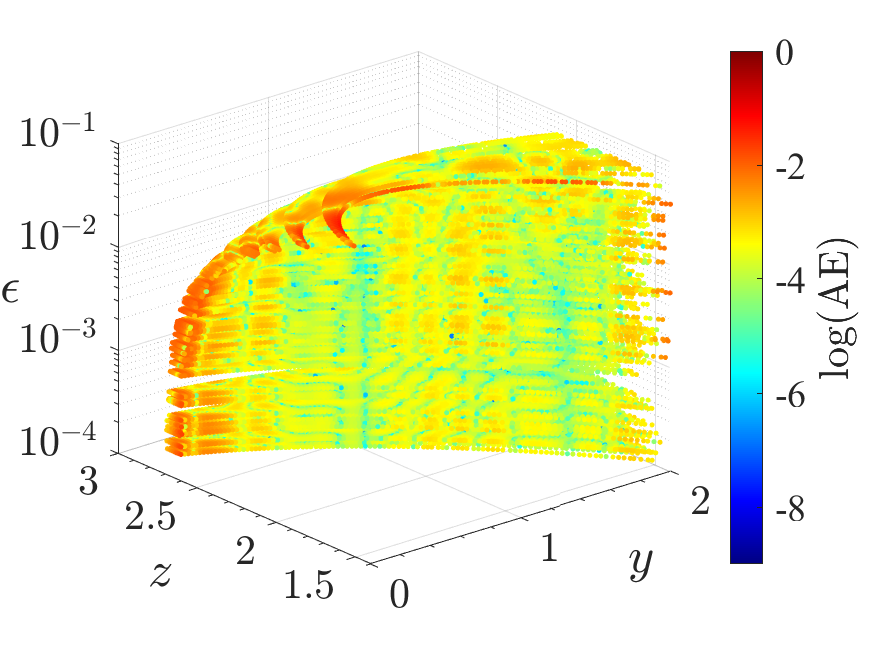}}
    \subfigure[sQSSA]{
    \includegraphics[width=0.3\textwidth]{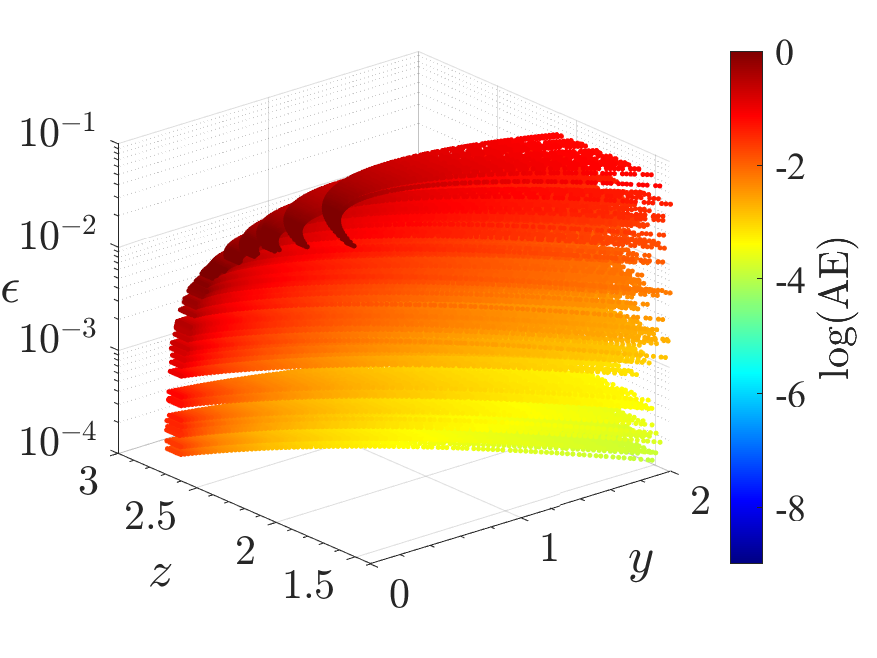}}
    \subfigure[GSPT $\mathcal{O}(\epsilon)$]{
    \includegraphics[width=0.3\textwidth]{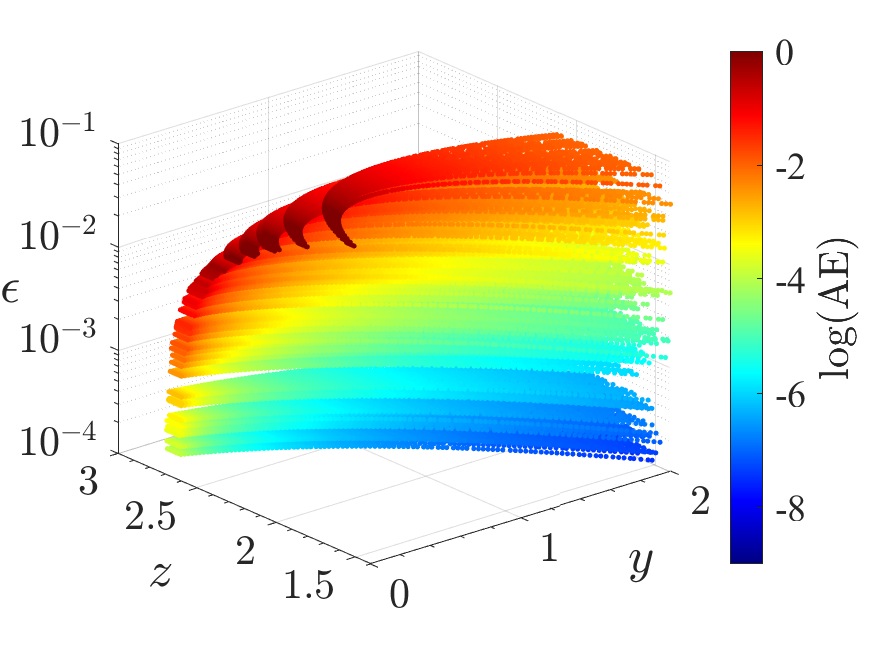}}
    \subfigure[GSPT $\mathcal{O}(\epsilon^2)$]{
    \includegraphics[width=0.3\textwidth]{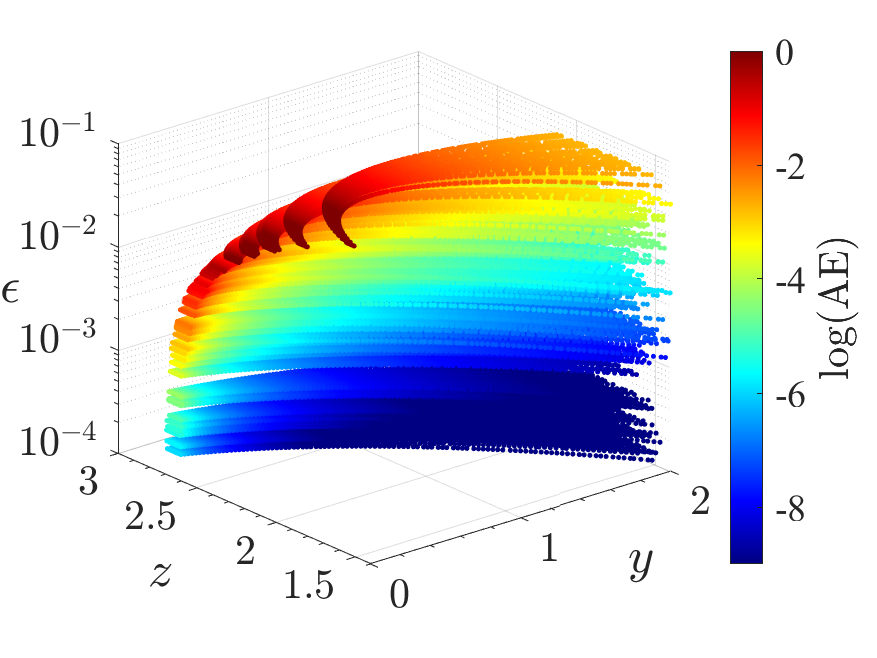}}
    \subfigure[CSP]{
    \includegraphics[width=0.3\textwidth]{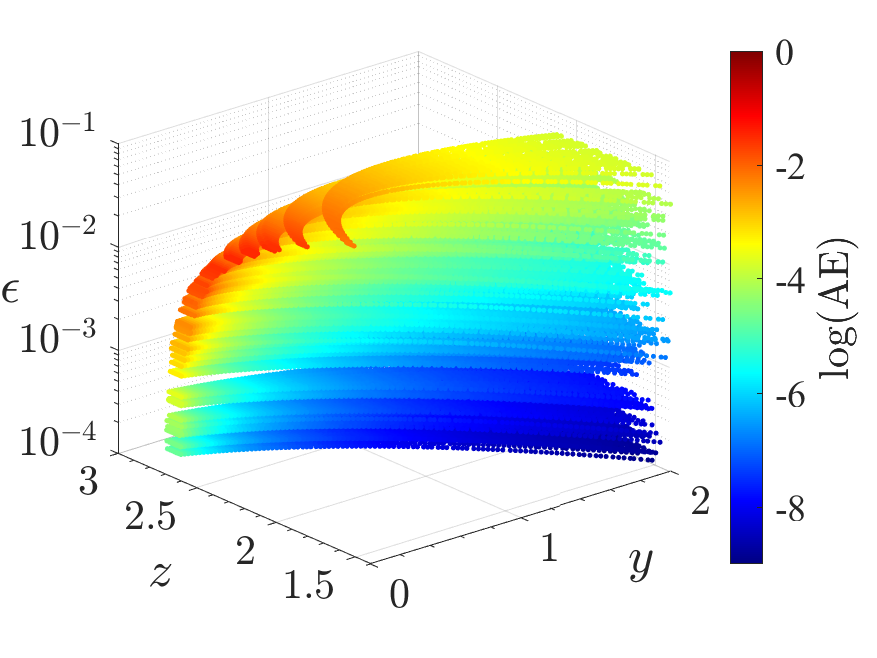}}
    \caption{TMDD system ~\eqref{eq:TMDD_SPA}. Absolute errors (AE) of the SIM approximation in comparison to the numerical solution $x_{i,j}$ of the TMDD slow subsystem in Eq.~\eqref{eq:TMDD_SPA}.~Panels (a) and (b) depict the $\lvert x_{i,j} - \mathcal{N}(\mathbf{y}_i,\epsilon_j) \rvert$ of the PIML schemes, while panels (c), (d), (e) and (f) depict the $\lvert x_{i,j} - h(\mathbf{y}_i,\epsilon_j) \rvert$ of the QSSA, GSPT $\mathcal{O}(\epsilon)$ and $\mathcal{O}(\epsilon^2)$ asymptotic expansions and the CSP with one iteration approximation in Eqs.~(\ref{eq:TMDD_GSPT}, \ref{eq:TMDD_CSP}), respectively.}
    \label{fig:TMDD_AE}
\end{figure}

 
It is interesting to note, that the GSPT-based  $\mathcal{O}({\epsilon}^2)$ asymptotic series expansion, results in worse approximations when compared to the GSPT-based  $\mathcal{O}({\epsilon})$ asymptotic series expansion.~For example, the overall  $l^\infty$ norm is $7.04E+00$ for the GSPT-based $\mathcal{O}(\epsilon)$ approximation, and $2.83E+01$ for the GSPT-based $\mathcal{O}(\epsilon^2)$ one.~This large error is attributed to the relatively larger values of $\epsilon$ for which the GSPT-based asymptotic series expansions is no more valid.~This result is shown in Fig.~\ref{fig:TMDD_AE}, which depicts the approximation accuracy in terms of $|x_{i,j}-\mathcal{N}(\mathbf{y}_{i},\epsilon_j)|$ (for the PIML) and $|x_{i,j}-h(\mathbf{y}_{i},\epsilon_j)|$ (for QSSA, GSPT and CSP with one iteration).~As it is further shown, the proposed PIML approach provides high approximation accuracy that is not affected by the magnitude of $\epsilon$. It outperforms QSSA for all values of $\epsilon$ and GSPT $\mathcal{O}(\epsilon)$ and $\mathcal{O}(\epsilon^2)$ and CSP one iteration approximations of the SIM for high values of $\epsilon$.~For smaller values of $\epsilon$, GSPT $\mathcal{O}(\epsilon^2)$ and CSP are more accurate, as expected.

Finally, it is important to revisit the performance of the PIML schemes close to the boundaries of the SIM; i.e., as $y\rightarrow 0.2$, where the condition $\epsilon/y \ll 1$ is not strongly satisfied.~At this region (see Fig.~\ref{fig:TMDD_AE} for low values of $y$), the PIML schemes result in much better SIM approximations compared to the QSSA, GSPT-based  $\mathcal{O}(\epsilon)$ and $\mathcal{O}({\epsilon}^2)$ asymptotic series expansions.~This large error arises due to the fact that the latter expansions are accurate when the conditions $k_2 \epsilon \ll y$ and $k_3 \epsilon x \ll y$ are additionally satisfied \citep{patsatzis2016asymptotic}; however, these conditions are violated at the boundary of the SIM.~The CSP approximation is not restricted by these conditions \cite{patsatzis2016asymptotic}. Still, while the  CSP scheme with one iteration results, for the particular problem,  in fair approximations, the PIML scheme outperforms it close to the boundary.

\subsection{The 3D Sel'kov model of glycolytic oscillations}

For our illustrations, we computed the SIM around the limit cycle in the domain $(y,z)\in \Omega = [0.2, 1.4]\times[0.3,2.1]$ for $\epsilon \in [10^{-4}, 10^{-1}]$.~Training, validation and test sets were sampled for this domain, a procedure which is trickier than the previous case studies, since for the 3D Sel'kov model we need to collect data both in the interior and exterior of the limit cycle  in the domain $\Omega$.~Thus, after considering $n_\epsilon=13$ logarithmically spaced values of $\epsilon_j \in [10^{-4}, 10^{-1}]$, we generated 25 trajectories for each $\epsilon_j$ with random initial conditions, chosen so that the $20$ of them evolve in the exterior of the limit cycle, while the rest in its interior.~In particular, we have randomly sampled the initial value of the fast variable  $x(0)$ from a uniform distribution in $[0,2]$, while the initial conditions for the slow variables have been randomly sampled from  uniform distributions as follows: (i) $y(0)\in [0.3,0.5] \cup [1.0, 1.2]$ and $z(0)\in [0.4,0.8] \cup [1.6, 2.0]$ for the exterior trajectories, and (ii) $y(0)\in [0.5,0.7]$ and $z(0)\in [1.2,1.4]$ for the interior ones.~From these trajectories, $n_y = 500$ points were sampled ($20$ equidistant in time points per trajectory) and the corresponding values of the slow variables $\mathbf{y}_i = [y_i, z_i]^\top$ for $i=1,\ldots,n_y$ were collected.~Numerical integration was stopped when two subsequent crosses of the periodic orbit with the Poincare section $y=0.7$ had lesser than $0.01$ distance, to avoid significant sampling on the limit cycle.

Table~\ref{tab:TLC_train} depicts the loss function $\lVert \boldsymbol{\mathcal{F}} \rVert^2_2$ for the training and validation sets of the proposed PIML schemes for all three differentiation schemes (SD, AD, FD).~The corresponding computational costs are also provided, obtained by averaging over 10 runs of different randomly sampled training and validations sets.~As shown, the training of SLFNNs using either of the three differentiation schemes resulted in a similar convergence accuracy, while the training of RPNNs resulted in better, by one order of magnitude, convergence.
\begin{table}[!h]
    \centering
    \footnotesize
    \begin{tabular}{l| c c | c c c}
    \toprule
    & \multicolumn{2}{|c}{Loss Function $\lVert \boldsymbol{\mathcal{F}} \rVert^2_2$}  & \multicolumn{3}{|c}{Computational times (s)} \\
    PIML scheme & Training & Validation & mean & min & max \\
    \midrule 
    SLFNN AD    &	5.07E$-$05	&	1.43E$-$05	&	3.79E$+$02	&	3.24E$+$02	&	4.22E$+$02	\\
    SLFNN FD    &	7.27E$-$05	&	2.20E$-$05	&	8.94E$+$01	&	7.64E$+$01	&	9.48E$+$01	\\
    SLFNN SD    &	4.30E$-$05	&	1.33E$-$05	&	3.67E$+$01	&	1.24E$+$01	&	4.76E$+$01	\\
    RPNN SD     &	1.18E$-$06  &	5.52E$-$07	&	3.31E$-$01	&	2.84E$-$01	&	4.34E$-$01	\\
    \bottomrule
    \end{tabular}
    \caption{3D Sel'kov system~\eqref{eq:ToyLC_SPA}. Loss function $\lVert \boldsymbol{\mathcal{F}} \rVert^2_2$ of the ML schemes for the training and validation sets for the PIML schemes using Automatic Differentiation (AD), Finite Differences (FD) and Symbolic Differentiation (SD). The corresponding  computational times (in seconds) are also given. The results are obtained by averaging over 10 runs.}
    \label{tab:TLC_train}
\end{table}
Regarding the computational costs in the training process of the SLFNNs, the SD scheme is $\sim3\times$ faster than  the FD scheme (computed with 20 parallel processors), which in turn is $\sim4\times$ faster than the AD scheme.~As expected, the training of the RPNNs is much more computationally efficient than that of SLFNNs: the faster SLFNN scheme is $\sim100\times$ slower than RPNNs.

For building the test set (the values of $[\mathbf{y}_i,\epsilon_j]^\top$ and the corresponding values of the fast variables $x_{i,j}$), we considered 100 trajectories of the 3D Sel'kov system with $50$ values of $\epsilon$ randomly sampled from a uniform distribution in $[10^{-4}, 10^{-1}]$.~As for the construction of the training set, for each $\epsilon_j$, we generated 25 random initial conditions (sampling from uniform distributions): 20  using $y(0)\in [0,0.2] \cup [1.3, 1.5]$ and $z(0)\in [0,0.5] \cup [2.1, 2.6]$ resulting in trajectories confined in the exterior of the limit cycle, and 5 initial conditions using $y(0)\in [0.5,0.7]$ and $z(0)\in [1.2,1.4]$ resulting in trajectories confined in the interior of the limit cycle; the fast variables were randomly chosen $x(0)\in[0,2]$.~From the resulting trajectories, we kept the data after $t=10 \epsilon$ in order for the trajectory to lie on the SIM and recorded 100 equidistant in time points per trajectory, including only the data in the desired domain $(y,z) \in \Omega = [0.2, 1.4]\times[0.3,2.1]$.~Again, integration was stopped when two subsequent crosses of the periodic orbit with the Poincare section $y=0.7$ had a distance lesser than $0.01$. 

\begin{table}[!h]
    \centering
    \footnotesize
    \begin{tabular}{l| c c c c | c c c c}
    \toprule
    & \multicolumn{4}{|c}{PIML} &  \multicolumn{4}{| c}{analytic GSPT approximations} \\
    \midrule
    Error  &  SLFNN AD &  SLFNN FD  & SLFNN SD  & RPNN SD &   QSSA & GSPT $\mathcal{O}(\epsilon)$ & GSPT $\mathcal{O}(\epsilon^2)$ & CSP \\
    \midrule 
    $l^2$  &	8.28E$-$02	&	9.11E$-$02	&	7.47E$-$02	&	3.32E$-$02	&	2.65E$+$00	&	1.05E$+$00	&	8.25E$-$02	&	2.45E$-$01	\\
    $l^{\infty}$    &	6.56E$-$03	&	6.69E$-$03	&	5.32E$-$03	&	3.10E$-$03	&	6.34E$-$02	&	2.36E$-$02	&	3.16E$-$03	&	8.86E$-$03	\\
    MSE &	5.62E$-$08	&	6.79E$-$08	&	5.66E$-$08	&	9.35E$-$09	&	5.62E$-$05	&	8.79E$-$06	&	5.44E$-$08	&	4.81E$-$07	\\
    \bottomrule
    \end{tabular}
    \caption{3D Sel'kov system~\eqref{eq:ToyLC_SPA}. SIM approximation accuracy in terms of the \emph{overall, i.e. over all simulations for all values of $\epsilon$}, $l^2$, $l^{\infty}$ and MSE approximation errors resulting from the PIML schemes (SLFNN and RPNN) and the analytic GSPT approximations; QSSA, GSPT $\mathcal{O}(\epsilon)$ and $\mathcal{O}(\epsilon^2)$ asymptotic expansions and the CSP with one iteration, SIM approximations in Eqs.~(\ref{eq:TLC_GSPT}, \ref{eq:TLC_CSP}), respectively.~The numerical accuracy of each SIM approximation is compared with the numerical solution $x_{i,j}$ of the Sel'kov 3D model in Eq.~\eqref{eq:ToyLC_SPA}.}
    \label{tab:TLC_acc}
\end{table}
In Table~\ref{tab:TLC_acc}, we report the \emph{overall, with respect to all values of $\epsilon$}, $l^2$, $l^{\infty}$ and MSE, $\lVert x_{i,j} - \mathcal{N}(\mathbf{y}_i, \epsilon_j) \rVert$ approximation errors, obtained by the PIML schemes. We also report the corresponding $\lVert x_{i,j} - h(\mathbf{y}_i, \epsilon_j) \rVert$ errors as obtained by the analytical QSSA, GSPT $\mathcal{O}(\epsilon)$ and $\mathcal{O}(\epsilon^2)$ asymptotic series expansions and the CSP with one iteration, SIM approximations in Eqs.~(\ref{eq:TLC_GSPT}, \ref{eq:TLC_CSP}). As shown, the PIML schemes provide a high numerical approximation accuracy. They outperform QSSA, GSPT $\mathcal{O(\epsilon)}$ and CSP with one iteration approximations, while, for any practical purposes, the PIML schemes result in similar approximation accuracy when compared to the GSPT $\mathcal{O}({\epsilon}^2)$.

\begin{figure}[ht]
    \centering
    \subfigure[PIML SLFNN]{
    \includegraphics[width=0.3\textwidth]{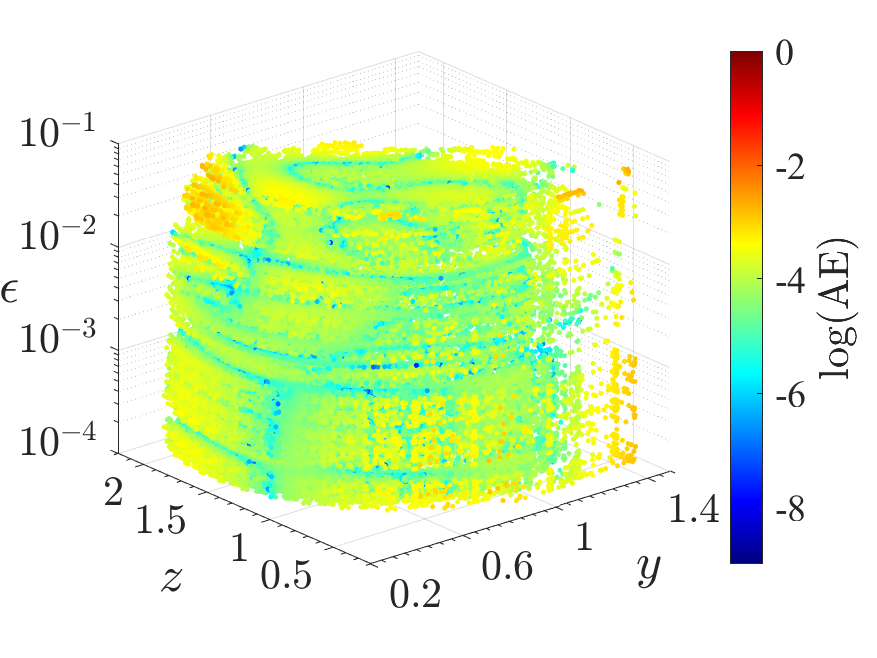}}
    \subfigure[PIML RPNN]{
    \includegraphics[width=0.3\textwidth]{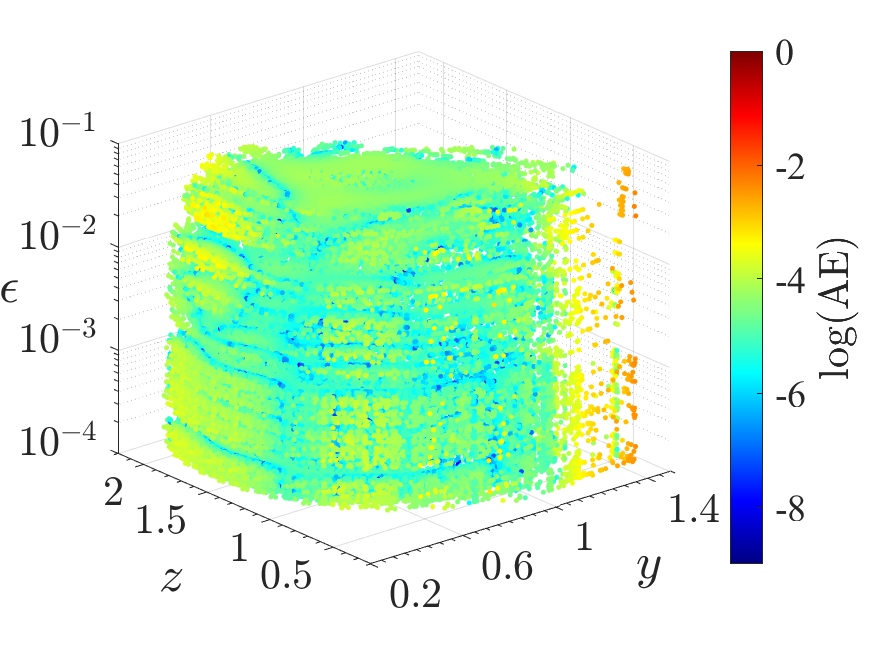}}
    \subfigure[sQSSA]{
    \includegraphics[width=0.3\textwidth]{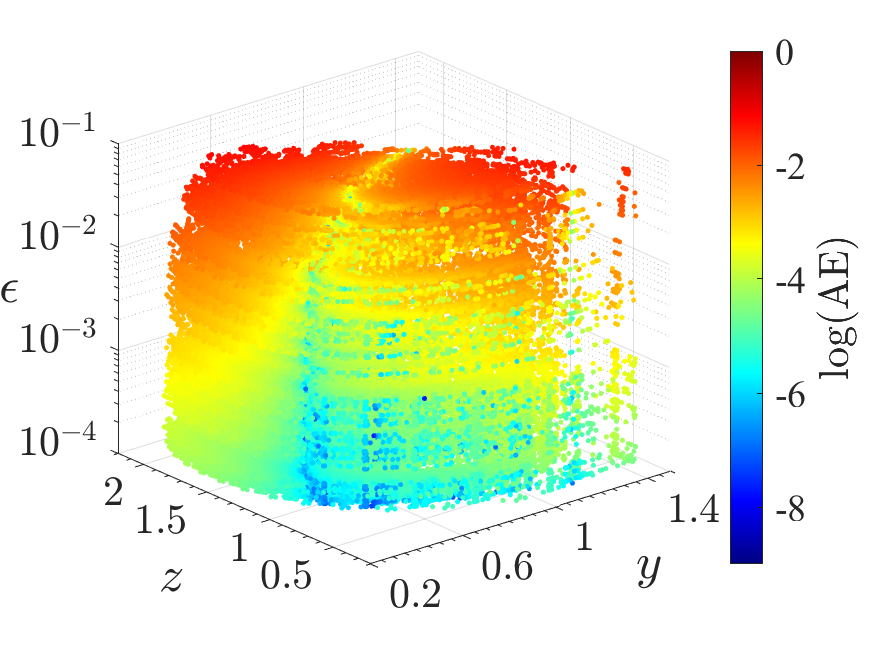}}
    \subfigure[GSPT $\mathcal{O}(\epsilon)$]{
    \includegraphics[width=0.3\textwidth]{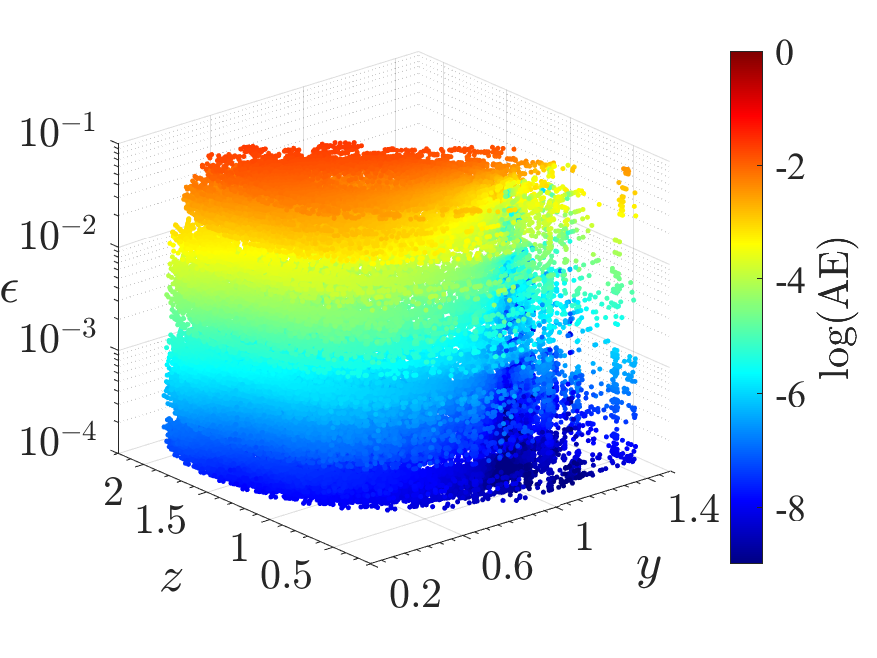}}
    \subfigure[GSPT $\mathcal{O}(\epsilon^2)$]{
    \includegraphics[width=0.3\textwidth]{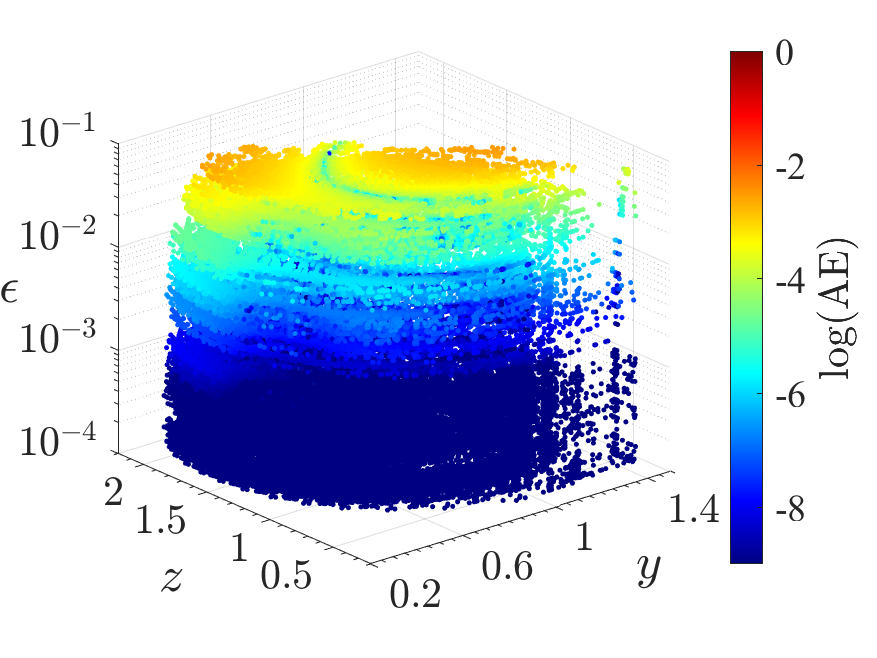}}
    \subfigure[CSP]{
    \includegraphics[width=0.3\textwidth]{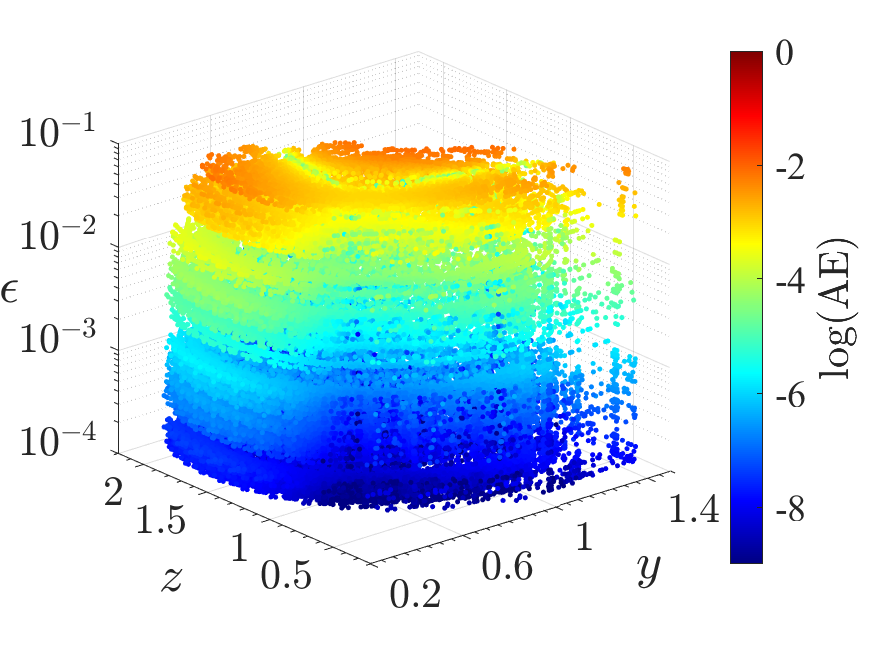}}
    \caption{3D Sel'kov system~\eqref{eq:ToyLC_SPA}. Absolute errors (AE) of the SIM approximation in comparison to the numerical solution $x_{i,j}$ of the 3D Sel'kov system in Eq.~\eqref{eq:ToyLC_SPA}.~Panels (a) and (b) depict the AE $\lvert x_{i,j} - \mathcal{N}(\mathbf{y}_i,\epsilon_j) \rvert$ of the PIML schemes, while panels (c), (d), (e) and (f) depict the AE  $\lvert x_{i,j} - h(\mathbf{y}_i,\epsilon_j) \rvert$ of the QSSA, GSPT $\mathcal{O}(\epsilon)$ and $\mathcal{O}(\epsilon^2)$ asymptotic expansions and the CSP with one iteration, approximations in Eqs.~(\ref{eq:TLC_GSPT}, \ref{eq:TLC_CSP}), respectively.}
    \label{fig:TLC_AE}
\end{figure}

Fig.~\ref{fig:TLC_AE}, depicts the numerical approximation accuracy in terms of 
$\lvert x_{i,j} - \mathcal{N}(\mathbf{y}_i, \epsilon_j) \rvert$ for the PIML schemes and $\lvert x_{i,j} - h(\mathbf{y}_i,\epsilon_j) \rvert$ for the GSPT expressions for every point $x_{i,j}$.~Again, as also in the other two benchmark problems, the approximation accuracy of the PIML schemes is not
affected by the magnitude of $\epsilon$.

\section{Conclusions}

We addressed a physics-informed machine learning approach, within the context of GSPT, for the approximation of SIMs of singularly perturbed dynamical systems. With the assumption of the local invariance of the SIM \cite{fenichel1979geometric}, the proposed PIML approach provides explicit functional forms of SIMs as solutions of the (partial) differential equation (PDE) corresponding to the invariance equation (IE). The proposed PIML approach results in high approximation accuracy for all values of the perturbation parameter $\epsilon$, while traditional GSPT-based methods, such as the asymptotic series expansion approach, fail for relatively large values of $\epsilon$. In addition, as shown in the TMDD benchmark problem, although the GSPT-based approaches, including CSP, lose accuracy close to the boundaries of the underlying SIM as expected \citep{maris2015hidden,jones1995geometric,kuehn2015multiple}, the proposed PIML approach successfully provides accurate SIM approximations there.~The former result stems from the fact that the proposed PIML approach is not limited only locally to small values of $\epsilon$, as the traditional asymptotic series expansion method is, because of the universal approximation properties of neural networks \cite{cybenko1989approximation} and random projection neural networks \cite{barron1993universal,pao1994learning,rahimi2007random,huang2006extreme,fabiani2022parsimonious}. \emph{This is of particular importance and interest, as there are many systems for which, while the gap between the fast and slow timescales is not that big, but still exists. In such cases, one can still construct ROMs (see for example the discussion for this particular issue in \cite{kaper2002asymptotic})}.

More advanced, state-of-the-art techniques, such as the CSP, can be employed when a higher approximation accuracy is being sought. However, they usually result in implicit SIM representations, which may lead to
extraneous explicit solutions (see for example the discussion in \cite{sobolev2012asymptotic}). Furthermore, from the construction of ROMs point of view, implicit forms of SIMs, even if derived in an analytical form, as e.g., in the form of holonomic constraints, may increase the index of the system of differential equations, thus making the task of their numerical integration more difficult. The proposed framework copes with singularly perturbed systems characterized by explicit timescale splitting.~Inevitably, the transformation of the original system to its fast and slow subsystems imposes the assumptions that the dimension of the SIM and the variables associated with the fast dynamics do not vary in the timeframe of interest.~Variations in time are dealt by SPT/GSPT by constructing different slow subsystems when one of the above assumptions changes in time  \citep{fenichel1979geometric,kuehn2015multiple,verhulst2005methods}.~The identification of such alterations is usually intuition-based, with the only notable exception being the systematic diagnostics toolset provided by CSP \citep{lam1989understanding,goussis1992study,valorani2001explicit}.

There are several directions that we aim to pursuit in future works. First, implementing the proposed scheme for finding SIMs of high dimensional systems and compare its performance with other state-of-the-art methods such as CSP but also ILDM.~We further aim to generalize the method for non-standard forms in which $\epsilon$ or the transformations to fast/slow subsystems are not known a-priori, or when the dimension of the SIM and the variables associated with the fast dynamics change in time. Another direction on which the proposed method can be extended with the aid of Equation-free framework \cite{kevrekidis2003equation} is that of the computation of stable, unstable and center invariant manifolds of large-scale microscopic simulators, where explicit equations for the emergent dynamics are not explicitly available \cite{gear2005projecting,siettos2014equation,siettos2022numerical}. Another important direction is that of the quantification of the validity of the SIM approximations provided by the PIML scheme and therefore the corresponding regions of validity of the resulting ROMs \cite{homescu2007error}, with respect to the uncertainty quantification in neural networks \cite{psaros2023uncertainty}.

\bibliographystyle{apalike}


\begin{thebibliography}{1}

\bibitem[Aston et~al., 2011]{aston2011mathematical}
Aston, P.~J., Derks, G., Raji, A., Agoram, B.~M., and van~der Graaf, P.~H.
  (2011).
\newblock Mathematical analysis of the pharmacokinetic--pharmacodynamic (pkpd)
  behaviour of monoclonal antibodies: predicting in vivo potency.
\newblock {\em Journal of theoretical biology}, 281(1):113--121.

\bibitem[Balasubramanian et~al., 2002]{balasubramanian2002isomap}
Balasubramanian, M., Schwartz, E.~L., Tenenbaum, J.~B., de~Silva, V., and
  Langford, J.~C. (2002).
\newblock The {I}somap algorithm and topological stability.
\newblock {\em Science}, 295(5552):7--7.

\bibitem[Barron, 1993]{barron1993universal}
Barron, A.~R. (1993).
\newblock Universal approximation bounds for superpositions of a sigmoidal
  function.
\newblock {\em IEEE Transactions on Information theory}, 39(3):930--945.

\bibitem[Baydin et~al., 2018]{baydin2018automatic}
Baydin, A.~G., Pearlmutter, B.~A., Radul, A.~A., and Siskind, J.~M. (2018).
\newblock Automatic differentiation in machine learning: a survey.
\newblock {\em Journal of Marchine Learning Research}, 18:1--43.

\bibitem[Bollt, 2007]{bollt2007attractor}
Bollt, E. (2007).
\newblock Attractor modeling and empirical nonlinear model reduction of
  dissipative dynamical systems.
\newblock {\em International Journal of Bifurcation and Chaos},
  17(04):1199--1219.

\bibitem[Bollt et~al., 2018]{bollt2018matching}
Bollt, E.~M., Li, Q., Dietrich, F., and Kevrekidis, I. (2018).
\newblock On matching, and even rectifying, dynamical systems through koopman
  operator eigenfunctions.
\newblock {\em SIAM Journal on Applied Dynamical Systems}, 17(2):1925--1960.

\bibitem[Bowen et~al., 1963]{bowen1963singular}
Bowen, J., Acrivos, A., and Oppenheim, A. (1963).
\newblock Singular perturbation refinement to quasi-steady state approximation
  in chemical kinetics.
\newblock {\em Chemical Engineering Science}, 18(3):177--188.

\bibitem[Calabr{\`o} et~al., 2021]{calabro2021extreme}
Calabr{\`o}, F., Fabiani, G., and Siettos, C. (2021).
\newblock Extreme learning machine collocation for the numerical solution of
  elliptic pdes with sharp gradients.
\newblock {\em Computer Methods in Applied Mechanics and Engineering},
  387:114188.

\bibitem[Chen and Ferguson, 2018]{chen2018molecular}
Chen, W. and Ferguson, A.~L. (2018).
\newblock Molecular enhanced sampling with autoencoders: On-the-fly collective
  variable discovery and accelerated free energy landscape exploration.
\newblock {\em Journal of computational chemistry}, 39(25):2079--2102.

\bibitem[Chen et~al., 2021]{chen2021physics}
Chen, W., Wang, Q., Hesthaven, J.~S., and Zhang, C. (2021).
\newblock Physics-informed machine learning for reduced-order modeling of
  nonlinear problems.
\newblock {\em Journal of computational physics}, 446:110666.

\bibitem[Coifman et~al., 2005]{coifman2005geometric}
Coifman, R.~R., Lafon, S., Lee, A.~B., Maggioni, M., Nadler, B., Warner, F.,
  and Zucker, S.~W. (2005).
\newblock Geometric diffusions as a tool for harmonic analysis and structure
  definition of data: diffusion maps.
\newblock {\em Proceedings of the National Academy of Sciences},
  102(21):7426--7431.

\bibitem[Cybenko, 1989]{cybenko1989approximation}
Cybenko, G. (1989).
\newblock Approximation by superpositions of a sigmoidal function.
\newblock {\em Mathematics of control, signals and systems}, 2(4):303--314.

\bibitem[Dong and Yang, 2022a]{dong2022numerical}
Dong, S. and Yang, J. (2022a).
\newblock Numerical approximation of partial differential equations by a
  variable projection method with artificial neural networks.
\newblock {\em Computer Methods in Applied Mechanics and Engineering},
  398:115284.

\bibitem[Dong and Yang, 2022b]{dong2022computing}
Dong, S. and Yang, J. (2022b).
\newblock On computing the hyperparameter of extreme learning machines:
  Algorithm and application to computational pdes, and comparison with
  classical and high-order finite elements.
\newblock {\em Journal of Computational Physics}, 463:111290.

\bibitem[Dsilva et~al., 2016]{dsilva2016data}
Dsilva, C.~J., Talmon, R., Gear, C.~W., Coifman, R.~R., and Kevrekidis, I.~G.
  (2016).
\newblock Data-driven reduction for a class of multiscale fast-slow stochastic
  dynamical systems.
\newblock {\em SIAM Journal on Applied Dynamical Systems}, 15(3):1327--1351.

\bibitem[Fabiani et~al., 2021]{fabiani2021numerical}
Fabiani, G., Calabr{\`o}, F., Russo, L., and Siettos, C. (2021).
\newblock Numerical solution and bifurcation analysis of nonlinear partial
  differential equations with extreme learning machines.
\newblock {\em Journal of Scientific Computing}, 89(2):1--35.

\bibitem[Fabiani et~al., 2023]{fabiani2022parsimonious}
Fabiani, G., Galaris, E., Russo, L., and Siettos, C. (2023).
\newblock Parsimonious physics-informed random projection neural networks for
  initial-value problems of odes and index-1 daes.
\newblock {\em Chaos}, 33:043128.

\bibitem[Fenichel, 1979]{fenichel1979geometric}
Fenichel, N. (1979).
\newblock Geometric singular perturbation theory for ordinary differential
  equations.
\newblock {\em Journal of differential equations}, 31(1):53--98.

\bibitem[Fraser, 1988]{fraser1988steady}
Fraser, S.~J. (1988).
\newblock The steady state and equilibrium approximations: A geometrical
  picture.
\newblock {\em The Journal of chemical physics}, 88(8):4732--4738.

\bibitem[Galaris et~al., 2022]{galaris2022numerical}
Galaris, E., Fabiani, G., Gallos, I., Kevrekidis, I., and Siettos, C. (2022).
\newblock Numerical bifurcation analysis of pdes from lattice boltzmann model
  simulations: a parsimonious machine learning approach.
\newblock {\em Journal of Scientific Computing}, 92(2):1--30.

\bibitem[Galassi et~al., 2022]{galassi2022adaptive}
Galassi, R.~M., Ciottoli, P.~P., Valorani, M., and Im, H.~G. (2022).
\newblock An adaptive time-integration scheme for stiff chemistry based on
  computational singular perturbation and artificial neural networks.
\newblock {\em Journal of Computational Physics}, 451:110875.

\bibitem[Gear et~al., 2005]{gear2005projecting}
Gear, C.~W., Kaper, T.~J., Kevrekidis, I.~G., and Zagaris, A. (2005).
\newblock Projecting to a slow manifold: Singularly perturbed systems and
  legacy codes.
\newblock {\em SIAM Journal on Applied Dynamical Systems}, 4(3):711--732.

\bibitem[Ginoux, 2021]{ginoux2021slow}
Ginoux, J.-M. (2021).
\newblock Slow invariant manifolds of slow--fast dynamical systems.
\newblock {\em International Journal of Bifurcation and Chaos}, 31(07):2150112.

\bibitem[Ginoux et~al., 2008]{ginoux2008slow}
Ginoux, J.-M., Rossetto, B., and Chua, L.~O. (2008).
\newblock Slow invariant manifolds as curvature of the flow of dynamical
  systems.
\newblock {\em International Journal of Bifurcation and Chaos},
  18(11):3409--3430.

\bibitem[Gorban and Karlin, 2003]{gorban2003method}
Gorban, A.~N. and Karlin, I.~V. (2003).
\newblock Method of invariant manifold for chemical kinetics.
\newblock {\em Chemical Engineering Science}, 58(21):4751--4768.

\bibitem[Goussis and Lam, 1992]{goussis1992study}
Goussis, D. and Lam, S. (1992).
\newblock A study of homogeneous methanol oxidation kinetics using csp.
\newblock In {\em Symposium (International) on Combustion}, volume~24, pages
  113--120. Elsevier.

\bibitem[Goussis, 2012]{goussis2012quasi}
Goussis, D.~A. (2012).
\newblock Quasi steady state and partial equilibrium approximations: their
  relation and their validity.
\newblock {\em Combustion Theory and Modelling}, 16(5):869--926.

\bibitem[Goussis and Valorani, 2006]{goussis2006efficient}
Goussis, D.~A. and Valorani, M. (2006).
\newblock An efficient iterative algorithm for the approximation of the fast
  and slow dynamics of stiff systems.
\newblock {\em Journal of Computational Physics}, 214(1):316--346.

\bibitem[Guckenheimer and Holmes, 2013]{guckenheimer2013nonlinear}
Guckenheimer, J. and Holmes, P. (2013).
\newblock {\em Nonlinear oscillations, dynamical systems, and bifurcations of
  vector fields}, volume~42.
\newblock Springer Science \& Business Media.

\bibitem[Hagan and Menhaj, 1994]{hagan1994training}
Hagan, M.~T. and Menhaj, M.~B. (1994).
\newblock Training feedforward networks with the marquardt algorithm.
\newblock {\em IEEE transactions on Neural Networks}, 5(6):989--993.

\bibitem[Homescu et~al., 2007]{homescu2007error}
Homescu, C., Petzold, L.~R., and Serban, R. (2007).
\newblock Error estimation for reduced-order models of dynamical systems.
\newblock {\em Siam Review}, 49(2):277--299.

\bibitem[Huang et~al., 2006]{huang2006extreme}
Huang, G.-B., Zhu, Q.-Y., and Siew, C.-K. (2006).
\newblock Extreme learning machine: theory and applications.
\newblock {\em Neurocomputing}, 70(1-3):489--501.

\bibitem[Igelnik and Pao, 1995]{igelnik1995stochastic}
Igelnik, B. and Pao, Y.-H. (1995).
\newblock Stochastic choice of basis functions in adaptive function
  approximation and the functional-link net.
\newblock {\em IEEE transactions on Neural Networks}, 6(6):1320--1329.

\bibitem[Jones, 1995]{jones1995geometric}
Jones, C.~K. (1995).
\newblock Geometric singular perturbation theory.
\newblock {\em Dynamical systems}, pages 44--118.

\bibitem[Kaper and Kaper, 2002]{kaper2002asymptotic}
Kaper, H.~G. and Kaper, T.~J. (2002).
\newblock Asymptotic analysis of two reduction methods for systems of chemical
  reactions.
\newblock {\em Physica D: Nonlinear Phenomena}, 165(1-2):66--93.

\bibitem[Kaper, 1999]{kaper1999introduction}
Kaper, T.~J. (1999).
\newblock An introduction to geometric methods and dynamical systems theory for
  singular perturbation problems.
\newblock In {\em Proceedings of Symposia in Applied Mathematics}, volume~56,
  pages 85--132. American Mathematical Society.

\bibitem[Karniadakis et~al., 2021]{karniadakis2021physics}
Karniadakis, G.~E., Kevrekidis, I.~G., Lu, L., Perdikaris, P., Wang, S., and
  Yang, L. (2021).
\newblock Physics-informed machine learning.
\newblock {\em Nature Reviews Physics}, 3(6):422--440.

\bibitem[Kevorkian and Cole, 2013]{kevorkian2013perturbation}
Kevorkian, J. and Cole, J.~D. (2013).
\newblock {\em Perturbation methods in applied mathematics}, volume~34.
\newblock Springer Science \& Business Media.

\bibitem[Kevrekidis et~al., 2003]{kevrekidis2003equation}
Kevrekidis, I.~G., Gear, C.~W., Hyman, J.~M., Kevrekidis, P.~G., Runborg, O.,
  Theodoropoulos, C., et~al. (2003).
\newblock Equation-free, coarse-grained multiscale computation: enabling
  microscopic simulators to perform system-level analysis.
\newblock {\em Commun. Math. Sci}, 1(4):715--762.

\bibitem[Kourdis and Goussis, 2013]{kourdis2013glycolysis}
Kourdis, P.~D. and Goussis, D.~A. (2013).
\newblock Glycolysis in saccharomyces cerevisiae: algorithmic exploration of
  robustness and origin of oscillations.
\newblock {\em Mathematical biosciences}, 243(2):190--214.

\bibitem[Kristiansen, 2019]{kristiansen2019geometric}
Kristiansen, K.~U. (2019).
\newblock Geometric singular perturbation analysis of a dynamical target
  mediated drug disposition model.
\newblock {\em Journal of mathematical biology}, 79(1):187--222.

\bibitem[Kuehn, 2015]{kuehn2015multiple}
Kuehn, C. (2015).
\newblock {\em Multiple time scale dynamics}, volume 191.
\newblock Springer.

\bibitem[Lam and Goussis, 1989]{lam1989understanding}
Lam, S.-H. and Goussis, D.~A. (1989).
\newblock Understanding complex chemical kinetics with computational singular
  perturbation.
\newblock In {\em Symposium (International) on Combustion}, volume~22, pages
  931--941. Elsevier.

\bibitem[Lee et~al., 2020]{lee2020coarse}
Lee, S., Kooshkbaghi, M., Spiliotis, K., Siettos, C.~I., and Kevrekidis, I.~G.
  (2020).
\newblock Coarse-scale pdes from fine-scale observations via machine learning.
\newblock {\em Chaos: An Interdisciplinary Journal of Nonlinear Science},
  30(1):013141.

\bibitem[Lee et~al., 2023]{lee2022learning}
Lee, S., Psarellis, Y.~M., Siettos, C.~I., and Kevrekidis, I.~G. (2023).
\newblock Learning black-and gray-box chemotactic pdes/closures from agent
  based monte carlo simulation data.
\newblock {\em Journal of Mathemtical Biology}, 87:15.

\bibitem[Levy, 1994]{levy1994pharmacologic}
Levy, G. (1994).
\newblock Pharmacologic target-mediated drug disposition.
\newblock {\em Clinical Pharmacology \& Therapeutics}, 56(3):248--252.

\bibitem[Linot and Graham, 2020]{linot2020deep}
Linot, A.~J. and Graham, M.~D. (2020).
\newblock Deep learning to discover and predict dynamics on an inertial
  manifold.
\newblock {\em Physical Review E}, 101(6):062209.

\bibitem[Lu et~al., 2021]{lu2021deepxde}
Lu, L., Meng, X., Mao, Z., and Karniadakis, G.~E. (2021).
\newblock Deepxde: A deep learning library for solving differential equations.
\newblock {\em SIAM Review}, 63(1):208--228.

\bibitem[Lusch et~al., 2018]{lusch2018deep}
Lusch, B., Kutz, J.~N., and Brunton, S.~L. (2018).
\newblock Deep learning for universal linear embeddings of nonlinear dynamics.
\newblock {\em Nature communications}, 9(1):4950.

\bibitem[Maas and Pope, 1992]{maas1992simplifying}
Maas, U. and Pope, S.~B. (1992).
\newblock Simplifying chemical kinetics: intrinsic low-dimensional manifolds in
  composition space.
\newblock {\em Combustion and flame}, 88(3-4):239--264.

\bibitem[Mager and Jusko, 2001]{mager2001general}
Mager, D.~E. and Jusko, W.~J. (2001).
\newblock General pharmacokinetic model for drugs exhibiting target-mediated
  drug disposition.
\newblock {\em Journal of pharmacokinetics and pharmacodynamics},
  28(6):507--532.

\bibitem[Maris and Goussis, 2015]{maris2015hidden}
Maris, D.~T. and Goussis, D.~A. (2015).
\newblock The “hidden” dynamics of the r{\"o}ssler attractor.
\newblock {\em Physica D: Nonlinear Phenomena}, 295:66--90.

\bibitem[Michaelis et~al., 1913]{michaelis1913kinetik}
Michaelis, L., Menten, M.~L., et~al. (1913).
\newblock Die kinetik der invertinwirkung.
\newblock {\em Biochem. z}, 49(333-369):352.

\bibitem[Pao et~al., 1994]{pao1994learning}
Pao, Y.-H., Park, G.-H., and Sobajic, D.~J. (1994).
\newblock Learning and generalization characteristics of the random vector
  functional-link net.
\newblock {\em Neurocomputing}, 6(2):163--180.

\bibitem[Papaioannou et~al., 2022]{papaioannou2022time}
Papaioannou, P.~G., Talmon, R., Kevrekidis, I.~G., and Siettos, C. (2022).
\newblock Time-series forecasting using manifold learning, radial basis
  function interpolation, and geometric harmonics.
\newblock {\em Chaos: An Interdisciplinary Journal of Nonlinear Science},
  32(8):083113.

\bibitem[Patsatzis and Goussis, 2019]{patsatzis2019new}
Patsatzis, D.~G. and Goussis, D.~A. (2019).
\newblock A new michaelis-menten equation valid everywhere multi-scale dynamics
  prevails.
\newblock {\em Mathematical biosciences}, 315:108220.

\bibitem[Patsatzis and Goussis, 2023]{patsatzis2023algorithmic}
Patsatzis, D.~G. and Goussis, D.~A. (2023).
\newblock Algorithmic criteria for the validity of quasi-steady state and
  partial equilibrium models: the michaelis–menten reaction mechanism.
\newblock {\em Journal of Mathematical Biology}, 87(27).

\bibitem[Patsatzis et~al., 2016]{patsatzis2016asymptotic}
Patsatzis, D.~G., Maris, D.~T., and Goussis, D.~A. (2016).
\newblock Asymptotic analysis of a target-mediated drug disposition model:
  algorithmic and traditional approaches.
\newblock {\em Bulletin of mathematical biology}, 78(6):1121--1161.

\bibitem[Patsatzis et~al., 2023]{patsatzis2023data}
Patsatzis, D.~G., Russo, L., Kevrekidis, I.~G., and Siettos, C. (2023).
\newblock Data-driven control of agent-based models: An equation/variable-free
  machine learning approach.
\newblock {\em Journal of Computational Physics}, 478:111953.

\bibitem[Peletier and Gabrielsson, 2009]{peletier2009dynamics}
Peletier, L.~A. and Gabrielsson, J. (2009).
\newblock Dynamics of target-mediated drug disposition.
\newblock {\em European Journal of Pharmaceutical Sciences}, 38(5):445--464.

\bibitem[Peletier and Gabrielsson, 2012]{peletier2012dynamics}
Peletier, L.~A. and Gabrielsson, J. (2012).
\newblock Dynamics of target-mediated drug disposition: characteristic profiles
  and parameter identification.
\newblock {\em Journal of pharmacokinetics and pharmacodynamics},
  39(5):429--451.

\bibitem[Psaros et~al., 2023]{psaros2023uncertainty}
Psaros, A.~F., Meng, X., Zou, Z., Guo, L., and Karniadakis, G.~E. (2023).
\newblock Uncertainty quantification in scientific machine learning: Methods,
  metrics, and comparisons.
\newblock {\em Journal of Computational Physics}, 477:111902.

\bibitem[Pye and Chance, 1966]{pye1966sustained}
Pye, K. and Chance, B. (1966).
\newblock Sustained sinusoidal oscillations of reduced pyridine nucleotide in a
  cell-free extract of saccharomyces carlsbergensis.
\newblock {\em Proceedings of the National Academy of Sciences},
  55(4):888--894.

\bibitem[Rahimi and Recht, 2007]{rahimi2007random}
Rahimi, A. and Recht, B. (2007).
\newblock Random features for large-scale kernel machines.
\newblock {\em Advances in neural information processing systems}, 20.

\bibitem[Raissi et~al., 2019]{raissi2019physics}
Raissi, M., Perdikaris, P., and Karniadakis, G.~E. (2019).
\newblock Physics-informed neural networks: A deep learning framework for
  solving forward and inverse problems involving nonlinear partial differential
  equations.
\newblock {\em Journal of Computational physics}, 378:686--707.

\bibitem[Roussel and Fraser, 1991]{roussel1991geometry}
Roussel, M.~R. and Fraser, S.~J. (1991).
\newblock On the geometry of transient relaxation.
\newblock {\em The Journal of chemical physics}, 94(11):7106--7113.

\bibitem[Roweis and Saul, 2000]{roweis2000nonlinear}
Roweis, S.~T. and Saul, L.~K. (2000).
\newblock Nonlinear dimensionality reduction by locally linear embedding.
\newblock {\em science}, 290(5500):2323--2326.

\bibitem[Roy et~al., 2011]{roy2011periodic}
Roy, T., Bhattacharjee, J., and Mallik, A. (2011).
\newblock Periodic orbits in glycolytic oscillators: From elliptic orbits to
  relaxation oscillations.
\newblock {\em The European Physical Journal E}, 34:1--8.

\bibitem[Rudin et~al., 1976]{rudin1976principles}
Rudin, W. et~al. (1976).
\newblock {\em Principles of mathematical analysis}, volume~3.
\newblock McGraw-hill New York.

\bibitem[Santos~Guti{\'e}rrez et~al., 2021]{santos2021reduced}
Santos~Guti{\'e}rrez, M., Lucarini, V., Chekroun, M.~D., and Ghil, M. (2021).
\newblock Reduced-order models for coupled dynamical systems: Data-driven
  methods and the koopman operator.
\newblock {\em Chaos: An Interdisciplinary Journal of Nonlinear Science},
  31(5):053116.

\bibitem[Schnell and Maini, 2000]{schnell2000enzyme}
Schnell, S. and Maini, P.~K. (2000).
\newblock Enzyme kinetics at high enzyme concentration.
\newblock {\em Bulletin of mathematical biology}, 62(3):483--499.

\bibitem[Segel and Slemrod, 1989]{segel1989quasi}
Segel, L.~A. and Slemrod, M. (1989).
\newblock The quasi-steady-state assumption: A case study in perturbation.
\newblock {\em SIAM Review}, 31(3):446--477.

\bibitem[Sel'Kov, 1968]{sel1968self}
Sel'Kov, E. (1968).
\newblock Self-oscillations in glycolysis 1. a simple kinetic model.
\newblock {\em European Journal of Biochemistry}, 4(1):79--86.

\bibitem[Siettos, 2014]{siettos2014equation}
Siettos, C. (2014).
\newblock Equation-free computation of coarse-grained center manifolds of
  microscopic simulators.
\newblock {\em Journal of Computational Dynamics}, 1(2):377--389.

\bibitem[Siettos and Russo, 2022]{siettos2022numerical}
Siettos, C. and Russo, L. (2022).
\newblock A numerical method for the approximation of stable and unstable
  manifolds of microscopic simulators.
\newblock {\em Numerical Algorithms}, 89(3):1335--1368.

\bibitem[Siettos and Bafas, 2002]{siettos2002semiglobal}
Siettos, C.~I. and Bafas, G.~V. (2002).
\newblock Semiglobal stabilization of nonlinear systems using fuzzy control and
  singular perturbation methods.
\newblock {\em Fuzzy Sets and Systems}, 129(3):275--294.

\bibitem[Singer et~al., 2009]{singer2009detecting}
Singer, A., Erban, R., Kevrekidis, I.~G., and Coifman, R.~R. (2009).
\newblock Detecting intrinsic slow variables in stochastic dynamical systems by
  anisotropic diffusion maps.
\newblock {\em Proceedings of the National Academy of Sciences},
  106(38):16090--16095.

\bibitem[Sobolev and Tropkina, 2012]{sobolev2012asymptotic}
Sobolev, V.~A. and Tropkina, E. (2012).
\newblock Asymptotic expansions of slow invariant manifolds and reduction of
  chemical kinetics models.
\newblock {\em Computational Mathematics and Mathematical Physics},
  52(1):75--89.

\bibitem[Tikhonov, 1952]{tikhonov1952systems}
Tikhonov, A.~N. (1952).
\newblock Systems of differential equations containing small parameters in the
  derivatives.
\newblock {\em Matematicheskii sbornik}, 73(3):575--586.

\bibitem[Valorani and Goussis, 2001]{valorani2001explicit}
Valorani, M. and Goussis, D.~A. (2001).
\newblock Explicit time-scale splitting algorithm for stiff problems:
  auto-ignition of gaseous mixtures behind a steady shock.
\newblock {\em Journal of Computational Physics}, 169(1):44--79.

\bibitem[Valorani et~al., 2005]{valorani2005higher}
Valorani, M., Goussis, D.~A., Creta, F., and Najm, H.~N. (2005).
\newblock Higher order corrections in the approximation of low-dimensional
  manifolds and the construction of simplified problems with the csp method.
\newblock {\em Journal of Computational Physics}, 209(2):754--786.

\bibitem[Valorani et~al., 2015]{valorani2015dynamical}
Valorani, M., Paolucci, S., Martelli, E., Grenga, T., and Ciottoli, P.~P.
  (2015).
\newblock Dynamical system analysis of ignition phenomena using the tangential
  stretching rate concept.
\newblock {\em Combustion and Flame}, 162(8):2963--2990.

\bibitem[van~der Graaf et~al., 2016]{van2016topics}
van~der Graaf, P.~H., Benson, N., and Peletier, L.~A. (2016).
\newblock Topics in mathematical pharmacology.
\newblock {\em Journal of Dynamics and Differential Equations},
  28(3):1337--1356.

\bibitem[Verhulst, 2005]{verhulst2005methods}
Verhulst, F. (2005).
\newblock {\em Methods and applications of singular perturbations}.
\newblock Springer.

\bibitem[Vlachas et~al., 2022]{vlachas2022multiscale}
Vlachas, P.~R., Arampatzis, G., Uhler, C., and Koumoutsakos, P. (2022).
\newblock Multiscale simulations of complex systems by learning their effective
  dynamics.
\newblock {\em Nature Machine Intelligence}, 4(4):359--366.

\bibitem[Wan and Sapsis, 2018]{wan2018machine}
Wan, Z.~Y. and Sapsis, T.~P. (2018).
\newblock Machine learning the kinematics of spherical particles in fluid
  flows.
\newblock {\em Journal of Fluid Mechanics}, 857:R2.

\bibitem[Wechselberger, 2020]{wechselberger2020geometric}
Wechselberger, M. (2020).
\newblock {\em Geometric singular perturbation theory beyond the standard
  form}, volume~6.
\newblock Springer.

\bibitem[Williams et~al., 2015]{williams2015data}
Williams, M.~O., Kevrekidis, I.~G., and Rowley, C.~W. (2015).
\newblock A data--driven approximation of the koopman operator: Extending
  dynamic mode decomposition.
\newblock {\em Journal of Nonlinear Science}, 25(6):1307--1346.

\bibitem[Zagaris et~al., 2009]{zagaris2009analysis}
Zagaris, A., Gear, C.~W., Kaper, T.~J., and Kevrekidis, Y.~G. (2009).
\newblock Analysis of the accuracy and convergence of equation-free projection
  to a slow manifold.
\newblock {\em ESAIM: Mathematical Modelling and Numerical Analysis},
  43(4):757--784.

\bibitem[Zagaris et~al., 2004]{zagaris2004analysis}
Zagaris, A., Kaper, H.~G., and Kaper, T.~J. (2004).
\newblock Analysis of the computational singular perturbation reduction method
  for chemical kinetics.
\newblock {\em Journal of Nonlinear Science}, 14(1):59--91.

\end{thebibliography}

%
%
%
\clearpage
\newpage
\renewcommand{\thetable}{A\arabic{table}}  
\renewcommand{\thefigure}{A\arabic{figure}}
\renewcommand{\theequation}{A\arabic{equation}}
\renewcommand{\thesection}{A} 
\setcounter{section}{0}
\setcounter{equation}{0}
\section*{Appendix}

\section{Derivation of the GSPT approximations of the SIM using the invariance equation}
\label{app:A}
Consider a general slow subsystem in the form of Eq.~\eqref{eq:SPslow}:
\begin{equation}
    \epsilon \dfrac{d \mathbf{x}}{dt} =  \mathbf{f}(\mathbf{x},\mathbf{y},\epsilon), \qquad
    \dfrac{d \mathbf{y}}{dt} =   \mathbf{g}(\mathbf{x},\mathbf{y},\epsilon)
    \label{eq:app1}
\end{equation}
for $(\mathbf{x},\mathbf{y}) \in \mathbb{R}^M \times \mathbb{R}^{N-M}$, $0 \le \epsilon \ll 1$ and $\mathbf{f}$, $\mathbf{g}$ sufficiently smooth functions and let $C_0$ denote the critical manifold $C_0=\{(\mathbf{x},\mathbf{y}) \in \mathbb{R}^N : \mathbf{f}(\mathbf{x},\mathbf{y},0)=\mathbf{0} \}$.~As discussed in Section \ref{sec:Methods}, under the assumptions of the Fenichel-Tikhonov theorem \citep{fenichel1979geometric,tikhonov1952systems}, the SIM $S_{\epsilon}$, that is $\mathcal{O}(\epsilon)$ close to $S_0 \subset C_0$, can be locally represented by the regular asymptotic expansion of Eq.~\eqref{eq:SIMrexp}:
\begin{equation}
    \mathbf{x}=\mathbf{h}_{\epsilon}(\mathbf{y}) = \mathbf{h}_0(\mathbf{y}) + \epsilon \mathbf{h}_1(\mathbf{y}) + \epsilon^2 \mathbf{h}_2(\mathbf{y}) +\ldots + \epsilon^j \mathbf{h}_j(\mathbf{y}) +  \mathcal{O}(\epsilon^{j+1})  
    \label{eq:app2}
\end{equation}
Due to the invariance of the SIM, the latter expression (by differentiation) results to the invariance equation in Eq.~\eqref{eq:Inv}:
\begin{equation}
    \mathbf{f}(\mathbf{h}(\mathbf{y},\epsilon),\mathbf{y},\epsilon) - \epsilon \nabla_\mathbf{y}\mathbf{h}_{\epsilon}(\mathbf{y}) \mathbf{g}(\mathbf{h}(\mathbf{y},\epsilon)) =0
    \label{eq:app3}
\end{equation} 
which can be used to calculate the terms $\mathbf{h}_j(\mathbf{y})$ of Eq.~\eqref{eq:app2} up to the desired order $j$.~The procedure begins by expressing the function $\mathbf{f}$, $\mathbf{g}$ as Taylor expansions near $\epsilon=0$, as:
\begin{align}
    \mathbf{f}(\mathbf{h}_{\epsilon}(\mathbf{y}),\mathbf{y},\epsilon) & = \mathbf{f}_0(\mathbf{y}) + \epsilon \mathbf{f}_1(\mathbf{y}) + \epsilon^2 \mathbf{f}_2(\mathbf{y}) + \ldots + \epsilon^j \mathbf{f}_j(\mathbf{y}) +  \mathcal{O}(\epsilon^{j+1}), \nonumber \\
    \mathbf{g}(\mathbf{h}_{\epsilon}(\mathbf{y}),\mathbf{y},\epsilon) & = \mathbf{g}_0(\mathbf{y}) + \epsilon \mathbf{g}_1(\mathbf{y}) + \epsilon^2 \mathbf{g}_2(\mathbf{y}) + \ldots + \epsilon^j \mathbf{g}_j(\mathbf{y}) +  \mathcal{O}(\epsilon^{j+1}),
    \label{eq:app4}
\end{align}
up to the desired order of accuracy $j$.~Then, the employment of the invariance equation in Eq.~\eqref{eq:app3} implies:
\begin{equation}
    \mathbf{f}_0(\mathbf{y}) + \epsilon \mathbf{f}_1(\mathbf{y}) + \epsilon^2 \mathbf{f}_2(\mathbf{y}) + \mathcal{O}(\epsilon^3) - \epsilon \left( \dfrac{d\mathbf{h}_0(\mathbf{y})}{d\mathbf{y}} + \epsilon \dfrac{d\mathbf{h}_1(\mathbf{y})}{d\mathbf{y}} + \epsilon^2 \dfrac{d\mathbf{h}_2(\mathbf{y})}{d\mathbf{y}} + \mathcal{O}(\epsilon^3) \right) (\mathbf{g}_0(\mathbf{y}) + \epsilon \mathbf{g}_1(\mathbf{y}) + \epsilon^2 \mathbf{g}_2(\mathbf{y}) + \mathcal{O}(\epsilon^3)) = 0,
    \label{eq:app5}
\end{equation}
where the $\mathcal{O}(\epsilon^3)$ and the higher order terms are not included for simplicity.~Now, equating the same order terms in Eq.~\eqref{eq:app5}, one retrieves a system of $j+1$ equations for determining the terms $\mathbf{h}_j(\mathbf{y})$ of Eq.~\eqref{eq:app2}.~For example, selecting $j=2$ one retrieves the system:
\begin{equation}
    \mathbf{f}_0(\mathbf{y}) = \mathbf{0}, \qquad \mathbf{f}_1 (\mathbf{y}) - \dfrac{d\mathbf{h}_0(\mathbf{y})}{d\mathbf{y}} \mathbf{g}_0(\mathbf{y}) = \mathbf{0}, \qquad \mathbf{f}_2(\mathbf{y})-\dfrac{d\mathbf{h}_0(\mathbf{y})}{d\mathbf{y}} \mathbf{g}_1(\mathbf{y})-\dfrac{d\mathbf{h}_1(\mathbf{y})}{d\mathbf{y}} \mathbf{g}_0(\mathbf{y}) = \mathbf{0}
    \label{eq:app6}
\end{equation}
where $\mathbf{f}_0 (\mathbf{y}) \equiv \mathbf{f}_0 (\mathbf{h}_0(\mathbf{y}),\mathbf{y})$, $\mathbf{g}_0 (\mathbf{y}) \equiv \mathbf{g}_0 (\mathbf{h}_0(\mathbf{y}),\mathbf{y})$, $\mathbf{f}_1 (\mathbf{y}) \equiv \mathbf{f}_1 (\mathbf{h}_0(\mathbf{y})+ \epsilon\mathbf{h}_1(\mathbf{y}),\mathbf{y})$, $\mathbf{g}_1 (\mathbf{y}) \equiv \mathbf{g}_1 (\mathbf{h}_0(\mathbf{y})+ \epsilon\mathbf{h}_1(\mathbf{y}),\mathbf{y})$ and $\mathbf{f}_2 (\mathbf{y}) \equiv \mathbf{f}_2 (\mathbf{h}_0(\mathbf{y})+ \epsilon\mathbf{h}_1(\mathbf{y})+\epsilon^2\mathbf{h}_2(\mathbf{y}),\mathbf{y})$, thus allowing for the determination of $\mathbf{h}_0(\mathbf{y})$ from the first equation in Eq.~\eqref{eq:app6}, $\mathbf{h}_1(\mathbf{y})$ from the second one and $\mathbf{h}_2(\mathbf{y})$ from the third one.~This procedure can be used for determining the terms $\mathbf{h}_j(\mathbf{y})$ of Eq.~\eqref{eq:app2} up to the desired order $j$.~Herein, we are interested in the zeroth, first and second order terms, $\mathbf{h}_0$, $\mathbf{h}_1$ and $\mathbf{h}_2$, of the SIM approximations for the MM, TMDD and Sel'kov 3D systems.~In what follows, we demonstrate the calculation of these terms for deriving the expressions in Eqs.~(\ref{eq:MM_GSPT}, \ref{eq:TMDD_GSPT}, \ref{eq:TLC_GSPT}) in order to prove the first part of Lemmas 1, 2 and 3.

\subsection{MM system: proof of Lemma 1}
\label{app:MM_SIM_GSPTexp}

Considering the MM subsystem in Eq.~\eqref{eq:MMss}, the functions $f(h_\epsilon(y),y)$ and $g(h_\epsilon(y),y)$ are written in the form of Eq.~\eqref{eq:app4}, as:
\begin{align}
    f(h_{\epsilon}(y),y,\epsilon) & = y - \dfrac{\kappa + y}{\kappa + 1} h_0(y) - \epsilon \dfrac{\kappa + y}{\kappa + 1} h_1(y) -\epsilon^2 \dfrac{\kappa + y}{\kappa + 1} h_2(y) + \mathcal{O}(\epsilon^3), \nonumber \\
    g(h_{\epsilon}(y),y,\epsilon) & = \dfrac{-(\kappa +1) y + (\kappa-\sigma + y) h_0(y)}{\sigma} + \epsilon \dfrac{(\kappa-\sigma + y) h_1(y)}{\sigma} + \epsilon^2 \dfrac{(\kappa-\sigma + y) h_2(y)}{\sigma}+  \mathcal{O}(\epsilon^3),
    \label{eq:appMM1}
\end{align}
where the third order terms are neglected.~Substitution of $f_0(y)$ in Eq.~\eqref{eq:appMM1} to the first equation in Eq.~\eqref{eq:app6}, yields:
\begin{equation}
    y - \dfrac{\kappa + y}{\kappa + 1} h_0(y) = 0 \Rightarrow h_0(y) = \dfrac{k+1}{k+y} y,
    \label{eq:MM_IE1}
\end{equation}
for obtaining the zeroth order term $h_0(y)$ of the asymptotic expansion in Eq.~\eqref{eq:app2}.~Given $h_0(y)$, substitution of $f_1(y)$ and $g_0(y)$ in Eq.~\eqref{eq:appMM1} to the second equation in Eq.~\eqref{eq:app6}, implies:
\begin{equation}
    \dfrac{\kappa + y}{\kappa + 1} h_1(y) - \dfrac{dh_0(y)}{dy} \dfrac{-(\kappa +1) y + (\kappa-\sigma + y) h_0(y)}{\sigma} = 0 \Rightarrow h_1(y) = \dfrac{\kappa (\kappa+1)^3 y}{(\kappa+y)^4}
    \label{eq:MM_IE2}
\end{equation}
for obtaining the first order term $h_1(y)$ of the asymptotic expansion in Eq.~\eqref{eq:app2}.~Now, given $h_0(y)$ and $h_1(y)$, substitution of $f_2(y)$ and $g_1(y)$ in Eq.~\eqref{eq:appMM1} to the third equation in Eq.~\eqref{eq:app6} yields:
\begin{align}
    & \dfrac{\kappa + y}{\kappa + 1} h_2(y) - \dfrac{dh_0(y)}{dy} \dfrac{(\kappa-\sigma + y) h_1(y)}{\sigma} -\dfrac{dh_1(y)}{dy} \dfrac{-(\kappa +1) y + (\kappa-\sigma + y) h_0(y)}{\sigma} = 0 \Rightarrow \nonumber \\
    & h_2(y) = -\dfrac{\kappa (\kappa+1)^5 y(\kappa^2+3\sigma y+\kappa (y-2\sigma))}{\sigma(\kappa+y)^7},
    \label{eq:MM_IE3}
\end{align}
for obtaining the second order term $h_2(y)$ of the asymptotic expansion in Eq.~\eqref{eq:app2}.~According to the expressions in Eqs.~(\ref{eq:MM_IE1}-\ref{eq:MM_IE3}), the $\mathcal{O}(\epsilon^2)$ regular expansion of the SIM for the MM system is given by the expression:
\begin{equation}
    x 
    = \dfrac{\kappa+1}{\kappa+y} y + \epsilon \dfrac{\kappa (\kappa+1)^3 y}{(\kappa+y)^4} - \epsilon^2 \dfrac{\kappa (\kappa+1)^5 y(\kappa^2+3\sigma y+\kappa (y-2\sigma))}{\sigma(\kappa+y)^7}+\mathcal{O}(\epsilon^3)
    \label{eq:appMM_GSPT}
\end{equation}
thus recovering the expression in Eq.~\eqref{eq:MM_GSPT}.

\subsection{TMDD system: proof of Lemma 2}
\label{app:TMDD_SIM_GSPTexp}

Considering the TMDD system in Eq.~\eqref{eq:TMDD_SPA}, the $2$-dim. SIM approximation $x=h_{\epsilon}(\mathbf{y})=h_{\epsilon}(y,z)$ is now a function of the slow variables $\mathbf{y} = [ y, z ]^\top$.~Hence, in the TMDD case, the functions in Eq.~\eqref{eq:app1} are $f(x,y,z,\epsilon): \mathbb{R}^3 \times \mathbb{R} \mapsto \mathbb{R}$ and $\mathbf{g}(x,y,z,\epsilon)=[g^{(1)}(x,y,z,\epsilon), g^{(2)}(x,y,z,\epsilon)]^\top: \mathbb{R}^3 \times \mathbb{R} \mapsto \mathbb{R}^2$, with $0 \le \epsilon \ll 1$.~Casting these functions to the form of Eq.~\eqref{eq:app4}, results to:
\begin{align}
    f(h_{\epsilon}(y,z),y,z,\epsilon) & = -h_0(y,z) y+k_1z+1 - \epsilon \left( h_1(y,z) y +k_2 h_0(y,z) \right) - \epsilon^2 \left( h_2(y,z) y +k_2 h_1(y,z) \right) + \mathcal{O}(\epsilon^3), \nonumber \\
    g^{(1)}(h_{\epsilon}(y,z),y,z,\epsilon) & = k_3(-h_0(y,z)y+k_1z)-k_4y - \epsilon k_3 h_1(y,z) y - \epsilon^2 k_3 h_2(y,z) y + \mathcal{O}(\epsilon^3), \nonumber \\
    g^{(2)}(h_{\epsilon}(y,z),y,z,\epsilon) & = k_2(h_0(y,z) y-k_1z)-z + \epsilon k_2 h_1(y,z) y + \epsilon^2 k_2 h_2(y,z) y + \mathcal{O}(\epsilon^3),
    \label{eq:appTMDD1}
\end{align}
where the $\mathcal{O}(\epsilon^3)$ terms are truncated.~Substitution of $f_0(y,z)$ in Eq.~\eqref{eq:appTMDD1} to the first equation in Eq.~\eqref{eq:app6}, yields:
\begin{equation}
   -h_0(y,z) y+k_1z+1 = 0 \Rightarrow h_0(y,z) = \dfrac{1+ k_1 z}{y},
    \label{eq:TMDD_IE1}
\end{equation}
for obtaining the zeroth order term $h_0(y,z)$ of the asymptotic expansion in Eq.~\eqref{eq:app2}.~Given $h_0(y,z)$, substitution of $f_1(y,z)$, $g^{(1)}_0(y,z)$ and $g^{(2)}_0(y,z)$ in Eq.~\eqref{eq:appTMDD1} to the second equation in Eq.~\eqref{eq:app6} implies:
\begin{align}
    & h_1(y,z) y + k_2 h_0(y,z) + \dfrac{\partial h_0(y,z)}{\partial y} \left( k_3(-h_0(y,z)y+k_1z)-k_4y \right) + \dfrac{\partial h_0(y,z)}{\partial z} \left( k_2(h_0(y,z) y-k_1z)-z \right) = 0 \Rightarrow \nonumber \\ & h_1(y,z) = -\dfrac{k_3(1+k_1)z+y(k_2+k_1k_2+k_4+k_1z(-1+k_2+k_4))}{y^3},    \label{eq:TMDD_IE2}
\end{align}
for obtaining the first order term $h_1(y,z)$ of the asymptotic expansion in Eq.~\eqref{eq:app2}.~Now, given $h_0(y,z)$ and $h_1(y,z)$, substitution of $f_2(y,z)$, $g^{(1)}_1(y,z)$ and $g^{(2)}_1(y,z)$ in Eq.~\eqref{eq:appTMDD1} to the third equation in Eq.~\eqref{eq:app6}, yields:
\begin{align}
    &  h_2(y,z) y +k_2 h_1(y,z) - \dfrac{\partial h_0(y,z)}{\partial y} k_3 h_1(y,z) y + \dfrac{\partial h_0(y,z)}{\partial z} k_2 h_1(y,z) y + \dfrac{\partial h_1(y,z)}{\partial y} \left(k_3(-h_0(y,z)y+k_1z)- \right. \nonumber \\
    &  \left.k_4y\right) + \dfrac{\partial h_1(y,z)}{\partial z} (k_2(h_0(y,z) y-k_1z)-z) = 0 \Rightarrow \nonumber \\ 
    & h_2(y,z) = \dfrac{k_3^2 (k_1 z+1) (k_1 z + 4) + y^2 ((k_2 + k_4) (k_2 + 2 k_4) + k_1 k_2 (3 k_2 + 4 k_4-1) + k_1 (k_2 + k_4-1) (k_2 + 2 k_4-1) z+}{y^5}  \nonumber \\
    & \dfrac{k_1^2 k_2 (k_2 + ( k_2 + k_4-1) z)) + k_3 y (6 k_4 + k_1 z (7 k_4 + k_1 (k_4-1) z-4) + k_2 (4 + k_1 (5 + z (5 + k_1 (2 + z)))))}{y^5},    \label{eq:TMDD_IE3}
\end{align}
resulting in the second order term $h_2(y,z)$ of the asymptotic expansion in Eq.~\eqref{eq:app2}.~Collecting the expressions $h_0(y,z)$, $h_1(y,z)$ and $h_2(y,z)$ in Eqs.~(\ref{eq:TMDD_IE1}-\ref{eq:TMDD_IE3}), the $\mathcal{O}(\epsilon^2)$ regular expansion of the SIM for the TMDD system is obtained, in the form:
\begin{equation}
    x = h_0(y,z)+\epsilon h_1(y,z)+\epsilon^2 h_2(y,z) + \mathcal{O}(\epsilon^3) 
    \label{eq:appTMDD_GSPT}
\end{equation}
as presented in Eq.~\eqref{eq:TMDD_GSPT}.

\subsection{Sel'kov 3D system: proof of Lemma 3}
\label{app:TLC_SIM_GSPTexp} 

Considering the Sel'kov 3D system in Eq.~\eqref{eq:ToyLC_SPA}, the functions in Eq.~\eqref{eq:app1} now become $f(x,y,z,\epsilon): \mathbb{R}^3 \times \mathbb{R} \mapsto \mathbb{R}$ and $\mathbf{g}(x,y,z,\epsilon)=[g^{(1)}(x,y,z,\epsilon), g^{(2)}(x,y,z,\epsilon)]^\top: \mathbb{R}^3 \times \mathbb{R} \mapsto \mathbb{R}^2$, with $0 \le \epsilon \ll 1$.~Casting these functions to the form of Eq.~\eqref{eq:app4}, results to:
\begin{align}
    f(h_{\epsilon}(y,z),y,z,\epsilon) & = y^2z-k h_0(y,z) y - \epsilon k h_1(y,z) y - \epsilon^2 k h_2(y,z) y + \mathcal{O}(\epsilon^3), \nonumber \\
    g^{(1)}(h_{\epsilon}(y,z),y,z,\epsilon) & = az + y^2z-y + \epsilon h_0(y,z) + \epsilon^2 h_1(y,z) + \mathcal{O}(\epsilon^3) \nonumber \\
    g^{(2)}(y,z) & = -a z -y^2z+b
    \label{eq:appTLC1}
\end{align}
where the $\mathcal{O}(\epsilon^3)$ terms are truncated.~Note that the $g^{(2)}(y,z)$ is independent of the SIM approximation and $\epsilon$ in this case study.~Substitution of $f_0(y,z)$ in Eq.~\eqref{eq:appTLC1} to the first equation in Eq.~\eqref{eq:app6}, yields:
\begin{equation}
   y^2z-k h_0(y,z) y = 0 \Rightarrow h_0(y,z) = \dfrac{y z}{k}, \qquad \text{for} \quad y \neq 0
    \label{eq:TLC_IE1}
\end{equation}
for obtaining the zeroth order term $h_0(y,z)$ of the asymptotic expansion in Eq.~\eqref{eq:app2}.~Given $h_0(y,z)$, substitution of $f_1(y,z)$, $g^{(1)}_0(y,z)$ and $g^{(2)}_0(y,z)$ in Eq.~\eqref{eq:appTLC1} to the second equation in Eq.~\eqref{eq:app6} yields:
\begin{align}
    & k h_1(y,z) y  + \dfrac{\partial h_0(y,z)}{\partial y} \left( az + y^2z-y \right) + \dfrac{\partial h_0(y,z)}{\partial z} \left( -az - y^2z + b \right) = 0 \Rightarrow \nonumber \\ & h_1(y,z) = \dfrac{z-b}{k^2} + \dfrac{z(a+y^2)(y-z)}{k^2y}    \label{eq:TLC_IE2}
\end{align}
for obtaining the first order term $h_1(y,z)$ of the asymptotic expansion in Eq.~\eqref{eq:app2}.~Now, given $h_0(y,z)$ and $h_1(y,z)$, substitution of $f_2(y,z)$ and $g^{(1)}_1(y,z)$ in Eq.~\eqref{eq:appTLC1} to the third equation in Eq.~\eqref{eq:app6}, yields:
\begin{align}
     & k h_2(y,z) y + \dfrac{\partial h_0(y,z)}{\partial y} h_0(y,z) + \dfrac{\partial h_1(y,z)}{\partial y} \left( az + y^2z -y \right) + \dfrac{\partial h_1(y,z)}{\partial z} \left( -az - y^2z -b \right)  = 0 \Rightarrow \nonumber \\ & h_2(y,z) =  \dfrac{b y (2 (a + y^2) z-y (1 + a + y^2) ) + z (a y (y + 2 y^3 + z - 6 y^2 z) + a^2 (y^2 - 2 y z - z^2) +  }{k^3 y^3}  \nonumber \\
     & \dfrac{y^3 (-2 z + y (3 + y^2 - 4 y z + z^2))) }{k^3 y^3},  \label{eq:TLC_IE3}
\end{align}
resulting to the second order term $h_2(y,z)$ of the asymptotic expansion in Eq.~\eqref{eq:app2}.~Collecting all terms $h_0(y,z)$, $h_1(y,z)$ and $h_2(y,z)$ in Eqs.~(\ref{eq:TLC_IE1}- \ref{eq:TLC_IE3}), the $\mathcal{O}(\epsilon^2)$ regular expansion of the SIM for the 3D Sel'kov system is recovered in the form: 
\begin{equation}
     x = h_0(y,z)+\epsilon h_1(y,z) + \epsilon^2 h_2(y,z) + \mathcal{O}(\epsilon^3),
    \label{eq:appTLC_GSPT}
\end{equation}
as presented in Eq.~\eqref{eq:TLC_GSPT}.

\renewcommand{\thesection}{B} 
\renewcommand{\theequation}{B\arabic{equation}}
\setcounter{equation}{0}
\section{Derivation of the GSPT approximations of the SIM using the CSP method with one iteration}
\label{app:B}

Computational Singular Perturbation (CSP) is an algorithmic methodology employed in the context of GSPT for the derivation of the SIM, as well as the identification of its geometrical properties \citep{lam1989understanding,goussis1992study,valorani2005higher,goussis2012quasi}.~As already discussed in Section~\ref{sec:Methods}, sophisticated computational methods as CSP, ILDM and TSR, provide iterative procedures to locally approximate the \emph{fast} and \emph{slow} subspaces resolving the \emph{tangent} space, along which the solution of the system evolves.~When the solution of the system evolves on the SIM, the projection of the vector field to the fast subspace is negligible and thus, it can be used to discover SIM approximations  \citep{kuehn2015multiple,kaper1999introduction,zagaris2004analysis}.~Hence, the aim of such methods is to approximate the basis vectors spanning the fast and slow subspaces via iterative procedures, that deliver increased approximation accuracy in every iteration.\par
CSP was originally developed for high-dimensional systems to provide numerical SIM approximations \citep{goussis1992study,lam1989understanding}.~However, when employed to low-dimensional systems, CSP may result in analytic SIM approximations, especially when a low number of iterations is performed.~Here, for a straightforward comparison with other analytic SIM approximations, we employ the CSP algorithmic procedure with one iteration \citep{valorani2005higher}, since the provided analytic SIM approximations are in an explicit form for the fast variable.~Additional CSP iterations may result in implicit forms, e.g., see the SIM approximation for the MM mechanism in Eq.~\eqref{eq:appMM_CSP2} when using two CSP iterations.~Note that for the employment of CSP, the original system in Eq.~\eqref{eq:gen} is not required to be written as a slow subsystem in the form of Eq.~\eqref{eq:SPslow}.~However, hereby we will adopt this slow subsystem formulation and briefly present the basic concepts of CSP and the resulting CSP-derived SIM approximations after one or more CSP iterations. 

Let the general slow subsystem in the form of Eq.~\eqref{eq:SPslow} be written in its matrix form:
\begin{equation}
    \dfrac{d}{dt} \begin{bmatrix} \mathbf{x} \\ \mathbf{y}
    \end{bmatrix} = \begin{bmatrix} \dfrac{1}{\epsilon} \mathbf{f}(\mathbf{x},\mathbf{y},\epsilon) \\ \mathbf{g}(\mathbf{x},\mathbf{y},\epsilon),
    \end{bmatrix} 
    \label{eq:CSPslow}
\end{equation}
for $(\mathbf{x},\mathbf{y}) \in \mathbb{R}^M \times \mathbb{R}^{N-M}$, $0 \le \epsilon \ll 1$ and $\mathbf{f}$, $\mathbf{g}$ sufficiently smooth functions.~The CSP iterative procedure requires the employment of two types of refinements, namely the $\mathbf{b}^r$- and $\mathbf{a}_r$-refinements; the former ensuring $\mathcal{O}(\epsilon)$ accuracy of the SIM approximation and the latter ensuring stability of the resulting slow system.~For initializing the algorithmic procedure, we consider the initial set of basis vectors:
\begin{equation}
    \mathbf{a}_M = \begin{bmatrix} \mathbf{I}^M_M \\[2pt] \mathbf{0}^{N-M}_M \end{bmatrix}, \qquad \mathbf{a}_{N-M} = \begin{bmatrix} \mathbf{0}^M_{N-M} \\[2pt] \mathbf{I}^{N-M}_{N-M} \end{bmatrix}, \qquad \mathbf{b}^M = \begin{bmatrix} \mathbf{I}^M_M & \mathbf{0}^M_{N-M} \end{bmatrix}, \qquad \mathbf{b}^{N-M} = \begin{bmatrix} \mathbf{0}^{N-M}_M & \mathbf{I}^{N-M}_{N-M} \end{bmatrix}
    \label{eq:InitBV}
\end{equation}
where $\mathbf{I}^r_r\in \mathbb{R}^{r\times r}$ and $\mathbf{0}^r_q \in \mathbb{R}^{r\times q}$ are the unitary and zero matrices, respectively.~The matrix $\mathbf{a}_M\in\mathbb{R}^{N \times M}$ in Eq.~\eqref{eq:InitBV} collects the $M$ in number $N$-dim. column basis vectors, which intend to approximate the $M$-dim. fast subspace, while the matrix $\mathbf{a}_{N-M}\in\mathbb{R}^{N \times (N-M)}$ collects the $(N-M)$ in number $N$-dim. column basis vectors, which intend to approximate the $(N-M)$-dim. slow subspace.~The matrices $\mathbf{b}^M\in\mathbb{R}^{M \times N}$ and $\mathbf{b}^{N-M}\in\mathbb{R}^{(N-M) \times N}$ in Eq.~\eqref{eq:InitBV} correspond to the matrices collecting the dual row basis vectors of the fast and slow subspaces respectively, since they satisfy the orthogonality conditions:
\begin{align}
    & \mathbf{b}^M \mathbf{a}_M = \mathbf{I}^M_M, \qquad \mathbf{b}^M \mathbf{a}_{N-M} = \mathbf{0}^M_{N-M}, \qquad \mathbf{b}^{N-M} \mathbf{a}_M = \mathbf{0}^{N-M}_M, \qquad \mathbf{b}^{N-M} \mathbf{a}_{N-M} = \mathbf{I}^{N-M}_{N-M}, \nonumber \\
    & \mathbf{a}_M \mathbf{b}^M + \mathbf{a}_{N-M} \mathbf{b}^{N-M} = \mathbf{I}^N_N
    \label{eq:Ortho}
\end{align}
Note here that the initial set of basis vectors in Eq.~\eqref{eq:InitBV} is in accordance to the system in Eq.~\eqref{eq:CSPslow}, since in the statevector $[\mathbf{x}, \mathbf{y}]^\top$ the first $M$ variables are considered fast, while the latter $N-M$ are considered slow.~Now, given the initial set of basis vectors in Eq.~\eqref{eq:InitBV}, the employment of the CSP $\mathbf{b}^r$-refinement results to a new basis, as:
\begin{equation}
    \mathbf{\tilde{a}}_M=\mathbf{a}_M, \quad \mathbf{\tilde{a}}_{N-M}=\left[\mathbf{I}^N_N-\mathbf{\tilde{a}}_M \mathbf{\tilde{b}}^M\right] \mathbf{a}_{N-M}, \quad
     \mathbf{\tilde{b}}^M = \boldsymbol{\tau}^M_M \mathbf{b}^M \mathbf{J}, \quad \mathbf{\tilde{b}}^{N-M}=\mathbf{b}^{N-M}
     \label{eq:Br1}
\end{equation}
where 
\begin{equation*}
    \boldsymbol{\tau}^M_M = \left( \mathbf{b}^M \mathbf{J} \mathbf{a}_M \right)^{-1} \in \mathbb{R}^{M \times M}
\end{equation*}
and $\mathbf{J}\in\mathbb{R}^{N \times N}$ is the Jacobian of the system in Eq.~\eqref{eq:CSPslow}, which can be expressed as:
\begin{equation}
    \mathbf{J} = \begin{bmatrix} \mathbf{J}^M_M & \mathbf{J}^M_{N-M} \\[2pt] \mathbf{J}^{N-M}_M & \mathbf{J}^{N-M}_{N-M} \end{bmatrix} = \begin{bmatrix} \dfrac{1}{\epsilon} \dfrac{\partial \mathbf{f}}{\partial \mathbf{x}} & \dfrac{1}{\epsilon} \dfrac{\partial \mathbf{f}}{\partial \mathbf{y}} \\[2pt]
    \dfrac{\partial \mathbf{g}}{\partial \mathbf{x}} & \dfrac{\partial \mathbf{g}}{\partial \mathbf{y}}
    \end{bmatrix}
    \label{eq:Jac}
\end{equation}
Thus, the resulting set of basis vectors Eq.~\eqref{eq:Br1} is: 
\begin{equation}
    \mathbf{\tilde{a}}_M=\begin{bmatrix} \mathbf{I}^M_M \\[2pt] \mathbf{0}^{N-M}_M \end{bmatrix}, \quad
    \mathbf{\tilde{a}}_{N-M}=\begin{bmatrix} -\mathbf{G}^M_{N-M} \\[2pt] \mathbf{I}^{N-M}_{N-M} \end{bmatrix}, \quad
     \mathbf{\tilde{b}}^M = \mathbf{Z}^M_M \begin{bmatrix} \mathbf{I}^M_M & \mathbf{G}^M_{N-M} \end{bmatrix}, \quad
     \mathbf{\tilde{b}}^{N-M}=\begin{bmatrix} \mathbf{0}^{N-M}_M & \mathbf{I}^{N-M}_{N-M} \end{bmatrix}
     \label{eq:Br2}
\end{equation}
where the included matrices are:
\begin{equation}
    \mathbf{Z}^M_M=\left( \mathbf{I}^M_M + \mathbf{G}^M_{N-M} \mathbf{G}^{N-M}_M \right)^{-1}, \quad \mathbf{G}^M_{N-M} = \left( \mathbf{J}^M_M\right)^{-1} \mathbf{J}^M_{N-M}, \quad \mathbf{G}^{N-M}_M = \mathbf{J}^{N-M}_M \left( \mathbf{J}^M_M\right)^{-1} 
    \label{eq:CSP_mat1}
\end{equation}
which are guaranteed to exist since $ \mathbf{J}^M_M$ is invertible, due to the invertibility of $\partial \mathbf{f}/\partial \mathbf{x}$ guaranteed by the normal hyperbolicity assumption of the SIM.~Now, given the set of basis vectors in Eq.~\eqref{eq:Br2}, the employment of the CSP $\mathbf{a}_r$-refinement results to a new basis, as:
\begin{equation}
    \mathbf{\hat{a}}_M= \mathbf{J} \mathbf{\tilde{a}}^M  \boldsymbol{\tilde{\tau}}^M_M , \quad
    \mathbf{\hat{a}}_{N-M}=\mathbf{\tilde{a}}_{N-M}, \quad
     \mathbf{\hat{b}}^M = \mathbf{\tilde{b}}^M,
     \quad \mathbf{\hat{b}}^{N-M}=\mathbf{\tilde{b}}^{N-M} \left[\mathbf{I}^N_N-\mathbf{\hat{a}}_M \mathbf{\hat{b}}^M\right]
     \label{eq:Ar1}
\end{equation}
where 
\begin{equation*}
    \boldsymbol{\tilde{\tau}}^M_M = \left( \mathbf{\tilde{b}}_M  \mathbf{J} \mathbf{\tilde{a}}^M \right)^{-1} \in \mathbb{R}^{M \times M}
\end{equation*}
Thus, the resulting set of basis vectors Eq.~\eqref{eq:Ar1} is: 
\begin{equation}
    \mathbf{\hat{a}}_M=\begin{bmatrix} \mathbf{I}^M_M \\[2pt] \mathbf{G}^{N-M}_M \end{bmatrix}, \quad
    \mathbf{\hat{a}}_{N-M}=\begin{bmatrix} -\mathbf{G}^M_{N-M} \\[2pt] \mathbf{I}^{N-M}_{N-M} \end{bmatrix}, \quad
     \mathbf{\hat{b}}^M = \mathbf{Z}^M_M \begin{bmatrix} \mathbf{I}^M_M & \mathbf{G}^M_{N-M} \end{bmatrix}, \quad
     \mathbf{\hat{b}}^{N-M}=\mathbf{Z}^{N-M}_{N-M}\begin{bmatrix} -\mathbf{G}^{N-M}_M & \mathbf{I}^{N-M}_{N-M} \end{bmatrix}
     \label{eq:Ar2}
\end{equation}
where
\begin{equation*}
    \mathbf{Z}^{N-M}_{N-M}=\left( \mathbf{I}^{N-M}_{N-M} + \mathbf{G}^{N-M}_M \mathbf{G}^M_{N-M} \right)^{-1}
\end{equation*}
Note that both sets of basis vectors in Eqs.~(\ref{eq:Br2}, \ref{eq:Ar2}) satisfy the orthogonality conditions in Eq.~\eqref{eq:Ortho}.~Given the set of basis vectors in Eq.~\eqref{eq:Ar2} the CSP reduced model after one iteration is:
\begin{equation}
    \mathbf{\hat{b}}^M \begin{bmatrix} \dfrac{1}{\epsilon} \mathbf{f}(\mathbf{x},\mathbf{y},\epsilon) \\ \mathbf{g}(\mathbf{x},\mathbf{y},\epsilon)
    \end{bmatrix} = \mathcal{O}(\epsilon), \qquad \dfrac{d}{dt} \begin{bmatrix} \mathbf{x} \\ \mathbf{y}
    \end{bmatrix} = \mathbf{\hat{a}}_{N-M} \mathbf{\hat{b}}^{N-M} \begin{bmatrix} \dfrac{1}{\epsilon} \mathbf{f}(\mathbf{x},\mathbf{y},\epsilon) \\ \mathbf{g}(\mathbf{x},\mathbf{y},\epsilon)
    \end{bmatrix} + \mathcal{O}(\epsilon)
    \label{eq:CSPRM}
\end{equation}
where the first $M$ algebraic expressions correspond to the CSP approximation of the SIM and the second system of $N$ ODEs corresponds to the slow system provided by CSP.~The $M$ equations of the SIM approximations and the latter $N-M$ differential equation for the slow variables $\mathbf{y}$ constitute the CSP reduced model after one iteration.~Being interested in the CSP approximation of the SIM, we will next employ the resulting from Eq.~\eqref{eq:CSPRM} expression:
\begin{equation}
    \mathbf{f}(\mathbf{x},\mathbf{y},\epsilon) + \epsilon \mathbf{G}^M_{N-M} \mathbf{g}(\mathbf{x},\mathbf{y},\epsilon) = \mathcal{O}(\epsilon^2)
    \label{eq:CSP_SIMgen}
\end{equation}
for deriving the expressions in Eqs.~(\ref{eq:MM_CSP}, \ref{eq:TMDD_CSP}, \ref{eq:TLC_CSP}) in order to prove the second part of Lemmas 1, 2 and 3.

For the employment of additional CSP iterations, further $\mathbf{b}^r$ and $\mathbf{a}_r$-refinements are employed to the set of basis vectors in Eq.~\eqref{eq:Ar2}.~In order to account for the curvature of the slow and fast subspaces, the additional refinements require the computation of the time derivatives of the basis vectors $d\mathbf{b}^r/dt$ and $d\mathbf{a}_r/dt$.~The latter typically involve the calculation of increasing, in order, derivatives of the Jacobian matrix.~A detailed presentation of the algorithmic procedure is shown in \citep{valorani2005higher}.~The resulting basis vectors can then be substituted to Eq.~\eqref{eq:CSPRM} for deriving the corresponding SIM approximation and slow system.~Here, we only perform CSP with two iterations for the MM system, for demonstrating that the resulting SIM approximation is in an implicit form w.r.t the fast variables.   

\subsection{MM system: proof of Lemma 1}
\label{app:MM_SIM_CSPexp}
Consider the MM subsystem in Eq.~\eqref{eq:MMss}, written in the matrix form of Eq.~\eqref{eq:CSPslow}, as:
\begin{equation}
    \dfrac{d}{dt} \begin{bmatrix} x \\ y
    \end{bmatrix} = \begin{bmatrix} \dfrac{1}{\epsilon} f(x,y,\epsilon) \\ g(x,y,\epsilon)
    \end{bmatrix} = \begin{bmatrix} \dfrac{1}{\epsilon} \left(  y - \dfrac{\kappa + y}{\kappa + 1} x \right) \\ \sigma^{-1} (-(\kappa +1) y + (\kappa-\sigma + y)x)
    \end{bmatrix} 
\end{equation}
Then, the Jacobian matrix and the scalar, in this case, ${G}^M_{N-M}$ in Eq.~\eqref{eq:CSP_mat1} take the form:
\begin{equation}
    \mathbf{J} = \begin{bmatrix} -(\epsilon (\kappa+1))^{-1}(\kappa+y) & (\epsilon (\kappa+1))^{-1} (\kappa+1-x) \\ \sigma^{-1} (\kappa-\sigma+y) & -\sigma^{-1} (\kappa+1-x) 
    \end{bmatrix}, \qquad
    G^M_{N-M} = -\dfrac{\kappa+1 - x}{\kappa + y}
\end{equation}
according to which, the SIM approximation provided by CSP in Eq.~\eqref{eq:CSP_SIMgen} is:
\begin{equation}
    y - \dfrac{x (\kappa + y)}{1 + \kappa} + 
 \epsilon \dfrac{(\kappa +1 - x) (\sigma x + y - x y + \kappa ( y-x))}{\sigma (\kappa + y)} = \mathcal{O}(\epsilon^2)
 \label{eq:MM_CSPSIMimpl}
\end{equation}
clearly written in an implicit form.~Solving for the fast variable $x$:
\begin{equation}
    x=h(y,\epsilon)=\dfrac{\sigma(\kappa+y)^2+\epsilon(\kappa+1)^2(\kappa-\sigma+2y)-\sqrt{\left(\epsilon(\kappa+1)^2(\kappa-\sigma)+\sigma(\kappa+y)^2 \right)^2+4\epsilon(\kappa+1)^2\sigma^2y(\kappa+y)}}{2 \epsilon (\kappa+1)(\kappa-\sigma+y)}
    \label{eq:appMM_CSP}
\end{equation}
recovers the explicit SIM approximation provided by CSP with one iteration for the MM subsystem, as presented in Eq.~\eqref{eq:MM_CSP}.

Note that a regular asymptotic expansion of the CSP-generated SIM approximation can be obtained with a Taylor expansion of Eq.~\eqref{eq:appMM_CSP} up to the first order around $\epsilon=0$, yielding: 
\begin{equation}
    x = h_0(y)+\epsilon h_1(y) + \mathcal{O}(\epsilon^2) = \dfrac{\kappa+1}{\kappa+y} y + \epsilon \dfrac{\kappa (\kappa+1)^3 y}{(\kappa+y)^4} + \mathcal{O}(\epsilon^2)
\end{equation}
This expression retrieves the $\mathcal{O}(\epsilon)$ SIM approximation derived on the basis of the invariance equation in Eq.~\eqref{eq:appMM_GSPT}, as expected by \citep{zagaris2004analysis,goussis2006efficient}.


Following the CSP algorithmic procedure described in \citep{valorani2005higher}, we performed an additional CSP iteration.~The resulting SIM approximation provided by CSP with two iterations yields:
\begin{align}
h(x,y)=
    & \dfrac{-\sigma (\kappa + y) ((1 + \kappa) y - x (\kappa + y)) (\sigma (\kappa + y)^2 +  \epsilon (1 + \kappa) (1 + \kappa - x) (\kappa - \sigma + y)) +   \epsilon (1 + \kappa) (-(1 + \kappa) y +}{\epsilon (1 + \kappa)^2 \sigma^2 (\kappa + y)^2 \left(\dfrac{\kappa + y}{\epsilon + \epsilon \kappa} +\right.} \nonumber \\
    & \dfrac{x (\kappa - \sigma + y)) (\epsilon (1 + \kappa) (1 + \kappa - x) (\kappa + \kappa^2 - \sigma x) + \sigma (\kappa + y) (\kappa (1 + \kappa - 2 x) + 2 (1 + \kappa - x) y))}{\left.\dfrac{(\kappa - \sigma + y) (\epsilon (1 + \kappa) (1 + \kappa - x) (\kappa + \kappa^2 - \sigma x) +  \sigma (\kappa + y) (\kappa (1 + \kappa - 2 x) + 2 (1 + \kappa - x) y))}{\sigma (\kappa + y) (\sigma (\kappa + y)^2 + \epsilon (1 + \kappa) (1 + \kappa - x) (\kappa - \sigma + y))}\right)} = \mathcal{O}(\epsilon^3) \label{eq:appMM_CSP2}
\end{align}
which is clearly a very complicated expression in an implicit form.

\subsection{TMDD system: proof of Lemma 2}
\label{app:TMDD_SIM_CSPexp}
Consider the TMDD subsystem in Eq.~\eqref{eq:TMDD_SPA}, written in the matrix form of Eq.~\eqref{eq:CSPslow}, as:
\begin{equation}
    \dfrac{d}{dt} \begin{bmatrix} x \\ y \\ z
    \end{bmatrix} = \begin{bmatrix} \dfrac{1}{\epsilon} f(x,y,z,\epsilon) \\ g^{(1)}(x,y,z,\epsilon) \\ g^{(2)}(x,y,z,\epsilon)
    \end{bmatrix} = \begin{bmatrix} \dfrac{1}{\epsilon} \left(  -xy+k_1z+1-\epsilon k_2 x \right) \\ k_3(-xy+k_1z)-k_4y \\ k_2(xy-k_1z)-z
    \end{bmatrix} 
\end{equation}
Then, the Jacobian matrix and the matrix $\mathbf{G}^M_{N-M}$ in Eq.~\eqref{eq:CSP_mat1} take the expressions:
\begin{equation}
    \mathbf{J} = \begin{bmatrix} -k_2 -\epsilon^{-1} y & -\epsilon^{-1} x & \epsilon^{-1} k_1 \\ -k_3 y & -k_4-k_3 x & k_1 k_3 \\ k_2 y & k_2 x & -1-k_1 k_2 
    \end{bmatrix},\qquad
    \mathbf{G}^M_{N-M} = \begin{bmatrix} \dfrac{x}{\epsilon k_2 +y} & -\dfrac{k_1}{\epsilon k_2 +y}\end{bmatrix}
\end{equation}
according to which, the SIM approximation provided by CSP in Eq.~\eqref{eq:CSP_SIMgen} is:
\begin{equation}
    1-xy+k_1z - \epsilon \dfrac{x (\epsilon k_2^2 + y(k_2 + k_1 k_2 + k_4 + k_3 x) ) + k_1 z (1 + k_1 k_2 + k_3 x)}{\epsilon k_2 + y}  = \mathcal{O}(\epsilon^2)
 \label{eq:TMDD_CSPSIMimpl}
\end{equation}
clearly written in an implicit form.~Solving for the fast variable $x$:
\begin{equation}
x = h(y,z,\epsilon) = -\dfrac{(y+\epsilon k_2)^2+\epsilon(y(k_1k_2+k_4)-k_1k_3z)}{2\epsilon k_3 y} \left(1-\sqrt{1+\dfrac{4 \epsilon k_3 y (\epsilon k_2 + y + k_1z (\epsilon  (1 + k_2 + k_1 k_2) + y))}{((y+\epsilon k_2)^2+\epsilon(y(k_1k_2+k_4)-k_1k_3z))^2}} \right)
    \label{eq:appTMDD_CSP}
\end{equation}
recovers the explicit SIM approximation provided by CSP with one iteration for the TMDD subsystem, as presented in Eq.~\eqref{eq:TMDD_CSP}.

Note that the regular asymptotic expansion of the CSP-generated SIM approximation, using a Taylor expansion up to the first order around $\epsilon=0$, in Eq.~\eqref{eq:appTMDD_CSP} yields: 
\begin{equation}
    x = h_0(y,z)+\epsilon h_1(y,z) + \mathcal{O}(\epsilon^2) = \dfrac{1+k_1z}{y} - \epsilon \dfrac{k_3(1+k_1)z+y(k_2+k_1k_2+k_4+k_1z(-1+k_2+k_4))}{y^3} + \mathcal{O}(\epsilon^2)
\end{equation}
which recovers the same expression with the $\mathcal{O}(\epsilon)$ SIM approximation derived on the basis of the invariance equation in Eq.~\eqref{eq:appTMDD_GSPT}, as expected by \citep{zagaris2004analysis,goussis2006efficient}.


\subsection{Sel'kov 3D system: proof of Lemma 3}
\label{app:TLC_SIM_CSPexp}
Consider the Sel'kov 3D system in Eq.~\eqref{eq:ToyLC_SPA}, written in the matrix form of Eq.~\eqref{eq:CSPslow}, as:
\begin{equation}
    \dfrac{d}{dt} \begin{bmatrix} x \\ y \\ z
    \end{bmatrix} = \begin{bmatrix} \dfrac{1}{\epsilon} f(x,y,z,\epsilon) \\ g^{(1)}(x,y,z,\epsilon) \\ g^{(2)}(x,y,z,\epsilon)
    \end{bmatrix} = \begin{bmatrix} \dfrac{1}{\epsilon} \left( y^2z-kxy \right) \\ az+y^2z-y+\epsilon x \\ -az-y^2z+b
    \end{bmatrix} 
\end{equation}
Then, the Jacobian matrix and the matrix $\mathbf{G}^M_{N-M}$ in Eq.~\eqref{eq:CSP_mat1} take the expressions:
\begin{equation}
    \mathbf{J} = \begin{bmatrix} -\epsilon^{-1} k y  & \epsilon^{-1} (2yz-kx) & \epsilon^{-1} y^2 \\ \epsilon & -1+2yz & a+y^2 \\ 0 & -2yz & -a-y^2 
    \end{bmatrix}
    , \qquad
    \mathbf{G}^M_{N-M} = \begin{bmatrix} -\dfrac{2yz-kx}{ky} & -\dfrac{y}{k}\end{bmatrix}
\end{equation}
according to which, the SIM approximation provided by CSP in Eq.~\eqref{eq:CSP_SIMgen} is:
\begin{equation}
     y^2z-kxy + \epsilon \dfrac{y^2(z(a+y^2)-b)+(2yz-kx)(y-z(a+y^2)-\epsilon x)}{ky}  = \mathcal{O}(\epsilon^2)
 \label{eq:TLC_CSPSIMimpl}
\end{equation}
clearly written in an implicit form.~Solving for the fast variable $x$:  
\begin{equation}
     x = h(y,z,\epsilon) = \left( 
     \dfrac{yz}{k}+\dfrac{y-z(a+y^2)}{2\epsilon}+\dfrac{ky^2}{2\epsilon^2}\right) \left(1 - \sqrt{ 1 + \dfrac{4 \epsilon^2 k y (\epsilon b y +z (2 \epsilon z (a + y^2) -y (k y + \epsilon (2 + a + y^2))))}{(k y (k y+\epsilon) + \epsilon z(2 \epsilon y - k (a + y^2)))^2}} \right)
    \label{eq:appTLC_CSP}
\end{equation}
recovers the explicit SIM approximation provided by CSP with one iteration for the Sel'kov 3D subsystem, as presented in Eq.~\eqref{eq:TLC_CSP}. 

As in all previous case studies, casting the CSP-generated SIM approximation in Eq.~\eqref{eq:appTLC_CSP} as a regular asymptotic expansion (Taylor expansion up to the first order around $\epsilon=0$) yields:
\begin{equation}
     x = h_0(y,z)+\epsilon h_1(y,z) + \mathcal{O}(\epsilon^2) = \dfrac{y z}{k} + \epsilon
    \left( \dfrac{z-b}{k^2} + \dfrac{z(a+y^2)(y-z)}{k^2y} \right) + \mathcal{O}(\epsilon^2)
\end{equation}
which recovers the same expression with the $\mathcal{O}(\epsilon)$ SIM approximation derived on the basis of the invariance equation in Eq.~\eqref{eq:appTLC_GSPT}.
\end{document}